\documentclass[a4paper, reqno, 11pt]{amsart}
\usepackage{color}
\makeatletter

\@addtoreset{equation}{section}
\makeatother
\usepackage{setspace}
\usepackage{xcolor}
\usepackage{here}
\usepackage{graphicx}
\usepackage{float}
\usepackage[T1]{fontenc}
\everymath{\displaystyle}
\usepackage{amsmath}
\usepackage{amssymb}
\usepackage{latexsym}
\usepackage{amsthm}
\usepackage[hidelinks]{hyperref}
\usepackage{mathtools}
\usepackage{mathrsfs}
\newtheorem{Thm}{Theorem}[section]

\newtheorem{Lem}[Thm]{Lemma}
\newtheorem{Prop}[Thm]{Proposition}

\newtheorem{Cor}[Thm]{Corollary}

\theoremstyle{definition}

\newtheorem{Rem}[Thm]{Remark}
\newtheorem{Exa}[Thm]{Example}

\begin{document}

\begin{abstract}
We consider conjugate equations driven by two families of linear fractional transformations on the unit interval, which we call the source system and the target system.
The solution is a topological conjugacy between two piecewise linear fractional interval maps, and we consider the Hausdorff dimension of the measure whose distribution function is the solution.
If both systems are hyperbolic in a certain sense, then we express this dimension as the ratio of two integrals with respect to the stationary distribution of an explicit Markov chain, and we show that the solution is a linear fractional transformation if finitely many explicit equations for the coefficients hold, and that otherwise it is a singular function and this dimension is smaller than $1$.
If only the source system is hyperbolic, then this dimension is $0$.
We also give closed formulas in terms of the Rogers dilogarithm in some cases.
\end{abstract}

\title[Measures induced by conjugacies]{Dimension formulas for measures induced by conjugacies between piecewise linear fractional maps}
\author{Kazuki Okamura}
\date{\today}
\address{Department of Mathematics, Faculty of Science, Shizuoka University, 836 Ohya, Suruga-ku, Shizuoka 422-8529, Japan}
\email{okamura.kazuki@shizuoka.ac.jp}
\keywords{Conjugate equations, de Rham's functional equations, piecewise linear fractional maps, topological conjugacy, Hausdorff dimension, singular functions, stationary measures, Rogers dilogarithm}
\subjclass[2020]{Primary 37E05, 28A80, 39B12; Secondary 37C45, 26A30, 11K55, 60J05, 60G30}
\maketitle

\section{Introduction}\label{sec:intro}

For an integer $N \ge 2$, we consider the conjugate equation
\[ g_i \circ \varphi = \varphi \circ f_i, \ \ \ 0 \le i \le N-1, \]
where $f = \{f_i\}_{0 \le i \le N-1}$ and $g = \{g_i\}_{0 \le i \le N-1}$ are two families of continuous strictly increasing maps from $[0,1]$ into itself such that $f_0 (0) = 0$, $f_{i-1} (1) = f_i (0)$ for $1 \le i \le N-1$ and $f_{N-1} (1) = 1$, and the same holds for $g$.
We call $f$ and $g$ the source system and the target system, respectively.
This equation contains de Rham's functional equations \cite{dR57}, and equations of this type have been considered by many authors; see \cite{BKK18, H85, O19, SB15, SB17, Z01} and the references therein.
Under mild assumptions, it has a unique continuous strictly increasing solution $\varphi : [0,1] \to [0,1]$, and $\varphi$ is often a singular function, that is, its derivative vanishes almost everywhere with respect to the Lebesgue measure.
There are two classical examples for $N = 2$.
If $f_0 (x) = x/2$, $f_1 (x) = (x+1)/2$, $g_0 (z) = pz$ and $g_1 (z) = (1-p) z + p$ for a number $p \in (0,1)$ with $p \ne 1/2$, then $\varphi$ is Lebesgue's singular function.
If $f_0 (x) = x/(x+1)$, $f_1 (x) = 1/(2-x)$, $g_0 (z) = z/2$ and $g_1 (z) = (z+1)/2$, then $\varphi$ is Minkowski's question-mark function \cite{M05}, which is a singular function \cite{D38, S43}; see also \cite{TU95}.

The above equation has a dynamical meaning.
Let $T$ and $\tau$ be the maps on $[0,1)$ whose inverse branches are $f_i$ and $g_i$, respectively, that is, $T \coloneqq  f_i^{-1}$ on $[f_i (0), f_i (1))$ and $\tau \coloneqq  g_i^{-1}$ on $[g_i (0), g_i (1))$ for $0 \le i \le N-1$.
Each of these maps has $N$ increasing branches which are onto $[0,1)$, and the solution $\varphi$ satisfies
\[ \varphi \circ T = \tau \circ \varphi \ \textup{ on } [0,1), \]
that is, $\varphi$ is a topological conjugacy between the interval maps $T$ and $\tau$.
For Lebesgue's singular function, $T$ is the doubling map and $\tau$ is a piecewise linear map.
For Minkowski's question-mark function, $T$ is a variant of the Farey map with two increasing branches, which is called a twisted Farey map in \cite{KS07}, and $\tau$ is the doubling map; see also \cite{P08}.
Let $\mu_{\varphi}$ be the Borel probability measure on $[0,1]$ whose distribution function is $\varphi$, and let $\ell$ be the Lebesgue measure.
Then, $\mu_{\varphi} (B) = \ell (\varphi (B))$ for every Borel subset $B$ of $[0,1]$, that is, $\mu_{\varphi}$ is the image measure of $\ell$ under $\varphi^{-1}$.
Hence, $\varphi$ is a singular function if and only if it maps a set of Lebesgue measure zero onto a set of full Lebesgue measure, and the Hausdorff dimension $\dim_H \mu_{\varphi}$ of $\mu_{\varphi}$ is the infimum of the Hausdorff dimensions of the Borel sets which $\varphi$ maps onto sets of full Lebesgue measure.

Thus the singularity of $\varphi$ is a question on the regularity of the conjugacy between two interval maps, which is a classical subject in dynamical systems.
Shub and Sullivan \cite{SS85} showed that two orientation preserving expanding $C^r$ maps of the circle, $r \ge 2$, which are conjugate by an absolutely continuous homeomorphism are $C^r$ conjugate.
For expanding piecewise $C^{1+\epsilon}$ maps of $[0,1]$ with finitely many increasing full branches, Jordan, Kesseb\"ohmer, Pollicott and Stratmann \cite{JKPS09} showed that the conjugacy is a $C^{1+\epsilon}$ diffeomorphism if the two maps are cohomologically dependent in the sense of \cite{JKPS09}, and that otherwise it is a singular function, for which they gave the Hausdorff dimension of the set of the points at which its derivative is $+\infty$; see also \cite{AJ26, KS08}.
Hofbauer and Raith \cite{HR92} showed that the Hausdorff dimension of an ergodic invariant measure with positive Lyapunov exponent for a piecewise monotonic interval map whose derivative is of bounded variation is the entropy divided by the Lyapunov exponent.
In general, it is not easy to verify the cohomological dependence or to compute these quantities for given maps.

One of the motivations of this paper is the explicit formula for the Hausdorff dimensions of self-similar measures.
If $f_i$ and $g_i$ are affine maps with slopes $r_i$ and $p_i$, respectively, for $0 \le i \le N-1$, then $\mu_{\varphi}$ is a self-similar measure, and
\[ \dim_H \mu_{\varphi} = \frac{\sum_{i=0}^{N-1} p_i \log p_i}{\sum_{i=0}^{N-1} p_i \log r_i}, \]
which is the entropy of $T$ with respect to $\mu_{\varphi}$ divided by its Lyapunov exponent.
We aim to obtain $\dim_H \mu_{\varphi}$ itself explicitly for maps which are not piecewise linear.

From these points of view, we investigate the case that both of $f$ and $g$ consist of linear fractional transformations.
We call such families LF systems; the definition is given in the next subsection.
Then, $T$ and $\tau$ are piecewise linear fractional maps.
This class of interval maps contains the maps $T$ and $\tau$ for the two classical examples, and it has been studied in metric number theory as fibred systems, for which explicit invariant densities and dual maps are known in several cases \cite{GH99, H02, S83, S95, S21, S24}.
We say that an LF system $h = \{h_i\}_{0 \le i \le N-1}$ is hyperbolic if $h_0^{\prime} (0) < 1$ and $h_{N-1}^{\prime} (1) < 1$, and that $h$ is parabolic otherwise, where $0$ and $1$ are the fixed points of $h_0$ and of $h_{N-1}$, respectively, and $h_i^{\prime} \le 1$ on $[0,1]$ for every $i$.
If $h$ is parabolic, then one of these fixed points is neutral, as for the above variant of the Farey map; see \cite{KMS16} for the Farey map.
If $h$ is hyperbolic, then the second iterate of the map whose inverse branches are $h_i$ is uniformly expanding, although the map itself is not in general (Remark \ref{rem:JKPS} (1)).

We now describe our main results for LF systems $f$ and $g$.
For $y > -1$, let $k(y,x) \coloneqq  (1+y)/(1+yx)^2$, $x \in [0,1]$, which is a probability density function on $[0,1]$.
For an LF system $h$, let $H_i^h$ be the linear fractional transformation given by the transposed matrix of a matrix which gives $h_i$, and let $G_i^h (y)$ be the integral of $k(y, \cdot)$ over the interval $[h_i (0), h_i (1)]$.
We call $H_i^h$ the dual maps of $h$.

The first result is a formula for $\dim_H \mu_{\varphi}$ (Theorem \ref{thm:dim}).
Assume that $f$ and $g$ are hyperbolic.
Then, there exist compact intervals $Y^g$ and $Y^f$ in $(-1, +\infty)$ and a unique Borel probability measure $\nu$ on $Y^g \times Y^f$ which is a stationary distribution of the Markov chain by which a point $(y,w)$ moves to $(H_i^g (y), H_i^f (w))$ with probability $G_i^g (y)$, $0 \le i \le N-1$, and
\[ \dim_H \mu_{\varphi} = \frac{\int_{Y^g \times Y^f} S(y) \, \nu(dy \, dw)}{\int_{Y^g \times Y^f} \Lambda (y,w) \, \nu(dy \, dw)}, \]
where $S(y)$ is the entropy of the probability vector $(G_i^g (y))_{0 \le i \le N-1}$, and $\Lambda (y,w)$ is the sum of $- G_i^g (y) \log G_i^f (w)$ over $0 \le i \le N-1$.
The numerator and the denominator are the entropy and the Lyapunov exponent of $T$ with respect to an ergodic invariant probability measure which is mutually absolutely continuous with $\mu_{\varphi}$, as in the formula in \cite{HR92}, which is not used in our proof (Subsection \ref{subsec:invariant} and Remark \ref{rem:related} (1)).
Both of them are integrals of explicit functions with respect to the single measure $\nu$, and the Lyapunov exponent is computed without differentiating $\varphi$.
If $f$ and $g$ consist of affine maps, then $\nu$ is the Dirac measure at the origin, and this is the above formula for self-similar measures.

The second result is a dichotomy for $\varphi$ (Theorem \ref{thm:sing}).
Assume that $f$ and $g$ are hyperbolic.
Then, $\varphi$ is a linear fractional transformation if $2N-2$ explicit equations for the coefficients of $f$ and $g$ hold, and otherwise $\varphi$ is a singular function and $\dim_H \mu_{\varphi} < 1$.
If $N = 2$, then these equations hold if and only if $f_0^{\prime} (0) = g_0^{\prime} (0)$ and $f_1^{\prime} (1) = g_1^{\prime} (1)$ (Remark \ref{rem:odds}).
This result is derived from the above formula: since $\Lambda (y,w) - S(y)$ is the relative entropy of $(G_i^g (y))_{0 \le i \le N-1}$ with respect to $(G_i^f (w))_{0 \le i \le N-1}$, we have that $\dim_H \mu_{\varphi} \le 1$, and the case that the equality holds is determined by algebraic computations for linear fractional transformations.
The approach of \cite{JKPS09}, which is based on the thermodynamic formalism and is different from ours, also gives a dichotomy between diffeomorphisms and singular functions, since \cite[Theorem 1.2]{JKPS09} applies to $T^2$ and $\tau^2$; the above equations decide the cohomological dependence of $T^2$ and $\tau^2$ in the sense of \cite{JKPS09} (Remark \ref{rem:JKPS}).

The third result concerns the case that the target system is parabolic (Theorem \ref{thm:para} and Corollary \ref{cor:positive}), which is outside the setting of \cite{JKPS09}, since no iterate of $\tau$ is uniformly expanding.
In this case, $\mu_{\varphi}$ does not have the Gibbs property with respect to $T$ for any function, whereas it has this property if $g$ is hyperbolic (Remark \ref{rem:Gibbs-property}); measures without the Gibbs property are harder to deal with.
If $f$ is hyperbolic and $g$ is parabolic, then $\dim_H \mu_{\varphi} = 0$.
Thus, for a hyperbolic LF system $f$, $\dim_H \mu_{\varphi} > 0$ if and only if $g$ is hyperbolic.
For example, the measure whose distribution function is the inverse function of Minkowski's question-mark function has Hausdorff dimension $0$ (Example \ref{exa:Minkowski}), which settles a question left open in \cite[Example 35]{O20}.
In contrast, the measure whose distribution function is Minkowski's question-mark function, for which the source system is parabolic and the target system is hyperbolic, has positive Hausdorff dimension by a formula of Kinney \cite{K60}; we do not consider this case in general (Remark \ref{rem:assumption} (3)).

The fourth result gives closed formulas (Theorem \ref{thm:tiling}).
Assume that $f$ and $g$ are hyperbolic, and that the images of a compact interval of positive length in $(-1, +\infty)$ under the dual maps of $g$ have pairwise disjoint interiors and their union is this interval.
Then, the numerator of the above formula, which is the entropy of $\tau$ with respect to its absolutely continuous invariant probability measure, is expressed by the Rogers dilogarithm $L$, and the denominator also has a closed form if the dual maps of $f$ have a common fixed point, which is the case if $f$ consists of affine maps (Corollary \ref{cor:common}).
If $N = 2$, then, for every hyperbolic LF system $g$, this entropy is
\[ \frac{2 \left( L(\kappa_0) + L(\kappa_1) - \pi^2 / 6 \right)}{\log (\kappa_0 / (1 - \kappa_0)) + \log (\kappa_1 / (1 - \kappa_1))}, \ \ \ \kappa_0 \coloneqq  g_0^{\prime} (0), \ \ \kappa_1 \coloneqq  g_1^{\prime} (1), \]
if $\kappa_0 + \kappa_1 \ne 1$, and it is $- \kappa_0 \log \kappa_0 - \kappa_1 \log \kappa_1$ if $\kappa_0 + \kappa_1 = 1$ (Remark \ref{rem:two}).
This is a counterpart of the value $\pi^2 / (6 \log 2)$ of the entropy of the continued fraction map with respect to the Gauss measure \cite{IK02}; see also \cite{H02}.
In particular, for the equation considered in \cite{A26, O14}, that is, for $N = 2$, $f_0 (x) = x/2$ and $f_1 (x) = (x+1)/2$, $\dim_H \mu_{\varphi}$ is this entropy divided by $\log 2$.
This solves a problem posed by Arcozzi \cite[Subsection 3.3.1 (i)]{A26}, which is equivalent to determining $\dim_H \mu_{\varphi}$ in the case that $\varphi$ is the normalized cumulative capacitary function of a dyadic Riesz capacity (Remark \ref{rem:Arcozzi}).

The idea of the proofs is a duality for linear fractional transformations: for an LF system $h$,
\[ h_i^{\prime} (x) \, k(y, h_i (x)) = G_i^h (y) \, k \! \left( H_i^h (y), x \right), \ \ \ x \in [0,1], \ \ y > -1 \]
(Lemma \ref{lem:id} (3)).
For $h = g$, this shows that the Perron--Frobenius operator of $\tau$ maps the density $k(y, \cdot)$ to the mixture of the densities $k(H_i^g (y), \cdot)$ with the weights $G_i^g (y)$, $0 \le i \le N-1$.
Hence, if $\nu^g$ is a stationary distribution of the Markov chain by which $y$ moves to $H_i^g (y)$ with probability $G_i^g (y)$, $0 \le i \le N-1$, then the integral of $k(y, \cdot)$ with respect to $\nu^g (dy)$ is an invariant density of $\tau$ (Proposition \ref{prop:tau}).
This extends the method of dual algorithms of Schweiger \cite{S83, S21}, in which the dual maps are assumed to form a fibred system on a set; see also \cite{P22}.
In general, the images of an interval under the dual maps may overlap, and $\nu^g$ is a stationary distribution of the same kind as the Blackwell measures \cite{BPS12, B57}.
The measure $\nu$ is the stationary distribution for the pairs of the dual variables of $g$ and of $f$ which are driven by the same digits, and it gives the Lyapunov exponent of $T$ by the above identity for $h = g$ and for $h = f$ (Proposition \ref{prop:average}).

We now compare our results with the previous results for the conjugate equation for LF systems.
Upper and lower bounds for $\dim_H \mu_{\varphi}$ and a necessary and sufficient condition for the singularity were given in \cite{O14, O20} in the case that $f$ consists of affine maps, and sufficient conditions for the singularity were given in \cite{O26} for a wider class of source systems.
The functions $G_i^g$ and $H_i^g$ for the target system were already used in these papers, in which the value of $\dim_H \mu_{\varphi}$ was not obtained in general.
The new ingredients here are the dual maps of the source system, the stationary distribution $\nu$, and the invariant measures of $\tau$ and $T$.
The formula for $\dim_H \mu_{\varphi}$, the necessary and sufficient condition for the singularity in the case that $f$ does not consist of affine maps, the result for parabolic target systems and the closed formulas are new; in particular, the bounds for $\dim_H \mu_{\varphi}$ in \cite{O14} are replaced by its value in closed form.
See Remarks \ref{rem:Blackwell}, \ref{rem:Furstenberg}, \ref{rem:projective}, \ref{rem:related} and \ref{rem:Arcozzi} for related results, including recent results of Arcozzi \cite{A26} for the equation in \cite{O14}.

\subsection{Framework and main results}

Let $N \ge 2$ and $\mathcal{I}_N \coloneqq  \{0, 1, \dots, N-1\}$.
We follow the terminology and the notation of \cite{O26}.
We say that a family $h = \{h_i\}_{i \in \mathcal{I}_N}$ of functions on $[0,1]$ is an \textit{LF system} on $[0,1]$ if
\[ h_i (x) = \frac{a_i x + b_i}{c_i x + d_i}, \ \ \ x \in [0,1], \ \ i \in \mathcal{I}_N, \]
where $a_i, b_i, c_i, d_i$, $i \in \mathcal{I}_N$, are real numbers satisfying the following conditions (1)--(3).

\medskip

\noindent (1) $c_i + d_i > 0$ and $d_i > 0$ for $i \in \mathcal{I}_N$. \\[1mm]
(2) $0 < a_i d_i - b_i c_i \le \min\{d_i^2, (c_i + d_i)^2\}$ for $i \in \mathcal{I}_N$. \\[1mm]
(3) $b_0 = 0$ and $\frac{a_i + b_i}{c_i + d_i} = \frac{b_{i+1}}{d_{i+1}}$ for $i \in \mathcal{I}_N$, where we let $\frac{b_N}{d_N} \coloneqq  1$.

\medskip

This is \cite[Definition 2.8]{O26}, in which there is one more condition (4), which follows from (1)--(3).
Since $h_i$ does not change if $a_i$, $b_i$, $c_i$ and $d_i$ are multiplied by a positive constant, {\it we assume that $d_i = 1$ for every $i \in \mathcal{I}_N$ throughout this paper}, and we let $b_N \coloneqq  1$.
We denote the numbers $a_i$, $b_i$ and $c_i$ for an LF system $h$ by $a_i^h$, $b_i^h$ and $c_i^h$, respectively.
Then, $h_i (0) = b_i^h$ and $h_i (1) = b_{i+1}^h$ for $i \in \mathcal{I}_N$, and $0 = b_0^h < b_1^h < \cdots < b_N^h = 1$ (Lemma \ref{lem:basic}).

Let $f = \{f_i\}_{i \in \mathcal{I}_N}$ and $g = \{g_i\}_{i \in \mathcal{I}_N}$ be LF systems on $[0,1]$.
By \cite[Lemma 2.3, Lemma 2.7 and Lemma 2.9]{O26}, there exists a unique continuous strictly increasing function $\varphi : [0,1] \to [0,1]$ such that
\begin{equation}\label{eq:conj}
g_i \circ \varphi = \varphi \circ f_i, \ \ \ i \in \mathcal{I}_N.
\end{equation}
Since $f_0 (1) = f_1 (0)$, by \eqref{eq:conj}, $g_0 (\varphi (1)) = \varphi (f_0 (1)) = \varphi (f_1 (0)) = g_1 (\varphi (0))$.
Since $g_0$ and $g_1$ are strictly increasing on $[0,1]$ and $g_0 (1) = g_1 (0)$ (Lemma \ref{lem:basic}), we have that $g_0 (\varphi (1)) \le g_0 (1) = g_1 (0) \le g_1 (\varphi (0))$.
Hence, the equalities hold in these inequalities, and therefore $\varphi (0) = 0$ and $\varphi (1) = 1$.
Let $\mu_{\varphi}$ be the Borel probability measure on $[0,1]$ whose distribution function is $\varphi$, and let
\[ \dim_H \mu_{\varphi} \coloneqq  \inf\left\{ \dim_H A : A \in \mathcal{B}([0,1]), \ \mu_{\varphi} (A) = 1 \right\}, \]
where $\mathcal{B}([0,1])$ is the Borel $\sigma$-algebra on $[0,1]$ and $\dim_H A$ denotes the Hausdorff dimension of $A$.

For an LF system $h$, $y > -1$ and $i \in \mathcal{I}_N$, let
\begin{equation}\label{eq:GH}
G_i^h (y) \coloneqq  \frac{(b_{i+1}^h - b_i^h)(y + 1)}{(b_{i+1}^h y + 1)(b_i^h y + 1)}, \ \ \ H_i^h (y) \coloneqq  \frac{a_i^h y + c_i^h}{b_i^h y + 1},
\end{equation}
and let $\mathbf{G}^h (y) \coloneqq  \left( G_0^h (y), \dots, G_{N-1}^h (y) \right)$.
We call the functions $H_i^h$, $i \in \mathcal{I}_N$, the dual maps of $h$.
We see in Lemma \ref{lem:id} below that $b_i^h y + 1 > 0$, $G_i^h (y) > 0$ and $G_0^h (y) + \cdots + G_{N-1}^h (y) = 1$ for $y > -1$.

We say that an LF system $h$ is \textit{hyperbolic} if $h_0^{\prime} (0) < 1$ and $h_{N-1}^{\prime} (1) < 1$, and that $h$ is \textit{parabolic} if it is not hyperbolic.
The points $0$ and $1$ are fixed points of $h_0$ and $h_{N-1}$, respectively.
Since $h_i^{\prime} \le 1$ on $[0,1]$ for every $i \in \mathcal{I}_N$ (Lemma \ref{lem:basic}), $h$ is parabolic if and only if $h_0^{\prime} (0) = 1$ or $h_{N-1}^{\prime} (1) = 1$.
An LF system $h$ is hyperbolic if and only if $a_0^h < 1$ and $b_{N-1}^h + c_{N-1}^h > 0$ (Lemma \ref{lem:basic}).
For a hyperbolic LF system $h$, let
\[ e_0^h \coloneqq  \frac{c_0^h}{1 - a_0^h}, \ \ \ e_i^h \coloneqq  \frac{a_i^h - 1 + \left( (1 - a_i^h)^2 + 4 b_i^h c_i^h \right)^{1/2}}{2 b_i^h}, \ \ 1 \le i \le N-1, \]
\[ \alpha_h \coloneqq  \min\left\{ 0, e_0^h, \dots, e_{N-1}^h \right\}, \ \ \ \beta_h \coloneqq  \max\left\{ 0, e_0^h, \dots, e_{N-1}^h \right\}, \]
and $Y^h \coloneqq  [\alpha_h, \beta_h]$.
For each $i \in \mathcal{I}_N$, $e_i^h$ is the unique fixed point of $H_i^h$ in $(-1, +\infty)$ (Lemma \ref{lem:fixed}).
Furthermore, $-1 < \alpha_h \le 0 \le \beta_h < +\infty$ and $H_i^h (Y^h) \subset Y^h$ for every $i \in \mathcal{I}_N$ (Lemma \ref{lem:Y}).

Let $\ell$ be the Lebesgue measure on $[0,1]$.
For $n \ge 0$, we call an element $\sigma = (\sigma_1, \dots, \sigma_n)$ of $(\mathcal{I}_N)^n$ a word of length $n$, where $(\mathcal{I}_N)^0$ consists of the empty word.
For an LF system $h$ and a word $\sigma = (\sigma_1, \dots, \sigma_n)$, let $h_{\sigma} \coloneqq  h_{\sigma_1} \circ \cdots \circ h_{\sigma_n}$, which is the identity map if $n = 0$.
For a word $\sigma$, let
\[ J(\sigma) \coloneqq  \left[ f_{\sigma} (0), f_{\sigma} (1) \right). \]
For every $n \ge 0$, the intervals $J(\sigma)$, $\sigma \in (\mathcal{I}_N)^n$, are pairwise disjoint and their union is $[0,1)$ (Lemma \ref{lem:cyl}).
For $n \ge 1$ and $x \in [0,1)$, let $J_n (x)$ be the interval $J(\sigma)$, $\sigma \in (\mathcal{I}_N)^n$, which contains $x$.
By \eqref{eq:conj}, $\mu_{\varphi} (J(\sigma)) = g_{\sigma} (1) - g_{\sigma} (0)$, and $\ell (J(\sigma)) = f_{\sigma} (1) - f_{\sigma} (0)$ (Lemma \ref{lem:conj}).

For $y > -1$ and $w > -1$, let
\[ S(y) \coloneqq  - \sum_{i \in \mathcal{I}_N} G_i^g (y) \log G_i^g (y), \ \ \ \Lambda (y,w) \coloneqq  - \sum_{i \in \mathcal{I}_N} G_i^g (y) \log G_i^f (w). \]
For probability vectors $\mathbf{p} = (p_i)_{i \in \mathcal{I}_N}$ and $\mathbf{q} = (q_i)_{i \in \mathcal{I}_N}$ with positive entries, let
\[ s_N (\mathbf{p} | \mathbf{q}) \coloneqq  \sum_{i \in \mathcal{I}_N} p_i \log \frac{q_i}{p_i}. \]
By Jensen's inequality, $s_N (\mathbf{p} | \mathbf{q}) \le 0$, and $s_N (\mathbf{p} | \mathbf{q}) = 0$ if and only if $\mathbf{p} = \mathbf{q}$.
We have that
\[ S(y) - \Lambda (y,w) = s_N \left( \mathbf{G}^g (y) | \mathbf{G}^f (w) \right), \ \ \ y, w > -1. \]

\begin{Thm}\label{thm:dim}
Assume that $f$ and $g$ are hyperbolic LF systems on $[0,1]$, and let $\varphi$ be the solution of \eqref{eq:conj}. Then,\\
(1) There exists a Borel probability measure $\nu$ on $Y^g \times Y^f$ such that
\begin{equation}\label{eq:stationary}
\int_{Y^g \times Y^f} \sum_{i \in \mathcal{I}_N} G_i^g (y) \, u \! \left( H_i^g (y), H_i^f (w) \right) \nu(dy \, dw) = \int_{Y^g \times Y^f} u \, d\nu
\end{equation}
for every continuous function $u$ on $Y^g \times Y^f$.\\
(2) For every Borel probability measure $\nu$ on $Y^g \times Y^f$ satisfying \eqref{eq:stationary} for every continuous function $u$ on $Y^g \times Y^f$,
\begin{align*}
\lim_{n \to \infty} \frac{- \log \mu_{\varphi} (J_n (x))}{n} &= \int_{Y^g \times Y^f} S(y) \, \nu(dy \, dw), \ \ \ \textup{$\mu_{\varphi}$-a.s.} \, x, \\
\lim_{n \to \infty} \frac{- \log \ell (J_n (x))}{n} &= \int_{Y^g \times Y^f} \Lambda (y,w) \, \nu(dy \, dw), \ \ \ \textup{$\mu_{\varphi}$-a.s.} \, x,
\end{align*}
and
\begin{align*}
\dim_H \mu_{\varphi} &= \frac{\int_{Y^g \times Y^f} S(y) \, \nu(dy \, dw)}{\int_{Y^g \times Y^f} \Lambda (y,w) \, \nu(dy \, dw)} \\
&= 1 + \frac{\int_{Y^g \times Y^f} s_N \left( \mathbf{G}^g (y) | \mathbf{G}^f (w) \right) \nu(dy \, dw)}{\int_{Y^g \times Y^f} \Lambda (y,w) \, \nu(dy \, dw)}.
\end{align*}
\end{Thm}

The equation \eqref{eq:stationary} states that $\nu$ is a stationary distribution of the transition kernel on $Y^g \times Y^f$ by which $(y,w)$ moves to $\left( H_i^g (y), H_i^f (w) \right)$ with probability $G_i^g (y)$, $i \in \mathcal{I}_N$.
The first coordinate $y$ and the second coordinate $w$ are variables for $g$ and for $f$, respectively, and they are driven by the same index $i$, whose probability is determined by $y$.
Such a measure $\nu$ is unique, and it is a stationary distribution of the same kind as the Blackwell measures \cite{B57} (Remark \ref{rem:Blackwell}).
By Theorem \ref{thm:dim} (2),
\[ \lim_{n \to \infty} \frac{1}{n} \log \frac{\mu_{\varphi} (J_n (x))}{\ell (J_n (x))} = - \int_{Y^g \times Y^f} s_N \left( \mathbf{G}^g (y) | \mathbf{G}^f (w) \right) \nu(dy \, dw), \ \ \ \textup{$\mu_{\varphi}$-a.s.} \, x. \]
Since $S$ and $\Lambda$ are explicit functions and $\nu$ is a stationary distribution of an explicit transition kernel, Theorem \ref{thm:dim} (2) is also useful for numerical computations of $\dim_H \mu_{\varphi}$: the integrals with respect to $\nu$ are replaced by finite sums obtained by iterating this transition kernel, the averages of which over the numbers of the iterations converge to these integrals, and a density of $\nu$ is not needed (Example \ref{exa:two} (2)).
If $f$ consists of affine maps, then $Y^f = \{0\}$, the numbers $G_i^f (0)$, $i \in \mathcal{I}_N$, are the slopes of $f_i$, $i \in \mathcal{I}_N$, and \eqref{eq:stationary} is an equation for a measure on $Y^g$ (Remark \ref{rem:affine}).
In this case, the Borel probability measure on $[0,1]$ whose distribution function is $\varphi^{-1}$ is the unique stationary measure on $[0,1]$ for the maps $g_i$, $i \in \mathcal{I}_N$, chosen independently with the probabilities given by these slopes (Remark \ref{rem:Furstenberg}).
It is the unique stationary measure on the projective line, which is the Furstenberg measure in the sense of \cite{HS17}, if and only if the dual maps of $g$ have no common fixed points in $(-1, +\infty)$ (Remark \ref{rem:projective}).
We do not assume any relation between $Y^f$ and $Y^g$.

We now consider the singularity of $\mu_{\varphi}$.
For $s > -1$, let
\[ \psi_s (x) \coloneqq  \frac{(1+s) x}{1 + s x}, \ \ \ x \in [0,1], \ \ \ \ \gamma_s (y) \coloneqq  (1+s)(1+y) - 1, \ \ \ y > -1. \]
The function $\psi_s$ is a linear fractional transformation which is a homeomorphism of $[0,1]$ such that $\psi_s (0) = 0$ and $\psi_s (1) = 1$, and $\gamma_s$ is an affine bijection of $(-1, +\infty)$ (Lemma \ref{lem:psi}).
Since $f_0 (1) = b_1^f$ and $g_0 (1) = b_1^g$ are in $(0,1)$, there exists a unique number $s_0 > -1$ such that $\psi_{s_0} (f_0 (1)) = g_0 (1)$, that is,
\[ 1 + s_0 = \frac{g_0 (1) \left( 1 - f_0 (1) \right)}{\left( 1 - g_0 (1) \right) f_0 (1)}. \]

\begin{Thm}\label{thm:sing}
Assume that $f$ and $g$ are hyperbolic LF systems on $[0,1]$, and let $\varphi$ be the solution of \eqref{eq:conj}. Then,\\
(1) If $g_i \circ \psi_{s_0} = \psi_{s_0} \circ f_i$ on $[0,1]$ for every $i \in \mathcal{I}_N$, then $\varphi = \psi_{s_0}$, and $\mu_{\varphi}$ is absolutely continuous with respect to the Lebesgue measure.\\
(2) Otherwise, $\dim_H \mu_{\varphi} < 1$, and $\mu_{\varphi}$ is singular with respect to the Lebesgue measure.
\end{Thm}

Thus, for hyperbolic LF systems $f$ and $g$, the solution $\varphi$ of \eqref{eq:conj} is either a linear fractional transformation or a singular function, and $\mu_{\varphi}$ is singular if and only if $g_i \circ \psi_{s_0} \ne \psi_{s_0} \circ f_i$ for some $i \in \mathcal{I}_N$.
See the beginning of \cite[Section 3]{O26} for singular functions.
See Remark \ref{rem:JKPS} for the relation between Theorem \ref{thm:sing} and the results in \cite{JKPS09}.
Assertion (1) holds for all LF systems $f$ and $g$.
The assumption of Theorem \ref{thm:sing} (1) holds if and only if
\[ b_i^g = \psi_{s_0} (b_i^f), \ 1 \le i \le N-1, \ \ \ (1 + s_0)(1 + c_i^g) = 1 + H_i^f (s_0), \ i \in \mathcal{I}_N  \]
(Proposition \ref{prop:equiv}).
Since the first equation for $i = 1$ holds by the definition of $s_0$, this is a system of $2N - 2$ equations, which are the first one for $2 \le i \le N-1$ and the second one for $i \in \mathcal{I}_N$.
These are equations for the $4N - 2$ numbers $b_i^h$, $1 \le i \le N-1$, and $c_i^h$, $i \in \mathcal{I}_N$, $h \in \{f, g\}$, by which the LF systems $f$ and $g$ are determined (Lemma \ref{lem:basic}).

We now consider the case that $g$ is parabolic.

\begin{Thm}\label{thm:para}
Let $f$ and $g$ be LF systems on $[0,1]$, and let $\varphi$ be the solution of \eqref{eq:conj}.
Assume that $g$ is parabolic. Then,\\
(1) It holds that
\[ \lim_{n \to \infty} \frac{- \log \mu_{\varphi} (J_n (x))}{n} = 0, \ \ \ \textup{$\mu_{\varphi}$-a.s.} \, x. \]
(2) If $f$ is hyperbolic, then $\dim_H \mu_{\varphi} = 0$.
\end{Thm}

\begin{Cor}\label{cor:positive}
Assume that $f$ is a hyperbolic LF system on $[0,1]$ and $g$ is an LF system on $[0,1]$, and let $\varphi$ be the solution of \eqref{eq:conj}.
Then, $\dim_H \mu_{\varphi} > 0$ if and only if $g$ is hyperbolic.
\end{Cor}

Since $\mu_{\varphi} (J(\sigma)) = g_{\sigma} (1) - g_{\sigma} (0)$ for every word $\sigma$, Theorem \ref{thm:para} (1) is an assertion for $g$ alone, and $f$ is not assumed to be hyperbolic in it.
By Theorem \ref{thm:para} (2), if $f$ is hyperbolic and $g$ is parabolic, then $\mu_{\varphi}$ is singular with respect to the Lebesgue measure, and hence the assumption of Theorem \ref{thm:sing} (1) does not hold.
Thus the assertions of Theorem \ref{thm:sing} hold for every hyperbolic LF system $f$ and every LF system $g$.
Theorem \ref{thm:para} (2) does not hold in general if $f$ is parabolic (Remark \ref{rem:assumption} (3)).
If $N = 2$, $f_0 (x) = x/2$, $f_1 (x) = (x+1)/2$, $g_0 (z) = z/(z+1)$ and $g_1 (z) = 1/(2-z)$, then $\varphi$ is the inverse function of Minkowski's question-mark function, and $\dim_H \mu_{\varphi} = 0$ (Example \ref{exa:Minkowski}).

We now consider expressions for the two limits in Theorem \ref{thm:dim} (2) which do not involve the measure $\nu$.
Let $T : [0,1) \to [0,1)$ and $X_1 : [0,1) \to \mathcal{I}_N$ be the maps defined by
\[ T(x) \coloneqq  f_i^{-1} (x), \ \ \ X_1 (x) \coloneqq  i, \ \ \ b_i^f \le x < b_{i+1}^f, \ \ i \in \mathcal{I}_N. \]
For $j \ge 0$, we denote the $j$-fold composition of $T$ by $T^j$, and, for $n \ge 1$, let $X_n \coloneqq  X_1 \circ T^{n-1}$.
Then, $J_n (x) = J((X_1 (x), \dots, X_n (x)))$ for $n \ge 1$ and $x \in [0,1)$ (Lemma \ref{lem:coding}).
Let
\[ \omega_f (x) \coloneqq  \log f_{X_1 (x)}^{\prime} (T(x)), \ \ \ x \in [0,1). \]
The map $T$ is differentiable on each open interval $(b_i^f, b_{i+1}^f)$, $i \in \mathcal{I}_N$, and $- \omega_f$ is the logarithm of its derivative there.
Let $L$ be the Rogers dilogarithm, that is,
\[ L(u) \coloneqq  \sum_{n=1}^{\infty} \frac{u^n}{n^2} + \frac{1}{2} \log u \log (1-u), \ \ \ 0 < u < 1. \]
See \cite{K95, Z07} for the properties of this function.
For $y > -1$ and $1 \le i \le N-1$, let
\[ \chi_i (y) \coloneqq  \frac{b_i^g \left( 1 + b_{i+1}^g y \right)}{b_{i+1}^g \left( 1 + b_i^g y \right)}. \]
Then, $0 < \chi_i (y) < 1$ (Lemma \ref{lem:chi}).

\begin{Thm}\label{thm:tiling}
Assume that $f$ and $g$ are hyperbolic LF systems on $[0,1]$, and let $\varphi$ be the solution of \eqref{eq:conj}.
Assume that there exist real numbers $e_{-}$ and $e_{+}$ with $-1 < e_{-} < e_{+}$ such that the closed intervals $H_i^g ([e_{-}, e_{+}])$, $i \in \mathcal{I}_N$, have pairwise disjoint interiors and their union is $[e_{-}, e_{+}]$. Then,\\
(1) It holds that
\[ \lim_{n \to \infty} \frac{- \log \mu_{\varphi} (J_n (x))}{n} = \frac{2 \sum_{i=1}^{N-1} \left( L(\chi_i (e_{+})) - L(\chi_i (e_{-})) \right)}{\log (1 + e_{+}) - \log (1 + e_{-})}, \ \ \ \textup{$\mu_{\varphi}$-a.s.} \, x. \]
(2) For $\mu_{\varphi}$-a.s.\ $x$,
\begin{align*}
&\lim_{n \to \infty} \frac{- \log \ell (J_n (x))}{n} \\
&= \frac{e_{+} - e_{-}}{\log (1 + e_{+}) - \log (1 + e_{-})} \int_{[0,1)} \frac{- \omega_f (x^{\prime})}{(1 + e_{-} \varphi (x^{\prime}))(1 + e_{+} \varphi (x^{\prime}))} \, \mu_{\varphi} (dx^{\prime}).
\end{align*}
(3) It holds that
\[ \dim_H \mu_{\varphi} = \frac{\displaystyle 2 \sum_{i=1}^{N-1} \left( L(\chi_i (e_{+})) - L(\chi_i (e_{-})) \right)}{\displaystyle (e_{+} - e_{-}) \int_{[0,1)} \frac{- \omega_f (x)}{(1 + e_{-} \varphi (x))(1 + e_{+} \varphi (x))} \, \mu_{\varphi} (dx)}. \]
\end{Thm}

The limit in (1) is determined by $b_i^g$, $1 \le i \le N-1$, $e_{-}$ and $e_{+}$.
The limit in (2) is the integral with respect to $\mu_{\varphi}$ of a function given by $\omega_f$ and $\varphi$, and we do not know whether it has a closed form in general.
It is the Lyapunov exponent of $T$ with respect to a $T$-invariant measure whose density with respect to $\mu_{\varphi}$ is given by $e_{-}$, $e_{+}$ and $\varphi$ (Remark \ref{rem:common} (4)).
If the dual maps of $f$ have a common fixed point, which is the case if $f$ consists of affine maps, then the limit in (2) has a closed form (Corollary \ref{cor:common} and Remark \ref{rem:common}).

The assumption of Theorem \ref{thm:tiling} is a condition for $g$ alone (Remark \ref{rem:tiling}).
It holds if $e_0^g \ne e_{N-1}^g$ and $H_{i-1}^g (e_{N-1}^g) = H_i^g (e_0^g)$ for every $1 \le i \le N-1$, where the closed interval with endpoints $e_0^g$ and $e_{N-1}^g$ is taken as $[e_{-}, e_{+}]$ (Corollary \ref{cor:natural}).
These identities hold if and only if $\Theta_g \circ g_i = H_i^g \circ \Theta_g$ for every $i \in \mathcal{I}_N$, where $\Theta_g$ is the function on $[0,1]$ defined in Subsection \ref{subsec:natural}, which is a linear fractional transformation if $e_0^g \ne e_{N-1}^g$ (Proposition \ref{prop:natural}).
If $N = 2$, then the identity $H_0^g (e_1^g) = H_1^g (e_0^g)$ always holds, and the limit in (1) is expressed by $g_0^{\prime} (0)$ and $g_1^{\prime} (1)$.
If $N = 2$ and $f_0^{\prime} (0) + f_1^{\prime} (1) = 1$, then $\dim_H \mu_{\varphi}$ is expressed by $f_0^{\prime} (0)$, $f_1^{\prime} (1)$, $g_0^{\prime} (0)$ and $g_1^{\prime} (1)$ (Remark \ref{rem:two}).
If $N \ge 3$, then these identities do not hold in general, and the assumption of Theorem \ref{thm:tiling} may hold even if they fail (Example \ref{exa:non-natural}).
If these identities hold for both of $f$ and $g$, then a measure $\nu$ satisfying \eqref{eq:stationary} is given in terms of $\varphi$, and it is concentrated on a curve in $Y^g \times Y^f$ (Proposition \ref{prop:nu-natural} and Remark \ref{rem:curve}).

Finally, if $f$ and $g$ are hyperbolic, then there exists a Borel probability measure on $[0,1)$ which is $T$-invariant, $T$-ergodic and mutually absolutely continuous with $\mu_{\varphi}$, and the formula in Theorem \ref{thm:dim} (2) states that $\dim_H \mu_{\varphi}$ is the entropy of $T$ with respect to this measure divided by its Lyapunov exponent (Proposition \ref{prop:hatmu} and Remark \ref{rem:entropy}).
See Remark \ref{rem:Gibbs-property} for the Gibbs property of $\mu_{\varphi}$ and of this measure.

This paper is organized as follows.
In Section \ref{sec:pre}, we state notation and lemmas.
In Section \ref{sec:dim}, we show Theorem \ref{thm:dim}, and give a $T$-invariant probability measure which is mutually absolutely continuous with $\mu_{\varphi}$.
In Section \ref{sec:sing}, we show Theorem \ref{thm:sing}.
In Section \ref{sec:para}, we show Theorem \ref{thm:para} and Corollary \ref{cor:positive}.
In Section \ref{sec:tiling}, we show Theorem \ref{thm:tiling} and consider the case that $g_i$, $i \in \mathcal{I}_N$, are conjugate to their dual maps.
In Section \ref{sec:rem}, we give remarks.

\section{Preliminaries}\label{sec:pre}

\subsection{A family of linear fractional transformations}

For $y > -1$ and $x \in [0,1]$, let
\[ k(y,x) \coloneqq  \frac{1+y}{(1+yx)^2}. \]

\begin{Lem}\label{lem:psi}
Let $s > -1$ and $y > -1$. Then,\\
(1) $\psi_s$ is continuously differentiable on $[0,1]$, $\psi_s^{\prime} (x) = k(s,x) > 0$ for $x \in [0,1]$, $\psi_s (0) = 0$ and $\psi_s (1) = 1$.\\
(2) $\gamma_s (y) > -1$ and $\psi_y \circ \psi_s = \psi_{\gamma_s (y)}$ on $[0,1]$.\\
(3) For $x \in (0,1)$,
\[ \frac{\psi_s (x)}{1 - \psi_s (x)} = \frac{(1+s) x}{1 - x}. \]
(4) If $\psi_y (x) = \psi_s (x)$ for some $x \in (0,1)$, then $y = s$.
\end{Lem}

\begin{proof}
(1) Since $s > -1$, $1 + sx > 0$ for $x \in [0,1]$.
The assertions follow from the definition of $\psi_s$.

(2) Since $1 + \gamma_s (y) = (1+s)(1+y) > 0$, the first assertion holds.
For $x \in [0,1]$,
\[ \psi_y (\psi_s (x)) = \frac{(1+y)(1+s) x}{1 + sx + y(1+s) x} = \frac{(1 + \gamma_s (y)) x}{1 + \gamma_s (y) x} = \psi_{\gamma_s (y)} (x). \]

(3) This follows from $1 - \psi_s (x) = (1-x)/(1+sx)$.

(4) This follows from (3).
\end{proof}

Let $s > -1$.
By Lemma \ref{lem:psi} (1), $\psi_s$ is a homeomorphism of $[0,1]$, and, by Lemma \ref{lem:psi} (2), its inverse function is $\psi_{s^{\prime}}$, where $s^{\prime} > -1$ is the number such that $(1+s)(1+s^{\prime}) = 1$.
Since $\gamma_s (y) + 1 = (1+s)(1+y)$ for $y > -1$, $\gamma_s$ is a strictly increasing affine bijection of $(-1, +\infty)$.

\subsection{LF systems}

\begin{Lem}\label{lem:basic}
Let $h = \{h_i\}_{i \in \mathcal{I}_N}$ be an LF system on $[0,1]$, and let $i \in \mathcal{I}_N$. Then,\\
(1) $c_i^h x + 1 > 0$ and $0 < h_i^{\prime} (x) \le 1$ for $x \in [0,1]$.\\
(2) $h_i (0) = b_i^h$ and $h_i (1) = b_{i+1}^h$, and $0 = b_0^h < b_1^h < \cdots < b_N^h = 1$.\\
(3) $a_i^h - b_i^h c_i^h = (b_{i+1}^h - b_i^h)(c_i^h + 1)$, and $h_i (x) = b_i^h + (b_{i+1}^h - b_i^h) \, \psi_{c_i^h} (x)$ for $x \in [0,1]$.\\
(4) $h_0^{\prime} (0) = a_0^h$ and $h_{N-1}^{\prime} (1) = (1 - b_{N-1}^h) / (c_{N-1}^h + 1)$.
\end{Lem}

\begin{proof}
We omit the superscript $h$.

(1) By the conditions (1) and (2) of LF systems for $d_i = 1$, $c_i + 1 > 0$ and $0 < a_i - b_i c_i \le \min\{1, (c_i + 1)^2\}$.
Since $c_i x + 1 \ge \min\{1, c_i + 1\} > 0$ for $x \in [0,1]$, we have that
\[ 0 < h_i^{\prime} (x) = \frac{a_i - b_i c_i}{(c_i x + 1)^2} \le 1, \ \ \ x \in [0,1]. \]

(2) We have that $h_i (0) = b_i$, and, by the condition (3) of LF systems, $b_0 = 0$ and $h_i (1) = (a_i + b_i)/(c_i + 1) = b_{i+1}$, where $b_N = 1$.
By (1), $h_i (0) < h_i (1)$, and hence $b_i < b_{i+1}$.

(3) Since $a_i + b_i = b_{i+1} (c_i + 1)$, the first assertion holds.
By this,
\[ h_i (x) - b_i = \frac{(a_i - b_i c_i) x}{c_i x + 1} = (b_{i+1} - b_i) \, \frac{(1 + c_i) x}{1 + c_i x} = (b_{i+1} - b_i) \, \psi_{c_i} (x). \]

(4) By (3) and Lemma \ref{lem:psi} (1), $h_i^{\prime} (x) = (b_{i+1} - b_i) \, k(c_i, x)$.
Since $b_0 = 0$ and $b_N = 1$, by (3) for $i = 0$, we have that $h_0^{\prime} (0) = b_1 (1 + c_0) = a_0$, and $h_{N-1}^{\prime} (1) = (1 - b_{N-1})/(1 + c_{N-1})$.
\end{proof}

By Lemma \ref{lem:basic} (4), an LF system $h$ is hyperbolic if and only if $a_0^h < 1$ and $b_{N-1}^h + c_{N-1}^h > 0$.
By Lemma \ref{lem:basic} (2) and (3), an LF system $h$ is determined by the numbers $b_i^h$ and $c_i^h$, $i \in \mathcal{I}_N$.

For an LF system $h$ and $i \in \mathcal{I}_N$, let $A_i^h \coloneqq  \begin{pmatrix} a_i^h & b_i^h \\ c_i^h & 1 \end{pmatrix}$.
For a word $\sigma = (\sigma_1, \dots, \sigma_n)$, let
\[ A_{\sigma}^h \coloneqq  A_{\sigma_1}^h \cdots A_{\sigma_n}^h, \ \ \ y_{\sigma}^h \coloneqq  \frac{c_{\sigma}^h}{d_{\sigma}^h}, \]
where $(c_{\sigma}^h, d_{\sigma}^h)$ is the second row of $A_{\sigma}^h$.
If $n = 0$, then $A_{\sigma}^h$ is the identity matrix and $y_{\sigma}^h = 0$.
We see in Lemma \ref{lem:id} (1) below that $d_{\sigma}^h > 0$ and $y_{\sigma}^h > -1$.
For a word $\sigma = (\sigma_1, \dots, \sigma_n)$ and $i \in \mathcal{I}_N$, we denote the word $(\sigma_1, \dots, \sigma_n, i)$ by $\sigma i$.

We recall that, for an LF system $h$, the functions $G_i^h$ and $H_i^h$, $i \in \mathcal{I}_N$, are defined in \eqref{eq:GH}.

\begin{Lem}\label{lem:id}
Let $h = \{h_i\}_{i \in \mathcal{I}_N}$ be an LF system on $[0,1]$.
Let $y > -1$, $x \in [0,1]$, $i \in \mathcal{I}_N$, and let $\sigma$ be a word. Then,\\
(1) $b_i^h y + 1 > 0$, $H_i^h (y) > -1$, $d_{\sigma}^h > 0$, $y_{\sigma}^h > -1$ and $y_{\sigma i}^h = H_i^h (y_{\sigma}^h)$.\\
(2) $h_{\sigma}^{\prime} (x) = \left( h_{\sigma} (1) - h_{\sigma} (0) \right) k(y_{\sigma}^h, x)$.\\
(3) $h_i^{\prime} (x) \, k(y, h_i (x)) = G_i^h (y) \, k(H_i^h (y), x)$.\\
(4) $\psi_y (h_i (x)) = \psi_y (b_i^h) + G_i^h (y) \, \psi_{H_i^h (y)} (x)$.\\
(5) $G_i^h (y) = \psi_y (b_{i+1}^h) - \psi_y (b_i^h) > 0$ and $G_0^h (y) + \cdots + G_{N-1}^h (y) = 1$.\\
(6) $h_{\sigma i} (1) - h_{\sigma i} (0) = G_i^h (y_{\sigma}^h) \left( h_{\sigma} (1) - h_{\sigma} (0) \right)$.
\end{Lem}

\begin{proof}
We omit the superscript $h$.

(1) By Lemma \ref{lem:basic} (2), $0 \le b_i < b_{i+1} \le 1$.
Since $y > -1$, we have that $b_i y + 1 > 0$ and $b_{i+1} y + 1 > 0$.
Since $a_i + b_i = b_{i+1} (c_i + 1)$,
\begin{equation}\label{eq:one-plus-H}
1 + H_i (y) = \frac{(a_i + b_i) y + c_i + 1}{b_i y + 1} = \frac{(c_i + 1)(1 + b_{i+1} y)}{1 + b_i y} > 0.
\end{equation}
Since
\[ (c_{\sigma i}, d_{\sigma i}) = (c_{\sigma}, d_{\sigma}) A_i = (a_i c_{\sigma} + c_i d_{\sigma}, b_i c_{\sigma} + d_{\sigma}), \]
we see by induction on the length of $\sigma$ that $d_{\sigma} > 0$, $y_{\sigma} > -1$ and $y_{\sigma i} = H_i (y_{\sigma})$.

(2) Let $(a_{\sigma}, b_{\sigma})$ be the first row of $A_{\sigma}$.
By (1), $c_{\sigma} x + d_{\sigma} = d_{\sigma} (1 + y_{\sigma} x) > 0$, and $h_{\sigma} (x) = (a_{\sigma} x + b_{\sigma}) / (c_{\sigma} x + d_{\sigma})$.
Hence,
\[ h_{\sigma}^{\prime} (x) = \frac{\det A_{\sigma}}{(c_{\sigma} x + d_{\sigma})^2}, \ \ \ h_{\sigma} (1) - h_{\sigma} (0) = \frac{\det A_{\sigma}}{d_{\sigma} (c_{\sigma} + d_{\sigma})}, \]
and therefore
\[ \frac{h_{\sigma}^{\prime} (x)}{h_{\sigma} (1) - h_{\sigma} (0)} = \frac{d_{\sigma} (c_{\sigma} + d_{\sigma})}{(c_{\sigma} x + d_{\sigma})^2} = k(y_{\sigma}, x). \]

(3) By \eqref{eq:one-plus-H},
\begin{equation}\label{eq:log-derivative}
\frac{d}{dy} \log \left( 1 + H_i (y) \right) = \frac{b_{i+1}}{1 + b_{i+1} y} - \frac{b_i}{1 + b_i y} = \frac{G_i (y)}{1+y}.
\end{equation}
Hence, $(1+y) H_i^{\prime} (y) = G_i (y) \left( 1 + H_i (y) \right)$.
Since
\[ H_i^{\prime} (y) = \frac{a_i - b_i c_i}{(b_i y + 1)^2}, \ \ \ 1 + y \, h_i (x) = \frac{(b_i y + 1)(1 + H_i (y) x)}{c_i x + 1}, \]
we have that
\begin{align*}
h_i^{\prime} (x) \, k(y, h_i (x)) &= \frac{(1+y)(a_i - b_i c_i)}{(b_i y + 1)^2 (1 + H_i (y) x)^2} \\
&= \frac{(1+y) H_i^{\prime} (y)}{(1 + H_i (y) x)^2} = G_i (y) \, k(H_i (y), x).
\end{align*}

(4) By Lemma \ref{lem:psi} (1), the derivatives of both sides with respect to $x$ are the two sides of (3), and, by Lemma \ref{lem:basic} (2), both sides are equal to $\psi_y (b_i)$ at $x = 0$.

(5) By (4) for $x = 1$ and Lemma \ref{lem:basic} (2), $\psi_y (b_{i+1}) - \psi_y (b_i) = G_i (y)$, and, by the definition of $G_i$ and (1), $G_i (y) > 0$.
Since $b_0 = 0$ and $b_N = 1$, the sum of $G_j (y)$ over $j \in \mathcal{I}_N$ is $\psi_y (1) - \psi_y (0) = 1$.

(6) By Lemma \ref{lem:basic} (2), (2), Lemma \ref{lem:psi} (1) and (5),
\begin{align*}
h_{\sigma i} (1) - h_{\sigma i} (0) &= h_{\sigma} (b_{i+1}) - h_{\sigma} (b_i) = \left( h_{\sigma} (1) - h_{\sigma} (0) \right) \int_{b_i}^{b_{i+1}} k(y_{\sigma}, z) \, dz \\
&= \left( h_{\sigma} (1) - h_{\sigma} (0) \right) \left( \psi_{y_{\sigma}} (b_{i+1}) - \psi_{y_{\sigma}} (b_i) \right) \\
&= G_i (y_{\sigma}) \left( h_{\sigma} (1) - h_{\sigma} (0) \right). \qedhere
\end{align*}
\end{proof}

Lemma \ref{lem:id} (5) for $h = g$ contains \cite[Lemma 3.1 (iv)]{O26}.
By Lemma \ref{lem:id} (1), for an LF system $h$ and a word $\sigma = (\sigma_1, \dots, \sigma_n)$,
\[ y_{\sigma}^h = H_{\sigma_n}^h \circ \cdots \circ H_{\sigma_1}^h (0), \]
and the order of this composition is the reverse of that in $h_{\sigma} = h_{\sigma_1} \circ \cdots \circ h_{\sigma_n}$.

\subsection{Fixed points of the dual maps}

Let $h = \{h_i\}_{i \in \mathcal{I}_N}$ be an LF system on $[0,1]$.
The dual maps $H_i^h$, $i \in \mathcal{I}_N$, of $h$ are the linear fractional transformations given by the transposed matrices of $A_i^h$, $i \in \mathcal{I}_N$.
By \eqref{eq:log-derivative} and Lemma \ref{lem:id} (1) and (5), $H_i^h$ is strictly increasing on $(-1, +\infty)$ for every $i \in \mathcal{I}_N$.

\begin{Lem}\label{lem:fixed}
Let $h = \{h_i\}_{i \in \mathcal{I}_N}$ be an LF system on $[0,1]$, and let $i \in \mathcal{I}_N$.
Assume that $h_0^{\prime} (0) < 1$ if $i = 0$, and that $h_{N-1}^{\prime} (1) < 1$ if $i = N-1$.
Then, there exists a unique number $e \in (-1, +\infty)$ such that $H_i^h (e) = e$, and it holds that $H_i^h (y) > y$ for $-1 < y < e$ and that $H_i^h (y) < y$ for $y > e$.
Furthermore,
\[ e = \frac{c_0^h}{1 - a_0^h} \ \textup{ if } i = 0, \ \ \ e = \frac{a_i^h - 1 + \left( (1 - a_i^h)^2 + 4 b_i^h c_i^h \right)^{1/2}}{2 b_i^h} \ \textup{ if } 1 \le i \le N-1, \]
and $e = c_{N-1}^h / b_{N-1}^h$ if $i = N-1$.
\end{Lem}

In the proof, we use that the sign of $H_i^h (y) - y$ is equal to that of a polynomial in $y$ of degree at most two.

\begin{proof}
We omit the superscript $h$.
Let
\[ Q_i (y) \coloneqq  - b_i y^2 + (a_i - 1) y + c_i, \ \ \ y \in \mathbb{R}. \]
Then, $H_i (y) - y = Q_i (y) / (b_i y + 1)$ for $y > -1$, and, by Lemma \ref{lem:id} (1), the sign of $H_i (y) - y$ is equal to that of $Q_i (y)$.
Since $a_i + b_i = b_{i+1} (c_i + 1)$, we have that $Q_i (-1) = (c_i + 1)(1 - b_{i+1})$.

Let $i = 0$.
Then, $b_0 = 0$, and, by the assumption and Lemma \ref{lem:basic} (4), $a_0 < 1$.
Hence, $Q_0$ is strictly decreasing, and it vanishes only at $c_0 / (1 - a_0)$.
Since $Q_0 (-1) = (c_0 + 1)(1 - b_1) > 0$, this number is larger than $-1$.

Let $1 \le i \le N-2$.
Then, $b_i > 0$ and $Q_i (-1) > 0$.
Since $Q_i$ is a quadratic polynomial whose leading coefficient is negative, $Q_i$ has two real zeros, the larger zero, which we denote by $e$, is in $(-1, +\infty)$, the smaller zero is in $(-\infty, -1)$, and it holds that $Q_i > 0$ on $(-1, e)$ and that $Q_i < 0$ on $(e, +\infty)$.
The larger zero is given by the formula in the assertion.

Let $i = N-1$.
Since $a_{N-1} + b_{N-1} = c_{N-1} + 1$,
\[ Q_{N-1} (y) = - (y+1)(b_{N-1} y - c_{N-1}). \]
By the assumption and Lemma \ref{lem:basic} (4), we have that $b_{N-1} + c_{N-1} > 0$, that is, $c_{N-1} / b_{N-1} > -1$.
Hence, for $y > -1$, $Q_{N-1} (y) > 0$ if $y < c_{N-1} / b_{N-1}$, and $Q_{N-1} (y) < 0$ if $y > c_{N-1} / b_{N-1}$.
Since $(1 - a_{N-1})^2 + 4 b_{N-1} c_{N-1} = (b_{N-1} + c_{N-1})^2$, the number $c_{N-1} / b_{N-1}$ is also given by the formula in the assertion for $i = N-1$.
\end{proof}

If $h$ is hyperbolic, then the assumption of Lemma \ref{lem:fixed} holds for every $i \in \mathcal{I}_N$, and the number $e$ in Lemma \ref{lem:fixed} is the number $e_i^h$ in Section \ref{sec:intro}.
If $h$ is parabolic and $i \in \mathcal{I}_N$ satisfies the assumption of Lemma \ref{lem:fixed}, then we also denote the number $e$ in Lemma \ref{lem:fixed} by $e_i^h$.

\subsection{Hyperbolic LF systems}

\begin{Lem}\label{lem:Y}
Let $h$ be a hyperbolic LF system on $[0,1]$. Then,\\
(1) $-1 < \alpha_h \le 0 \le \beta_h < +\infty$.\\
(2) $H_i^h (Y^h) \subset Y^h$ for every $i \in \mathcal{I}_N$.
\end{Lem}

In the proof, we use Lemma \ref{lem:fixed} and that the dual maps are strictly increasing.

\begin{proof}
(1) By Lemma \ref{lem:fixed}, $e_i^h > -1$ for every $i \in \mathcal{I}_N$.
Hence, the assertion follows from the definitions of $\alpha_h$ and $\beta_h$.

(2) Let $i \in \mathcal{I}_N$ and $y \in Y^h$.
By (1) and the definitions of $\alpha_h$ and $\beta_h$, $-1 < \alpha_h \le e_i^h \le \beta_h$, and hence, by Lemma \ref{lem:fixed}, $H_i^h (\alpha_h) \ge \alpha_h$ and $H_i^h (\beta_h) \le \beta_h$.
Since $H_i^h$ is strictly increasing on $(-1, +\infty)$,
\[ \alpha_h \le H_i^h (\alpha_h) \le H_i^h (y) \le H_i^h (\beta_h) \le \beta_h. \qedhere \]
\end{proof}

For $h = g$, $Y^g$ is the set $Y$ in \cite[Section 3]{O26}, and the assertions of Lemma \ref{lem:Y} are contained in \cite[Lemma 3.1 and Lemma 3.14 (v)]{O26}.

For a hyperbolic LF system $h$, let
\[ C_h \coloneqq  \max\left\{ 1 + \beta_h, \frac{1}{1 + \alpha_h} \right\}, \]
\[ G_{\min}^h \coloneqq  \min_{i \in \mathcal{I}_N} \min_{y \in Y^h} G_i^h (y), \ \ \ G_{\max}^h \coloneqq  \max_{i \in \mathcal{I}_N} \max_{y \in Y^h} G_i^h (y). \]

\begin{Lem}\label{lem:hyp}
Let $h$ be a hyperbolic LF system on $[0,1]$.
Let $\sigma$ be a word of length $n \ge 0$, $i \in \mathcal{I}_N$ and $x \in [0,1]$. Then,\\
(1) $y_{\sigma}^h \in Y^h$.\\
(2) $C_h^{-1} \le k(y,x) \le C_h$ for $y \in Y^h$.\\
(3) $C_h^{-1} \left( h_{\sigma} (1) - h_{\sigma} (0) \right) \le h_{\sigma}^{\prime} (x) \le C_h \left( h_{\sigma} (1) - h_{\sigma} (0) \right)$.\\
(4) $0 < G_{\min}^h \le G_{\max}^h < 1$,
\[ h_{\sigma i} (1) - h_{\sigma i} (0) \ge G_{\min}^h \left( h_{\sigma} (1) - h_{\sigma} (0) \right), \]
and
\[ \left( G_{\min}^h \right)^n \le h_{\sigma} (1) - h_{\sigma} (0) \le \left( G_{\max}^h \right)^n. \]
\end{Lem}

\begin{proof}
(1) By Lemma \ref{lem:Y}, $0 \in Y^h$ and $H_j^h (Y^h) \subset Y^h$ for every $j \in \mathcal{I}_N$.
Hence, the assertion follows from Lemma \ref{lem:id} (1) and induction on $n$.

(2) If $y \ge 0$, then $1 \le 1 + yx \le 1 + y$ and hence $(1+y)^{-1} \le k(y,x) \le 1 + y$.
If $-1 < y < 0$, then $1 + y \le 1 + yx \le 1$ and hence $1 + y \le k(y,x) \le (1+y)^{-1}$.
Since $-1 < \alpha_h \le y \le \beta_h$ for $y \in Y^h$ by Lemma \ref{lem:Y} (1), the assertion holds.

(3) This follows from (1), (2) and Lemma \ref{lem:id} (2).

(4) By Lemma \ref{lem:id} (5), $G_j^h$, $j \in \mathcal{I}_N$, are positive and continuous on the compact set $Y^h$ and their sum is $1$.
Since $N \ge 2$, we have that $0 < G_{\min}^h \le G_{\max}^h < 1$.
The inequalities follow from this, (1), Lemma \ref{lem:id} (6) and induction on $n$.
\end{proof}

\subsection{Intervals and maps for the conjugate equation}

\begin{Lem}\label{lem:cyl}
Let $h = \{h_i\}_{i \in \mathcal{I}_N}$ be an LF system on $[0,1]$. Then,\\
(1) For every word $\sigma$, the interval $[h_{\sigma} (0), h_{\sigma} (1))$ is the disjoint union of the intervals $[h_{\sigma i} (0), h_{\sigma i} (1))$, $i \in \mathcal{I}_N$.\\
(2) For every $n \ge 0$, the intervals $[h_{\sigma} (0), h_{\sigma} (1))$, $\sigma \in (\mathcal{I}_N)^n$, are pairwise disjoint and their union is $[0,1)$.\\
(3) For two words $\sigma$ and $\sigma^{\prime}$, the intervals $[h_{\sigma} (0), h_{\sigma} (1))$ and $[h_{\sigma^{\prime}} (0), h_{\sigma^{\prime}} (1))$ are disjoint or one of them contains the other.
\end{Lem}

\begin{proof}
By Lemma \ref{lem:basic}, $h_{\sigma}$ is continuous and strictly increasing on $[0,1]$, $h_{\sigma i} (0) = h_{\sigma} (b_i^h)$ and $h_{\sigma i} (1) = h_{\sigma} (b_{i+1}^h)$, and $0 = b_0^h < b_1^h < \cdots < b_N^h = 1$.
Hence (1) holds.
Assertions (2) and (3) follow from (1) and induction on the lengths of words.
\end{proof}

{\it In the rest of this section, let $f$ and $g$ be LF systems on $[0,1]$, and let $\varphi$ be the solution of \eqref{eq:conj}.} 
We recall that $J(\sigma) = [f_{\sigma}(0), f_{\sigma}(1))$ for a word $\sigma$.

\begin{Lem}\label{lem:conj}
(1) $\varphi$ is a homeomorphism of $[0,1]$, and $\varphi (f_{\sigma} (x)) = g_{\sigma} (\varphi (x))$ for every word $\sigma$ and every $x \in [0,1]$.\\
(2) $\varphi (J(\sigma)) = [g_{\sigma} (0), g_{\sigma} (1))$, $\mu_{\varphi} (J(\sigma)) = g_{\sigma} (1) - g_{\sigma} (0)$ and $\ell (J(\sigma)) = f_{\sigma} (1) - f_{\sigma} (0)$ for every word $\sigma$.\\
(3) $\mu_{\varphi} (B) = \ell (\varphi (B))$ for every Borel subset $B$ of $[0,1]$.
\end{Lem}

\begin{proof}
(1) Since $\varphi$ is continuous and strictly increasing, $\varphi (0) = 0$ and $\varphi (1) = 1$, the first assertion holds.
The second assertion follows from \eqref{eq:conj} and induction on the length of $\sigma$.

(2) By (1), $\varphi (f_{\sigma} (0)) = g_{\sigma} (0)$ and $\varphi (f_{\sigma} (1)) = g_{\sigma} (1)$.
Since $\varphi$ is the distribution function of $\mu_{\varphi}$ and is continuous, the assertions hold.

(3) Since $\mu_{\varphi} ([0,x]) = \varphi (x) = \ell (\varphi ([0,x]))$ for every $x \in [0,1]$, the assertion follows from the Dynkin $\pi$-$\lambda$ theorem.
\end{proof}

By Lemma \ref{lem:conj} (2) and Lemma \ref{lem:id} (6) for $h = g$ and for $h = f$, for every word $\sigma$ and every $i \in \mathcal{I}_N$,
\begin{equation}\label{eq:next-letter}
\frac{\mu_{\varphi} (J(\sigma i))}{\mu_{\varphi} (J(\sigma))} = G_i^g (y_{\sigma}^g), \ \ \ \frac{\ell (J(\sigma i))}{\ell (J(\sigma))} = G_i^f (y_{\sigma}^f).
\end{equation}
These are the conditional probabilities of $J(\sigma i)$ given $J(\sigma)$ under $\mu_{\varphi}$ and under $\ell$, respectively.
In Theorem \ref{thm:dim} (2), the point $x$ is chosen according to $\mu_{\varphi}$.
Hence, the letter $i$ follows the word $\sigma$ with probability $G_i^g (y_{\sigma}^g)$, and then $- \log \mu_{\varphi} (J(\sigma))$ and $- \log \ell (J(\sigma))$ increase by $- \log G_i^g (y_{\sigma}^g)$ and $- \log G_i^f (y_{\sigma}^f)$, respectively.
The averages of these two increments with respect to the probabilities $G_i^g (y_{\sigma}^g)$, $i \in \mathcal{I}_N$, are $S(y_{\sigma}^g)$ and $\Lambda (y_{\sigma}^g, y_{\sigma}^f)$, respectively.
This is the reason why the pairs $(y_{\sigma}^g, y_{\sigma}^f)$ and the functions $S$ and $\Lambda$ appear in Theorem \ref{thm:dim}.

By Lemma \ref{lem:cyl} for $h = f$, for every $n \ge 0$, the intervals $J(\sigma)$, $\sigma \in (\mathcal{I}_N)^n$, are pairwise disjoint and their union is $[0,1)$.
Let $T$, $X_1$ and $X_n$, $n \ge 1$, be the maps in Section \ref{sec:intro}, and let $\tau : [0,1) \to [0,1)$ be the map defined by
\[ \tau (z) \coloneqq  g_i^{-1} (z), \ \ \ b_i^g \le z < b_{i+1}^g, \ \ i \in \mathcal{I}_N. \]
For $j \ge 0$, we denote the $j$-fold composition of $\tau$ by $\tau^j$. 
We recall that, for $n \ge 1$ and $x \in [0,1)$, $J_n(x)$ is the interval $J(\sigma)$, $\sigma \in (\mathcal{I}_N)^n$, which contains $x$. 

\begin{Lem}\label{lem:coding}
Let $x \in [0,1)$ and $n \ge 1$.
Let $\sigma_j \coloneqq  X_j (x)$ for $1 \le j \le n$, and let $\sigma \coloneqq  (\sigma_1, \dots, \sigma_n)$. Then,\\
(1) $T^{j-1} (x) = f_{\sigma_j} (T^j (x))$ for $1 \le j \le n$, $x = f_{\sigma} (T^n (x))$ and $J_n (x) = J(\sigma)$.\\
(2) $\varphi (T^j (x)) = \tau^j (\varphi (x))$ for $0 \le j \le n$, and $\tau^{j-1} (\varphi (x)) = g_{\sigma_j} \left( \tau^j (\varphi (x)) \right)$ for $1 \le j \le n$.\\
(3) $\tau^n (g_{\sigma^{\prime}} (z)) = z$ for every word $\sigma^{\prime}$ of length $n$ and every $z \in [0,1)$.
\end{Lem}

\begin{proof}
(1) Let $1 \le j \le n$.
Since $X_1 (T^{j-1} (x)) = \sigma_j$, by the definition of $T$, $T^j (x) = f_{\sigma_j}^{-1} (T^{j-1} (x))$.
Hence the first assertion holds, and $x = f_{\sigma} (T^n (x))$.
Since $T^n (x) \in [0,1)$, we have that $x \in J(\sigma)$, and hence $J_n (x) = J(\sigma)$.

(2) Let $x^{\prime} \in [0,1)$ and $i \coloneqq  X_1 (x^{\prime})$.
Then, $x^{\prime} = f_i (T(x^{\prime}))$, and, by \eqref{eq:conj}, $\varphi (x^{\prime}) = g_i (\varphi (T(x^{\prime})))$.
Since $\varphi (T(x^{\prime})) \in [0,1)$, by Lemma \ref{lem:basic} (2), $\varphi (x^{\prime}) \in [b_i^g, b_{i+1}^g)$, and hence $\tau (\varphi (x^{\prime})) = \varphi (T(x^{\prime}))$.
By this for $x^{\prime} = T^{j-1} (x)$, $1 \le j \le n$, and induction on $j$, the first assertion holds.
By (1) and \eqref{eq:conj}, $\varphi (T^{j-1} (x)) = g_{\sigma_j} (\varphi (T^j (x)))$ for $1 \le j \le n$.
Hence, by the first assertion, the second assertion holds.

(3) Let $i \in \mathcal{I}_N$ and $z \in [0,1)$.
By Lemma \ref{lem:basic} (2), $g_i (z) \in [b_i^g, b_{i+1}^g)$, and hence $\tau (g_i (z)) = z$.
The assertion follows from this and induction on $n$.
\end{proof}

The notation for the two sides of \eqref{eq:conj} is summarized as follows, where $\sigma$ is a word.

\begin{center}
\begin{tabular}{l|l|l}
 & source side ($f$) & target side ($g$) \\ \hline
variable & $x$ & $z = \varphi (x)$ \\
map & $T$ & $\tau$ \\
inverse branches & $f_i$ & $g_i$ \\
interval & $J(\sigma) = [f_{\sigma} (0), f_{\sigma} (1))$ & $\varphi (J(\sigma)) = [g_{\sigma} (0), g_{\sigma} (1))$ \\
dual variable & $w$, \ $y_{\sigma}^f$ & $y$, \ $y_{\sigma}^g$ \\
measures & $\ell$, \ $\mu_{\varphi}$, \ $\widehat{\mu}$ & $\mu_{\varphi^{-1}}$, \ $\ell$, \ $m$
\end{tabular}
\end{center}

In the last row, the three measures on the target side are the image measures under $\varphi$ of the three measures on the source side, in this order, where $\ell$ is the Lebesgue measure on each side.
The measures $m$, $\widehat{\mu}$ and $\mu_{\varphi^{-1}}$ are defined in Subsection \ref{subsec:tau}, Subsection \ref{subsec:invariant} and Remark \ref{rem:inverse}, respectively.
Here, the dual variables are the variables for the dual maps of $f$ and of $g$.
The measure $\nu$ in Theorem \ref{thm:dim} is a measure for the pairs $(y,w)$ of the dual variables, and it is different from these measures.

\subsection{Hausdorff dimensions}

The following corresponds to \cite[Lemma 7]{O20}; see also \cite[Lemma 3.5]{O26}.
By Lemmas \ref{lem:cyl} and \ref{lem:conj}, $0 < \mu_{\varphi} (J_n (x)) < 1$ and $0 < \ell (J_n (x)) < 1$ for $n \ge 1$ and $x \in [0,1)$.

\begin{Lem}\label{lem:dimension}
Assume that $f$ is hyperbolic.
Let $c \ge 0$. Then,\\
(1) If
\[ \limsup_{n \to \infty} \frac{\log \mu_{\varphi} (J_n (x))}{\log \ell (J_n (x))} \le c, \ \ \ \textup{$\mu_{\varphi}$-a.s.} \, x, \]
then, there exists a Borel set $K_0$ such that $\mu_{\varphi} (K_0) = 1$ and $\dim_H K_0 \le c$.\\
(2) If
\[ \liminf_{n \to \infty} \frac{\log \mu_{\varphi} (J_n (x))}{\log \ell (J_n (x))} \ge c, \ \ \ \textup{$\mu_{\varphi}$-a.s.} \, x, \]
then, $\mu_{\varphi} (K) = 0$ for every Borel set $K$ with $\dim_H K < c$.
\end{Lem}

For (2), we adapt the proof of the lower bound in Billingsley's lemma \cite[Lemma 1.4.1]{BP17}, and we choose the intervals $J_n (x)$ as in the proof of \cite[Theorem 2.2.2]{BP17}.

\begin{proof}
By Lemma \ref{lem:hyp} (4) for $h = f$ and Lemma \ref{lem:conj} (2), for every word $\sigma$ of length $n$ and every $i \in \mathcal{I}_N$,
\begin{equation}\label{eq:length}
\ell (J(\sigma i)) \ge G_{\min}^f \, \ell (J(\sigma)), \ \ \ \left( G_{\min}^f \right)^n \le \ell (J(\sigma)) \le \left( G_{\max}^f \right)^n.
\end{equation}

(1) Let $K_0$ be the set of the points $x \in [0,1)$ for which the limit superior in the assumption is at most $c$.
It is a Borel set, since, for each $n$, the ratio is constant on each $J(\sigma)$, $\sigma \in (\mathcal{I}_N)^n$, and $\mu_{\varphi} (K_0) = 1$ by the assumption.
Let $s > c$, and, for $m \ge 1$, let
\[ E_m \coloneqq  \left\{ x \in [0,1) : \mu_{\varphi} (J_n (x)) \ge \ell (J_n (x))^s \textup{ for every } n \ge m \right\}. \]
Since $\log \ell (J_n (x)) < 0$, $K_0$ is contained in the union of the sets $E_m$, $m \ge 1$.
For each $n \ge m$, the intervals $J(\sigma)$, $\sigma \in (\mathcal{I}_N)^n$, intersecting $E_m$ cover $E_m$, and the sum of $\ell (J(\sigma))^s$ over them is at most that of $\mu_{\varphi} (J(\sigma))$, which is at most $1$.
By \eqref{eq:length}, the maximum of their lengths converges to $0$ as $n \to \infty$, and hence $\dim_H E_m \le s$.
Thus $\dim_H K_0 \le s$ for every $s > c$, and (1) holds.

(2) Let $K$ be a nonempty Borel set with $\dim_H K < c$, and let $\dim_H K < s < c$.
As in (1), by the assumption, the $\mu_{\varphi}$-measure of the union of the sets
\[ F_m \coloneqq  \left\{ x \in [0,1) : \mu_{\varphi} (J_n (x)) \le \ell (J_n (x))^s \textup{ for every } n \ge m \right\}, \ \ \ m \ge 1, \]
is $1$.
Let $m \ge 1$, and let $I$ be an interval of length $\delta$ with $0 < \delta < \left( G_{\min}^f \right)^m$.
For each $x \in [0,1)$, we take the interval $J_n (x)$ for the smallest $n \ge 1$ such that $\ell (J_n (x)) < \delta$.
By Lemma \ref{lem:cyl} and \eqref{eq:length}, the distinct intervals taken in this manner are pairwise disjoint, and each of them is $J(\sigma)$ for a word $\sigma$ of length larger than $m$ and has length at least $G_{\min}^f \, \delta$.
Those intersecting $I$ are contained in the interval of length $3 \delta$ whose midpoint is that of $I$, and hence their number is at most $3 / G_{\min}^f$.
The $\mu_{\varphi}$-measure of each of them intersecting $F_m$ is at most the $s$-th power of its length, which is smaller than $\delta^s$.
Hence, $\mu_{\varphi} (I \cap F_m) \le 3 \delta^s / G_{\min}^f$.
Since $\dim_H K < s$, for every $\epsilon > 0$, $K$ is covered by countably many intervals of positive lengths smaller than $\left( G_{\min}^f \right)^m$ such that the sum of the $s$-th powers of their lengths is smaller than $\epsilon$.
Thus $\mu_{\varphi} (K \cap F_m) = 0$ for every $m \ge 1$, and (2) holds.
\end{proof}

\section{Proof of Theorem \ref{thm:dim}}\label{sec:dim}

{\it Throughout this section, we assume that $f$ and $g$ are hyperbolic LF systems on $[0,1]$, and let $\varphi$ be the solution of \eqref{eq:conj}.}

The proof consists of four steps.
We first show Theorem \ref{thm:dim} (1) (Lemma \ref{lem:stationary}).
Next, from the first marginal of a measure $\nu$ satisfying \eqref{eq:stationary}, we construct a Borel probability measure $m$ on $[0,1)$ which is invariant and ergodic with respect to the map $\tau$ in Section \ref{sec:pre} and is mutually absolutely continuous with $\ell$ (Proposition \ref{prop:tau}).
Then, by the Birkhoff ergodic theorem for $\tau$ and $m$, we show that the two limits in Theorem \ref{thm:dim} (2) exist and are given by integrals with respect to $m$ (Proposition \ref{prop:limits}).
Finally, we show that these integrals are equal to the integrals of $S$ and of $\Lambda$ with respect to $\nu$ (Proposition \ref{prop:average}), and obtain $\dim_H \mu_{\varphi}$ by Lemma \ref{lem:dimension}.
After the proof, we give a $T$-invariant measure with respect to which these two limits are the entropy and the Lyapunov exponent of $T$ (Subsection \ref{subsec:invariant}).

\subsection{A stationary measure}

\begin{Lem}\label{lem:stationary}
There exists a Borel probability measure $\nu$ on $Y^g \times Y^f$ such that \eqref{eq:stationary} holds for every continuous function $u$ on $Y^g \times Y^f$.
\end{Lem}

The following proof is the argument in the proof of the Krylov--Bogolyubov theorem \cite[Corollary 6.9.1]{W82}. 
We recall that, for $h \in \{f, g\}$ and a word $\sigma = (\sigma_1, \dots, \sigma_n)$, $y_{\sigma}^h = H_{\sigma_n}^h \circ \cdots \circ H_{\sigma_1}^h (0)$, which is $0$ if $n = 0$. 
We take a weak limit of the averages of the distributions of the points $(y_{\sigma}^g, y_{\sigma}^f)$, where the words $\sigma$ of each length are chosen with the probabilities $\mu_{\varphi} (J(\sigma))$.

\begin{proof}
For a continuous function $u$ on $Y^g \times Y^f$, let
\[ P u (y,w) \coloneqq  \sum_{i \in \mathcal{I}_N} G_i^g (y) \, u \! \left( H_i^g (y), H_i^f (w) \right), \ \ \ (y,w) \in Y^g \times Y^f. \]
Since $H_i^g (Y^g) \subset Y^g$ and $H_i^f (Y^f) \subset Y^f$, $P u$ is a continuous function on $Y^g \times Y^f$.
For $n \ge 1$, let $\nu_n$ be the Borel probability measure on $Y^g \times Y^f$ defined by
\[ \int_{Y^g \times Y^f} u \, d\nu_n \coloneqq  \frac{1}{n} \sum_{j=0}^{n-1} \sum_{\sigma \in (\mathcal{I}_N)^j} \mu_{\varphi} (J(\sigma)) \, u \! \left( y_{\sigma}^g, y_{\sigma}^f \right). \]
This is well-defined, since $(y_{\sigma}^g, y_{\sigma}^f) \in Y^g \times Y^f$ by Lemma \ref{lem:hyp} (1) and the sum of $\mu_{\varphi} (J(\sigma))$ over $\sigma \in (\mathcal{I}_N)^j$ is $1$ for every $j \ge 0$.
By \eqref{eq:next-letter}, $\mu_{\varphi} (J(\sigma i)) = G_i^g (y_{\sigma}^g) \, \mu_{\varphi} (J(\sigma))$ for every word $\sigma$ and every $i \in \mathcal{I}_N$.
By this and Lemma \ref{lem:id} (1), for every $j \ge 0$,
\begin{align*}
\sum_{\sigma \in (\mathcal{I}_N)^j} \mu_{\varphi} (J(\sigma)) \, P u \! \left( y_{\sigma}^g, y_{\sigma}^f \right) &= \sum_{\sigma \in (\mathcal{I}_N)^j} \sum_{i \in \mathcal{I}_N} \mu_{\varphi} (J(\sigma i)) \, u \! \left( y_{\sigma i}^g, y_{\sigma i}^f \right) \\
&= \sum_{\sigma \in (\mathcal{I}_N)^{j+1}} \mu_{\varphi} (J(\sigma)) \, u \! \left( y_{\sigma}^g, y_{\sigma}^f \right).
\end{align*}
Hence,
\[ \int_{Y^g \times Y^f} P u \, d\nu_n - \int_{Y^g \times Y^f} u \, d\nu_n = \frac{1}{n} \left( \sum_{\sigma \in (\mathcal{I}_N)^n} \mu_{\varphi} (J(\sigma)) \, u \! \left( y_{\sigma}^g, y_{\sigma}^f \right) - u(0,0) \right),\] 
and the absolute value of the right hand side is at most $2 \| u \| / n$, where $\| u \|$ denotes the maximum of $|u|$ on $Y^g \times Y^f$.
Since $Y^g \times Y^f$ is compact, by Prokhorov's theorem, there exist a subsequence $(\nu_{n_k})_{k \ge 1}$ and a Borel probability measure $\nu$ on $Y^g \times Y^f$ such that $(\nu_{n_k})_{k \ge 1}$ converges weakly to $\nu$.
By letting $k \to \infty$, we have that
\[ \int_{Y^g \times Y^f} P u \, d\nu = \int_{Y^g \times Y^f} u \, d\nu \]
for every continuous function $u$ on $Y^g \times Y^f$.
\end{proof}

\subsection{An invariant measure for a piecewise linear fractional map}\label{subsec:tau}

{\it In the rest of this section, let $\nu$ be a Borel probability measure on $Y^g \times Y^f$ satisfying \eqref{eq:stationary} for every continuous function $u$ on $Y^g \times Y^f$.}

Let $\nu^g$ be the first marginal of $\nu$, that is, the Borel probability measure on $Y^g$ defined by $\nu^g (B) \coloneqq  \nu (B \times Y^f)$ for Borel subsets $B$ of $Y^g$.
Then,
\[ \int_{Y^g \times Y^f} v(y) \, \nu(dy \, dw) = \int_{Y^g} v \, d\nu^g \]
for every bounded Borel measurable function $v$ on $Y^g$.
We use $\nu^g$ for the integrals with respect to $\nu$ of functions which do not depend on the second coordinate $w$.
By \eqref{eq:stationary} for the functions $(y,w) \mapsto v(y)$,
\begin{equation}\label{eq:stationary-marginal}
\int_{Y^g} \sum_{i \in \mathcal{I}_N} G_i^g (y) \, v(H_i^g (y)) \, \nu^g (dy) = \int_{Y^g} v \, d\nu^g
\end{equation}
for every continuous function $v$ on $Y^g$.
This equation does not involve $f$, and it states that $\nu^g$ is a stationary distribution of the transition kernel on $Y^g$ by which $y$ moves to $H_i^g (y)$ with probability $G_i^g (y)$, $i \in \mathcal{I}_N$.
In contrast, \eqref{eq:stationary} for the functions $(y,w) \mapsto v(w)$ involves both coordinates, since the weights $G_i^g (y)$ depend on the first coordinate $y$, and it is not an equation for the second marginal of $\nu$ alone.

Define a function $\rho : [0,1] \to (0, +\infty)$ by
\[ \rho (z) \coloneqq  \int_{Y^g} k(y,z) \, \nu^g (dy). \]

\begin{Lem}\label{lem:rho}
(1) $C_g^{-1} \le \rho (z) \le C_g$ for $z \in [0,1]$.\\
(2) The integral of $\rho$ over $[0,1]$ with respect to the Lebesgue measure is $1$.\\
(3) For $x \in [0,1]$,
\[ \sum_{i \in \mathcal{I}_N} g_i^{\prime} (x) \, \rho (g_i (x)) = \rho (x). \]
\end{Lem}

\begin{proof}
(1) This follows from Lemma \ref{lem:hyp} (2) for $h = g$.

(2) By Lemma \ref{lem:psi} (1),
\begin{equation}\label{eq:k-int}
\int_0^1 k(y,z) \, dz = \psi_y (1) - \psi_y (0) = 1, \ \ \ y > -1,
\end{equation}
and hence, by Fubini's theorem, the assertion holds.

(3) By the definition of $\rho$, Lemma \ref{lem:id} (3) for $h = g$, and \eqref{eq:stationary-marginal} for the continuous function $y \mapsto k(y,x)$ on $Y^g$,
\begin{align*}
\sum_{i \in \mathcal{I}_N} g_i^{\prime} (x) \, \rho (g_i (x)) &= \int_{Y^g} \sum_{i \in \mathcal{I}_N} G_i^g (y) \, k(H_i^g (y), x) \, \nu^g (dy) \\
&= \int_{Y^g} k(y,x) \, \nu^g (dy) = \rho (x). \qedhere
\end{align*}
\end{proof}

Let $m$ be the Borel measure on $[0,1)$ defined by $m(dz) \coloneqq  \rho (z) \, dz$.
By Lemma \ref{lem:rho} (2), $m$ is a probability measure, and, by Lemma \ref{lem:rho} (1), $m$ and $\ell$ are mutually absolutely continuous.

\begin{Prop}\label{prop:tau}
$m$ is invariant and ergodic with respect to $\tau$.
\end{Prop}

The invariance follows from Lemma \ref{lem:rho} (3) and a change of variables.
In the proof of the ergodicity, for a Borel set $B$ such that $\tau^{-1} (B) = B$, we show by Lemma \ref{lem:hyp} (3) for $h = g$ that the conditional expectation of the indicator function of $B$ with respect to $\ell$ given the intervals $[g_{\sigma} (0), g_{\sigma} (1))$, $\sigma \in (\mathcal{I}_N)^n$, is at least $C_g^{-1} \ell (B)$, and we apply the martingale convergence theorem.
The general case is reduced to this case.

\begin{proof}
We show the invariance.
Let $v$ be a bounded Borel measurable function on $[0,1)$.
By Lemma \ref{lem:basic} (2) and Lemma \ref{lem:rho} (3),
\begin{align*}
\int_0^1 v(\tau (z)) \, \rho (z) \, dz &= \sum_{i \in \mathcal{I}_N} \int_{b_i^g}^{b_{i+1}^g} v\left( g_i^{-1} (z) \right) \rho (z) \, dz \\
&= \sum_{i \in \mathcal{I}_N} \int_0^1 v(x) \, \rho (g_i (x)) \, g_i^{\prime} (x) \, dx = \int_0^1 v(x) \, \rho (x) \, dx.
\end{align*}

We next show the ergodicity.
Let $B$ be a Borel subset of $[0,1)$ such that $m\left( \tau^{-1} (B) \triangle B \right) = 0$.
We show that $m(B) \in \{0,1\}$.

We first consider the case that $\tau^{-1} (B) = B$.
For a word $\sigma$, let $C(\sigma) \coloneqq  [g_{\sigma} (0), g_{\sigma} (1))$.
Let $n \ge 1$ and let $\sigma$ be a word of length $n$.
By Lemma \ref{lem:coding} (3), $\tau^n (g_{\sigma} (z)) = z$ for $z \in [0,1)$.
Since $\tau^{-n} (B) = B$, we have that $g_{\sigma} (B) \subset B$.
Since $g_{\sigma} ([0,1)) = C(\sigma)$, we have that $g_{\sigma} (B) \subset B \cap C(\sigma)$.
By this and Lemma \ref{lem:hyp} (3) for $h = g$,
\[ \ell (B \cap C(\sigma)) \ge \ell (g_{\sigma} (B)) = \int_B g_{\sigma}^{\prime} (z) \, dz \ge C_g^{-1} \, \ell (C(\sigma)) \, \ell (B). \]
Let $\mathcal{G}_n$ be the $\sigma$-algebra on $[0,1)$ generated by $\{C(\sigma^{\prime}) : \sigma^{\prime} \in (\mathcal{I}_N)^n\}$.
Then, by Lemma \ref{lem:cyl} for $h = g$, the conditional expectation of the indicator function $1_B$ of $B$ with respect to $\ell$ satisfies that $E^{\ell} [1_B | \mathcal{G}_n] \ge C_g^{-1} \ell(B)$ on $[0,1)$, and $\mathcal{G}_n \subset \mathcal{G}_{n+1}$.
By Lemma \ref{lem:hyp} (4) for $h = g$, the lengths of $C(\sigma^{\prime})$, $\sigma^{\prime} \in (\mathcal{I}_N)^n$, are at most $\left( G_{\max}^g \right)^n$, and hence the Borel $\sigma$-algebra on $[0,1)$ is generated by $\mathcal{G}_n$, $n \ge 1$.
By the martingale convergence theorem,
\[ \lim_{n \to \infty} E^{\ell} [1_B | \mathcal{G}_n] = 1_B, \textup{ $\ell$-a.e.} \]
Hence, $1_B \ge C_g^{-1} \ell(B)$, $\ell$-a.e.
Therefore, if $\ell (B) > 0$, then $\ell (B) = 1$.
Since $m$ and $\ell$ are mutually absolutely continuous, $m(B) \in \{0,1\}$.

We now consider the general case.
Let
\[ B_{\infty} \coloneqq  \bigcap_{k \ge 1} \bigcup_{n \ge k} \tau^{-n} (B). \]
Then, $\tau^{-1} (B_{\infty}) = B_{\infty}$.
By the invariance of $m$, for every $n \ge 1$,
\[ m\left( \tau^{-n} (B) \triangle B \right) \le \sum_{j=0}^{n-1} m\left( \tau^{-j} \left( \tau^{-1} (B) \triangle B \right) \right) = 0. \]
Since $B_{\infty} \triangle B$ is contained in the union of the sets $\tau^{-n} (B) \triangle B$, $n \ge 1$, we have that $m(B_{\infty} \triangle B) = 0$.
Hence, by the above case, $m(B) = m(B_{\infty}) \in \{0,1\}$.
\end{proof}

\subsection{The two limits}

Let $\Psi_g$ and $\Psi_f$ be the functions on $[0,1)$ defined by
\[ \Psi_g (z) \coloneqq  - \log g_i^{\prime} (\tau (z)), \ \ \ \Psi_f (z) \coloneqq  - \log f_i^{\prime} \left( \varphi^{-1} (\tau (z)) \right) \]
for $z \in [b_i^g, b_{i+1}^g)$ and $i \in \mathcal{I}_N$.
By Lemma \ref{lem:basic} (1), $f_i^{\prime}$ and $g_i^{\prime}$ are positive and continuous on $[0,1]$, and hence $\Psi_g$ and $\Psi_f$ are bounded Borel measurable functions.
These are functions of the variable $z$ on the target side, and they are expressed in terms of the variable $x$ on the source side as follows.
Let $x \in [0,1)$ and $i \coloneqq  X_1 (x)$.
By Lemma \ref{lem:coding} (2), $\tau (\varphi (x)) = \varphi (T(x))$ and $\varphi (x) = g_i (\tau (\varphi (x)))$, and hence, by Lemma \ref{lem:basic} (2), $\varphi (x) \in [b_i^g, b_{i+1}^g)$.
Hence, by the definitions of $\Psi_g$ and $\Psi_f$ and the definition of $\omega_f$ in Section \ref{sec:intro},
\begin{equation}\label{eq:Psi-omega}
\Psi_g (\varphi (x)) = - \log g_{X_1 (x)}^{\prime} (\varphi (T(x))), \ \ \ \Psi_f (\varphi (x)) = - \omega_f (x), \ \ \ x \in [0,1).
\end{equation}
Let
\begin{equation}\label{eq:lambda}
\lambda_g \coloneqq  \int_{[0,1)} \Psi_g \, dm, \ \ \ \lambda_f \coloneqq  \int_{[0,1)} \Psi_f \, dm.
\end{equation}

\begin{Prop}\label{prop:limits}
It holds that
\[ \lim_{n \to \infty} \frac{- \log \mu_{\varphi} (J_n (x))}{n} = \lambda_g, \ \ \ \lim_{n \to \infty} \frac{- \log \ell (J_n (x))}{n} = \lambda_f, \ \ \ \textup{$\mu_{\varphi}$-a.s.} \, x. \]
\end{Prop}

In the proof, by Lemma \ref{lem:hyp} (3), we compare the logarithms of $\mu_{\varphi} (J_n (x))$ and of $\ell (J_n (x))$ with the sums of $- \Psi_g$ and of $- \Psi_f$ along the orbit of $\varphi (x)$ under $\tau$, respectively, and we apply the Birkhoff ergodic theorem for $\tau$ and $m$.
Since $m$ and $\ell$ are mutually absolutely continuous and $\ell$ is the image measure of $\mu_{\varphi}$ under $\varphi$, the limits hold for $\mu_{\varphi}$-a.s.\ $x$.

\begin{proof}
By Proposition \ref{prop:tau} and the Birkhoff ergodic theorem, there exists a Borel subset $\mathcal{N}$ of $[0,1)$ such that $m (\mathcal{N}) = 0$ and
\begin{equation}\label{eq:Birkhoff}
\lim_{n \to \infty} \frac{1}{n} \sum_{j=0}^{n-1} \Psi_g (\tau^j (z)) = \lambda_g, \ \ \ \lim_{n \to \infty} \frac{1}{n} \sum_{j=0}^{n-1} \Psi_f (\tau^j (z)) = \lambda_f
\end{equation}
for every $z \in [0,1) \setminus \mathcal{N}$.
Since $m$ and $\ell$ are mutually absolutely continuous, $\ell (\mathcal{N}) = 0$.

Let $x \in [0,1)$, $n \ge 1$ and $\sigma = (\sigma_1, \dots, \sigma_n) \coloneqq  (X_1 (x), \dots, X_n (x))$.
Let $x_j \coloneqq  T^j (x)$ and $z_j \coloneqq  \tau^j (\varphi (x))$ for $0 \le j \le n$.
By Lemma \ref{lem:coding}, for $1 \le j \le n$,
\[ x_{j-1} = f_{\sigma_j} (x_j), \ \ \ z_{j-1} = g_{\sigma_j} (z_j), \ \ \ z_j = \varphi (x_j). \]
Since $z_j \in [0,1)$, by Lemma \ref{lem:basic} (2), $z_{j-1} \in [b_{\sigma_j}^g, b_{\sigma_j + 1}^g)$, and hence
\[ \Psi_g (z_{j-1}) = - \log g_{\sigma_j}^{\prime} (z_j), \ \ \ \Psi_f (z_{j-1}) = - \log f_{\sigma_j}^{\prime} (x_j), \ \ \ 1 \le j \le n. \]
Since $x_j = f_{\sigma_{j+1}} \circ \cdots \circ f_{\sigma_n} (x_n)$ and $z_j = g_{\sigma_{j+1}} \circ \cdots \circ g_{\sigma_n} (z_n)$ for $0 \le j \le n$, by the chain rule,
\begin{align*}
\sum_{j=0}^{n-1} \Psi_g (z_j) &= - \sum_{j=1}^{n} \log g_{\sigma_j}^{\prime} (z_j) = - \log g_{\sigma}^{\prime} (z_n), \\
\sum_{j=0}^{n-1} \Psi_f (z_j) &= - \sum_{j=1}^{n} \log f_{\sigma_j}^{\prime} (x_j) = - \log f_{\sigma}^{\prime} (x_n).
\end{align*}
By Lemma \ref{lem:coding} (1) and Lemma \ref{lem:conj} (2), $\mu_{\varphi} (J_n (x)) = g_{\sigma} (1) - g_{\sigma} (0)$ and $\ell (J_n (x)) = f_{\sigma} (1) - f_{\sigma} (0)$.
Hence, by Lemma \ref{lem:hyp} (3) for $h = g$ and for $h = f$,
\begin{align*}
\left| \sum_{j=0}^{n-1} \Psi_g \left( \tau^j (\varphi (x)) \right) + \log \mu_{\varphi} (J_n (x)) \right| &\le \log C_g, \\
\left| \sum_{j=0}^{n-1} \Psi_f \left( \tau^j (\varphi (x)) \right) + \log \ell (J_n (x)) \right| &\le \log C_f,
\end{align*}
and the right hand sides do not depend on $n$ or $x$.
By Lemma \ref{lem:conj} (3), $\mu_{\varphi} \left( \varphi^{-1} (\mathcal{N}) \right) = \ell (\mathcal{N}) = 0$.
The assertion follows from these and \eqref{eq:Birkhoff}.
\end{proof}

\subsection{Expressions for the two limits}

\begin{Prop}\label{prop:average}
It holds that
\[ \lambda_g = \int_{Y^g} S \, d\nu^g, \ \ \ \lambda_f = \int_{Y^g \times Y^f} \Lambda \, d\nu. \]
\end{Prop}

\begin{proof}
For $(y,w) \in Y^g \times Y^f$, let
\[ E(y,w) \coloneqq  \int_0^1 k(y,z) \log k\left( w, \varphi^{-1} (z) \right) dz. \]
By Lemma \ref{lem:hyp} (2) for $h = g$ and for $h = f$, the integrand is bounded and is continuous with respect to $(y,w)$, and hence $E$ is continuous on $Y^g \times Y^f$.
Let $P$ be the operator in the proof of Lemma \ref{lem:stationary}.
We show that, for $(y,w) \in Y^g \times Y^f$,
\[ \int_0^1 k(y,z) \, \Psi_f (z) \, dz = \Lambda (y,w) + E(y,w) - P E (y,w). \]
Let $i \in \mathcal{I}_N$, $x \in [0,1)$ and $z \coloneqq  g_i (x)$.
Then, $\tau (z) = x$ and, by \eqref{eq:conj}, $f_i \left( \varphi^{-1} (x) \right) = \varphi^{-1} (z)$, and hence, by Lemma \ref{lem:id} (3) for $h = f$ at the point $\varphi^{-1} (x)$,
\[ \Psi_f (z) = - \log G_i^f (w) + \log k\left( w, \varphi^{-1} (z) \right) - \log k\left( H_i^f (w), \varphi^{-1} (x) \right). \]
We multiply both sides by $k(y,z)$, integrate them with respect to $z$ over $[b_i^g, b_{i+1}^g)$, and take the sum over $i \in \mathcal{I}_N$.
By Lemma \ref{lem:id} (3) for $h = g$, under the change of variables $z = g_i (x)$, $k(y,z) \, dz = G_i^g (y) \, k(H_i^g (y), x) \, dx$, and hence, by \eqref{eq:k-int}, the integral of $k(y,z)$ with respect to $z$ over $[b_i^g, b_{i+1}^g)$ is $G_i^g (y)$.
Thus the three terms of the right hand side give $\Lambda (y,w)$, $E(y,w)$ and $- P E (y,w)$, respectively.

We integrate the above identity with respect to $\nu (dy \, dw)$.
By Fubini's theorem, the definition of $\rho$ and \eqref{eq:lambda}, the integral of the left hand side is $\lambda_f$, and, by \eqref{eq:stationary} for $u = E$, the integrals of $E$ and of $P E$ are equal to each other.
Thus the second identity holds.
The first identity is shown in the same manner, where $f$, $\varphi$, $w$, $\nu$ and \eqref{eq:stationary} are replaced by $g$, the identity map, $y$, $\nu^g$ and \eqref{eq:stationary-marginal}, respectively, and then $\Lambda (y,w)$ is replaced by $S(y)$.
\end{proof}

\begin{proof}[Proof of Theorem \ref{thm:dim}]
Assertion (1) is Lemma \ref{lem:stationary}.
We show (2).
Let $\nu$ be a Borel probability measure on $Y^g \times Y^f$ satisfying \eqref{eq:stationary} for every continuous function $u$ on $Y^g \times Y^f$, and let $\nu^g$ be its first marginal.
Then, the integral of $S$ with respect to $\nu^g$ is equal to the integral of $S(y)$ with respect to $\nu (dy \, dw)$, and hence the two limits in (2) follow from Propositions \ref{prop:limits} and \ref{prop:average}.
By Lemma \ref{lem:id} (5) and Lemma \ref{lem:hyp} (4), $\Lambda (y,w) \ge - \log G_{\max}^f > 0$ for $(y,w) \in Y^g \times Y^f$, and hence the integral of $\Lambda$ with respect to $\nu$ is positive.
Let
\[ c \coloneqq  \frac{\int_{Y^g \times Y^f} S(y) \, \nu(dy \, dw)}{\int_{Y^g \times Y^f} \Lambda (y,w) \, \nu(dy \, dw)}. \]
Then,
\[ \lim_{n \to \infty} \frac{\log \mu_{\varphi} (J_n (x))}{\log \ell (J_n (x))} = c, \ \ \ \textup{$\mu_{\varphi}$-a.s.} \, x. \]
By this and Lemma \ref{lem:dimension} (1), $\dim_H \mu_{\varphi} \le c$.
By Lemma \ref{lem:dimension} (2), if $K$ is a Borel set such that $\mu_{\varphi} (K) = 1$, then $\dim_H K \ge c$.
Hence, $\dim_H \mu_{\varphi} = c$.
Since $S(y) - \Lambda (y,w) = s_N \left( \mathbf{G}^g (y) | \mathbf{G}^f (w) \right)$, we also have that
\[ c = 1 + \frac{\int_{Y^g \times Y^f} s_N \left( \mathbf{G}^g (y) | \mathbf{G}^f (w) \right) \nu(dy \, dw)}{\int_{Y^g \times Y^f} \Lambda (y,w) \, \nu(dy \, dw)}. \qedhere \]
\end{proof}

\subsection{An invariant measure on the source side}\label{subsec:invariant}

We regard $\mu_{\varphi}$ also as a Borel probability measure on $[0,1)$, since $\mu_{\varphi} (\{1\}) = 0$ by the continuity of $\varphi$.
By Lemma \ref{lem:conj}, the restriction of $\varphi$ to $[0,1)$ is a homeomorphism of $[0,1)$, and the image measure of $\mu_{\varphi}$ under $\varphi$ is $\ell$.
Let $\rho$ and $m$ be as in Subsection \ref{subsec:tau}, and let $\widehat{\mu}$ be the Borel probability measure on $[0,1)$ defined by
\[ \widehat{\mu} (B) \coloneqq  m(\varphi (B)) = \int_B \rho (\varphi (x)) \, \mu_{\varphi} (dx) \]
for Borel subsets $B$ of $[0,1)$.

\begin{Prop}\label{prop:hatmu}
(1) $\widehat{\mu}$ and $\mu_{\varphi}$ are mutually absolutely continuous, and
\[ \frac{d\widehat{\mu}}{d\mu_{\varphi}} = \rho \circ \varphi, \ \ \ C_g^{-1} \le \rho \circ \varphi \le C_g. \]
(2) $\widehat{\mu}$ is $T$-invariant and $T$-ergodic.
\end{Prop}

Assertion (2) follows from Proposition \ref{prop:tau}, since $\varphi \circ T = \tau \circ \varphi$ on $[0,1)$.

\begin{proof}
(1) This follows from the definition of $\widehat{\mu}$ and Lemma \ref{lem:rho} (1).

(2) By Lemma \ref{lem:coding} (2), $\varphi (T(x)) = \tau (\varphi (x))$ for $x \in [0,1)$.
Since the restriction of $\varphi$ to $[0,1)$ is a bijection, $\varphi \left( T^{-1} (B) \right) = \tau^{-1} (\varphi (B))$ for every subset $B$ of $[0,1)$.
Hence, since $m$ is invariant with respect to $\tau$ by Proposition \ref{prop:tau}, for every Borel subset $B$ of $[0,1)$,
\[ \widehat{\mu} \left( T^{-1} (B) \right) = m \left( \tau^{-1} (\varphi (B)) \right) = m(\varphi (B)) = \widehat{\mu} (B). \]
Let $B$ be a Borel subset of $[0,1)$ such that $\widehat{\mu} \left( T^{-1} (B) \triangle B \right) = 0$.
Then, $m \left( \tau^{-1} (\varphi (B)) \triangle \varphi (B) \right) = 0$, and, since $m$ is ergodic with respect to $\tau$ by Proposition \ref{prop:tau}, $\widehat{\mu} (B) = m(\varphi (B)) \in \{0,1\}$.
\end{proof}

\begin{Rem}[Entropy and Lyapunov exponent]\label{rem:entropy}
Let
\[ \omega_g (x) \coloneqq  \log g_{X_1 (x)}^{\prime} (\varphi (T(x))), \ \ \ x \in [0,1). \]
By the definition of $\widehat{\mu}$, $m$ is the image measure of $\widehat{\mu}$ under $\varphi$.
Hence, by \eqref{eq:lambda}, \eqref{eq:Psi-omega} and the definition of $\omega_g$,
\[ \lambda_g = - \int_{[0,1)} \omega_g \, d\widehat{\mu}, \ \ \ \lambda_f = - \int_{[0,1)} \omega_f \, d\widehat{\mu}, \]
and, by Propositions \ref{prop:limits} and \ref{prop:average} and Theorem \ref{thm:dim} (2),
\[ \dim_H \mu_{\varphi} = \frac{\int_{[0,1)} \omega_g \, d\widehat{\mu}}{\int_{[0,1)} \omega_f \, d\widehat{\mu}}. \]
Since $\mu_{\varphi}$ has no atoms, by Proposition \ref{prop:hatmu} (1), the set of the points $b_i^f$, $i \in \mathcal{I}_N$, is a $\widehat{\mu}$-null set.
Since $- \omega_f$ is the logarithm of the derivative of $T$ on each open interval $(b_i^f, b_{i+1}^f)$, $i \in \mathcal{I}_N$, $\lambda_f$ is the Lyapunov exponent of $T$ with respect to $\widehat{\mu}$.
By Proposition \ref{prop:hatmu} (1) and Proposition \ref{prop:limits},
\[ \lim_{n \to \infty} \frac{- \log \widehat{\mu} (J_n (x))}{n} = \lambda_g, \ \ \ \textup{$\widehat{\mu}$-a.s.} \, x. \]
Let $\mathcal{P}$ be the partition of $[0,1)$ into the intervals $[b_i^f, b_{i+1}^f)$, $i \in \mathcal{I}_N$.
By Lemma \ref{lem:coding} (1), for every $n \ge 1$, the common refinement of the partitions $\{ T^{-j} (A) : A \in \mathcal{P} \}$, $0 \le j \le n-1$, consists of the intervals $J(\sigma)$, $\sigma \in (\mathcal{I}_N)^n$, whose lengths are at most $\left( G_{\max}^f \right)^n$ by \eqref{eq:length}.
Hence, $\mathcal{P}$ is a generator, and, by Proposition \ref{prop:hatmu} (2), the Shannon--McMillan--Breiman theorem and the Kolmogorov--Sinai theorem, $\lambda_g$ is the entropy of $T$ with respect to $\widehat{\mu}$.
Thus the formula in Theorem \ref{thm:dim} (2) states that $\dim_H \mu_{\varphi}$ is the entropy of $T$ with respect to $\widehat{\mu}$ divided by its Lyapunov exponent.
\end{Rem}

\section{Proof of Theorem \ref{thm:sing}}\label{sec:sing}

We first express the assumption of Theorem \ref{thm:sing} (1) in terms of $a_i^h, b_i^h, c_i^h$ and in terms of $G_i^h, H_i^h$, $h \in \{f, g\}$.

\begin{Prop}\label{prop:equiv}
Let $f$ and $g$ be LF systems on $[0,1]$, and let $s > -1$.
Then, the following are equivalent:\\
(i) $g_i \circ \psi_s = \psi_s \circ f_i$ on $[0,1]$ for every $i \in \mathcal{I}_N$.\\
(ii) $G_i^g (y) = G_i^f (\gamma_s (y))$ and $\gamma_s \left( H_i^g (y) \right) = H_i^f (\gamma_s (y))$ for every $i \in \mathcal{I}_N$ and every $y > -1$.\\
(iii) There exists $y > -1$ such that $G_i^g (y) = G_i^f (\gamma_s (y))$ and $\gamma_s \left( H_i^g (y) \right) = H_i^f (\gamma_s (y))$ for every $i \in \mathcal{I}_N$.\\
(iv) $b_i^g = \psi_s (b_i^f)$ and $(1+s)(1 + c_i^g) = 1 + H_i^f (s)$ for $i \in \mathcal{I}_N$.\\
(v) For every $i \in \mathcal{I}_N$,
\[ \frac{b_i^g}{1 - b_i^g} = (1+s) \, \frac{b_i^f}{1 - b_i^f}, \ \ \ \frac{a_i^g + b_i^g}{1 - b_i^g} = \frac{a_i^f + b_i^f}{1 - b_i^f}. \]
\end{Prop}

\begin{proof}
Let $y > -1$ and $w \coloneqq  \gamma_s (y)$.
By Lemma \ref{lem:psi} (2) and Lemma \ref{lem:id} (4) for $h = g$ and for $h = f$, for each $i \in \mathcal{I}_N$, the compositions of $\psi_y$ with the two sides of (i) are
\begin{align*}
\psi_y \circ g_i \circ \psi_s &= \psi_y (b_i^g) + G_i^g (y) \, \psi_{\gamma_s (H_i^g (y))}, \\
\psi_y \circ \psi_s \circ f_i &= \psi_w (b_i^f) + G_i^f (w) \, \psi_{H_i^f (w)}.
\end{align*}
Since $G_i^g (y) > 0$, by Lemma \ref{lem:psi} (1) and (4) and the values of the right hand sides at $0$, $1$ and $1/2$, they are equal to each other on $[0,1]$ if and only if $\psi_y (b_i^g) = \psi_w (b_i^f)$, $G_i^g (y) = G_i^f (w)$ and $\gamma_s \left( H_i^g (y) \right) = H_i^f (w)$.
By Lemma \ref{lem:id} (5), $\psi_y (b_i^g)$ and $\psi_w (b_i^f)$ are the sums of $G_j^g (y)$ and of $G_j^f (w)$ over $0 \le j \le i-1$, respectively.
Hence, since $\psi_y$ is injective, (i) holds if and only if $G_i^g (y) = G_i^f (w)$ and $\gamma_s \left( H_i^g (y) \right) = H_i^f (w)$ for every $i \in \mathcal{I}_N$.
Since $y > -1$ is arbitrary, (i), (ii) and (iii) are equivalent.

Let $y = 0$.
Then, $w = s$, $H_i^g (0) = c_i^g$ and $\gamma_s (c_i^g) + 1 = (1+s)(1 + c_i^g)$.
By Lemma \ref{lem:id} (5), since $b_0^h = 0$ and $b_N^h = 1$ for $h \in \{f, g\}$, $G_i^g (0) = G_i^f (s)$ for every $i \in \mathcal{I}_N$ if and only if $b_i^g = \psi_s (b_i^f)$ for every $i \in \mathcal{I}_N$.
Hence, by the above for $y = 0$, (i) and (iv) are equivalent.

We show that (iv) and (v) are equivalent.
By Lemma \ref{lem:psi} (3), for $1 \le i \le N-1$, $b_i^g = \psi_s (b_i^f)$ if and only if the first equation in (v) holds, and both of them hold for $i = 0$.
Assume that $b_i^g = \psi_s (b_i^f)$ for every $i \in \mathcal{I}_N$, which also holds for $i = N$.
Then, for $i \in \mathcal{I}_N$, since $a_i^h + b_i^h = b_{i+1}^h (1 + c_i^h) > 0$ for $h \in \{f, g\}$, by the definition of $\psi_s$ and \eqref{eq:one-plus-H} for $h = f$ and $y = s$,
\[ \frac{(a_i^g + b_i^g)(1 - b_i^f)}{(a_i^f + b_i^f)(1 - b_i^g)} = \frac{(1+s)(1 + c_i^g)(1 + s b_i^f)}{(1 + c_i^f)(1 + s b_{i+1}^f)} = \frac{(1+s)(1 + c_i^g)}{1 + H_i^f (s)}. \]
Hence, the second equation in (v) holds if and only if the second condition in (iv) holds.
\end{proof}

Let $f$, $g$ and $s$ be as in Proposition \ref{prop:equiv}.
By \eqref{eq:one-plus-H}, the second condition in (iv) is that
\[ (1+s) \left( 1 + s b_i^f \right) \left( 1 + c_i^g \right) = \left( 1 + s b_{i+1}^f \right) \left( 1 + c_i^f \right), \ \ \ i \in \mathcal{I}_N. \]
If (i) holds, then, by (iv), $\psi_s (b_1^f) = b_1^g$, and hence $s = s_0$.
Hence, the assumption of Theorem \ref{thm:sing} (1) holds if and only if (iv) holds for $s = s_0$, and also if and only if (v) holds for $s = s_0$.

We now show Theorem \ref{thm:sing}.
In the proof of (2), we show that if $\dim_H \mu_{\varphi} = 1$, then (iii) in Proposition \ref{prop:equiv} holds for $s = s_0$.
By Theorem \ref{thm:dim} (2), the integral with respect to $\nu$ of the relative entropy of $\mathbf{G}^g (y)$ with respect to $\mathbf{G}^f (w)$ is $0$, and, by \eqref{eq:stationary}, so is the integral of the relative entropy at the points $\left( H_i^g (y), H_i^f (w) \right)$ for each $i \in \mathcal{I}_N$.
A point at which all of these relative entropies are $0$ satisfies the conditions in (iii).

\begin{proof}[Proof of Theorem \ref{thm:sing}]
(1) Assume that $g_i \circ \psi_{s_0} = \psi_{s_0} \circ f_i$ on $[0,1]$ for every $i \in \mathcal{I}_N$.
By Lemma \ref{lem:psi} (1), $\psi_{s_0}$ is a continuous and strictly increasing function from $[0,1]$ to $[0,1]$, and it satisfies \eqref{eq:conj}.
Hence, by the uniqueness of the solution of \eqref{eq:conj}, $\varphi = \psi_{s_0}$.
Since $\psi_{s_0}$ is continuously differentiable on $[0,1]$, $\mu_{\varphi}$ is absolutely continuous with respect to the Lebesgue measure.

(2) If $\dim_H \mu_{\varphi} < 1$, then there exists a Borel set $K$ such that $\mu_{\varphi} (K) = 1$ and $\dim_H K < 1$, and then $\ell (K) = 0$.
Hence, it suffices to show that if $\dim_H \mu_{\varphi} = 1$, then the assumption of (1) holds.

Assume that $\dim_H \mu_{\varphi} = 1$.
By Theorem \ref{thm:dim} (1), there exists a Borel probability measure $\nu$ on $Y^g \times Y^f$ satisfying \eqref{eq:stationary} for every continuous function $u$ on $Y^g \times Y^f$.
Let
\[ D(y,w) \coloneqq  - s_N \left( \mathbf{G}^g (y) | \mathbf{G}^f (w) \right), \ \ \ (y,w) \in Y^g \times Y^f. \]
Then, $D$ is nonnegative and continuous on $Y^g \times Y^f$, and $D(y,w) = 0$ if and only if $G_j^g (y) = G_j^f (w)$ for every $j \in \mathcal{I}_N$.
By Theorem \ref{thm:dim} (2) and \eqref{eq:stationary} for $u = D$,
\[ \int_{Y^g \times Y^f} D \, d\nu = 0, \ \ \ \int_{Y^g \times Y^f} \sum_{i \in \mathcal{I}_N} G_i^g (y) \, D\left( H_i^g (y), H_i^f (w) \right) \nu(dy \, dw) = 0. \]
Since $D \ge 0$ and $G_i^g (y) > 0$ for $y \in Y^g$ and $i \in \mathcal{I}_N$, each of the $N+1$ functions $(y,w) \mapsto D(y,w)$ and $(y,w) \mapsto D\left( H_i^g (y), H_i^f (w) \right)$, $i \in \mathcal{I}_N$, is $0$, $\nu$-a.e.
Since these functions are finitely many, there exists a point $(y,w) \in Y^g \times Y^f$ at which all of them are $0$, that is,
\[ D(y,w) = 0, \ \ \ D\left( H_i^g (y), H_i^f (w) \right) = 0, \ \ i \in \mathcal{I}_N. \]

We show that $w^{\prime} = \gamma_{s_0} (y^{\prime})$ for every point $(y^{\prime}, w^{\prime}) \in Y^g \times Y^f$ such that $D(y^{\prime}, w^{\prime}) = 0$, which we apply to the $N+1$ points above.
Let $(y^{\prime}, w^{\prime}) \in Y^g \times Y^f$ be a point such that $D(y^{\prime}, w^{\prime}) = 0$, and let $s > -1$ be the number such that $1 + w^{\prime} = (1+s)(1 + y^{\prime})$, that is, $w^{\prime} = \gamma_s (y^{\prime})$.
By Lemma \ref{lem:id} (5) and Lemma \ref{lem:psi} (2),
\[ \psi_{y^{\prime}} (b_1^g) = G_0^g (y^{\prime}) = G_0^f (w^{\prime}) = \psi_{w^{\prime}} (b_1^f) = \psi_{y^{\prime}} (\psi_s (b_1^f)). \]
Since $\psi_{y^{\prime}}$ is injective, $b_1^g = \psi_s (b_1^f)$, and hence, by Lemma \ref{lem:psi} (4) and the definition of $s_0$, $s = s_0$.
Thus $w^{\prime} = \gamma_{s_0} (y^{\prime})$.

By this for the points $(y,w)$ and $\left( H_i^g (y), H_i^f (w) \right)$, $i \in \mathcal{I}_N$, we have that $w = \gamma_{s_0} (y)$, and
\[ G_i^g (y) = G_i^f (\gamma_{s_0} (y)), \ \ \ H_i^f (\gamma_{s_0} (y)) = \gamma_{s_0} \left( H_i^g (y) \right), \ \ \ i \in \mathcal{I}_N. \]
Thus (iii) in Proposition \ref{prop:equiv} holds for $s = s_0$, and, by Proposition \ref{prop:equiv}, the assumption of (1) holds.
\end{proof}

\begin{Rem}[The conditions in Proposition \ref{prop:equiv} (v)]\label{rem:odds}
Let $f$ and $g$ be LF systems on $[0,1]$, and let $s > -1$.\\
(1) For $x \in [0,1)$, let $\widetilde{x} \coloneqq  x/(1-x)$, and, for an LF system $h$ and $i \in \mathcal{I}_N$, let
\[ \widetilde{A}_i^h \coloneqq  \begin{pmatrix} a_i^h + b_i^h & b_i^h \\ 1 + c_i^h - a_i^h - b_i^h & 1 - b_i^h \end{pmatrix}. \]
Since $a_i^h x + b_i^h = (a_i^h + b_i^h) x + b_i^h (1-x)$ and $(c_i^h x + 1) - (a_i^h x + b_i^h) = (1 + c_i^h - a_i^h - b_i^h) x + (1 - b_i^h)(1-x)$, we have that
\[ \frac{h_i (x)}{1 - h_i (x)} = \frac{(a_i^h + b_i^h) \, \widetilde{x} + b_i^h}{(1 + c_i^h - a_i^h - b_i^h) \, \widetilde{x} + 1 - b_i^h}, \ \ \ x \in [0,1), \]
that is, in the coordinate $\widetilde{x}$, $h_i$ is the linear fractional transformation given by $\widetilde{A}_i^h$.
By the condition (3) of LF systems, the first column of $\widetilde{A}_i^h$ is ${}^t \left( b_{i+1}^h (1 + c_i^h), (1 - b_{i+1}^h)(1 + c_i^h) \right)$, and hence the entries of $\widetilde{A}_i^h$ are nonnegative.
The determinant of $\widetilde{A}_i^h$ is $a_i^h - b_i^h c_i^h > 0$.
By Lemma \ref{lem:psi} (3), $\psi_s$ multiplies $\widetilde{x}$ by $1+s$.
Two $2 \times 2$ real matrices with non-zero determinants give the same linear fractional transformation on an interval if and only if they are constant multiples of each other.
Hence, (i) in Proposition \ref{prop:equiv} holds if and only if, for every $i \in \mathcal{I}_N$, $\widetilde{A}_i^g$ is a constant multiple of the matrix obtained from $\widetilde{A}_i^f$ by multiplying its $(1,2)$ entry by $1+s$ and dividing its $(2,1)$ entry by $1+s$, which is the conjugate of $\widetilde{A}_i^f$ by the diagonal matrix with diagonal entries $1+s$ and $1$.
The two equations in (v) are the identities for the ratios of the $(1,2)$ entries and of the $(1,1)$ entries to the $(2,2)$ entries.
The identity for the ratios of the $(2,1)$ entries to the $(2,2)$ entries follows from the equations in (v) and the condition (3) of LF systems.
Furthermore, by \eqref{eq:one-plus-H}, for $y > -1$, $1 + H_i^h (y)$ is the value at $1 + y$ of the linear fractional transformation given by the transposed matrix of $\widetilde{A}_i^h$, and $\gamma_s$ multiplies $1 + y$ by $1 + s$.
Thus the second identity in (ii), for every $y > -1$, is the transposed form of the same condition.\\
(2) The conditions in (v) are symmetric with respect to $f$ and $g$.
Indeed, (v) holds if and only if (v) holds for $g$, $f$ and $s^{\prime}$ in place of $f$, $g$ and $s$, respectively, where $s^{\prime} > -1$ is the number such that $(1+s)(1+s^{\prime}) = 1$, for which $\psi_{s^{\prime}}$ is the inverse function of $\psi_s$.
The second equation in (v) does not involve $s$.\\
(3) Since $a_{N-1}^h + b_{N-1}^h = 1 + c_{N-1}^h$ for $h \in \{f, g\}$, by Lemma \ref{lem:basic} (4), the second equation in (v) for $i = 0$ is that $g_0^{\prime} (0) = f_0^{\prime} (0)$, and that for $i = N-1$ is that $g_{N-1}^{\prime} (1) = f_{N-1}^{\prime} (1)$.
Hence, if $f_0^{\prime} (0) \ne g_0^{\prime} (0)$ or $f_{N-1}^{\prime} (1) \ne g_{N-1}^{\prime} (1)$, then the assumption of Theorem \ref{thm:sing} (1) does not hold.
In particular, it does not hold if one of $f$ and $g$ is hyperbolic and the other is parabolic.

Let $N = 2$.
Then, the first equation in (v) for $s = s_0$ holds for $i = 0$, and it holds for $i = 1$ by the definition of $s_0$ and Lemma \ref{lem:psi} (3).
Hence, the assumption of Theorem \ref{thm:sing} (1) holds if and only if $f_0^{\prime} (0) = g_0^{\prime} (0)$ and $f_1^{\prime} (1) = g_1^{\prime} (1)$.
If $f_0 (x) = x/2$ and $f_1 (x) = (x+1)/2$, then $f_0^{\prime} (0) = f_1^{\prime} (1) = 1/2$, and the conditions that $g_0^{\prime} (0) = 1/2$ and that $g_1^{\prime} (1) = 1/2$ are the conditions (i) and (ii) in \cite[Theorem 1.2]{O14} for the matrices $A_0^g$ and $A_1^g$, respectively.
Indeed, since $a_1^g + b_1^g = 1 + c_1^g$, the conditions (i) and (ii) in \cite[Theorem 1.2]{O14} are equivalent to
\[ (1 - 2 a_0^g)(1 + c_0^g - a_0^g) = 0, \ \ \ (1 - b_1^g)(1 - 2 b_1^g - c_1^g) = 0, \]
respectively, and we have that $1 + c_0^g - a_0^g = (1 + c_0^g)(1 - b_1^g) > 0$, $1 - b_1^g > 0$, $g_0^{\prime} (0) = a_0^g$ and $g_1^{\prime} (1) = (1 - b_1^g)/(1 + c_1^g)$.
\end{Rem}

\begin{Exa}\label{exa:two}
Let $N = 2$, and let
\[ f_0 (x) \coloneqq  \frac{3x}{2x+4}, \ \ \ f_1 (x) \coloneqq  \frac{2x+1}{x+2}, \ \ \ x \in [0,1]. \]
Then, $f = \{f_0, f_1\}$ is a hyperbolic LF system on $[0,1]$ for which $(a_0^f, b_0^f, c_0^f) = (3/4, 0, 1/2)$ and $(a_1^f, b_1^f, c_1^f) = (1, 1/2, 1/2)$, and, by Lemma \ref{lem:basic} (4), $f_0^{\prime} (0) = 3/4$ and $f_1^{\prime} (1) = 1/3$.\\
(1) Let
\[ g_0 (z) \coloneqq  \frac{3z}{z+4}, \ \ \ g_1 (z) \coloneqq  \frac{3z+3}{z+5}, \ \ \ z \in [0,1]. \]
Then, $g = \{g_0, g_1\}$ is a hyperbolic LF system on $[0,1]$ for which $(a_0^g, b_0^g, c_0^g) = (3/4, 0, 1/4)$ and $(a_1^g, b_1^g, c_1^g) = (3/5, 3/5, 1/5)$, and $g_0^{\prime} (0) = 3/4$ and $g_1^{\prime} (1) = 1/3$.
Hence, by Remark \ref{rem:odds} (3), the assumption of Theorem \ref{thm:sing} (1) holds.
Since $f_0 (1) = 1/2$ and $g_0 (1) = 3/5$, we have that $1 + s_0 = 3/2$, and, by Theorem \ref{thm:sing} (1),
\[ \varphi (x) = \psi_{1/2} (x) = \frac{3x}{x+2}, \ \ \ x \in [0,1]. \]
Indeed, $g_0 (\varphi (x)) = 9x/(7x+8) = \varphi (f_0 (x))$ and $g_1 (\varphi (x)) = (6x+3)/(4x+5) = \varphi (f_1 (x))$ for $x \in [0,1]$.\\
(2) Let
\[ g_0 (z) \coloneqq  \frac{3z}{z+4}, \ \ \ g_1 (z) \coloneqq  \frac{7z+3}{5z+5}, \ \ \ z \in [0,1]. \]
Then, $g = \{g_0, g_1\}$ is a hyperbolic LF system on $[0,1]$ for which $(a_0^g, b_0^g, c_0^g) = (3/4, 0, 1/4)$ and $(a_1^g, b_1^g, c_1^g) = (7/5, 3/5, 1)$, and $g_0^{\prime} (0) = 3/4$ and $g_1^{\prime} (1) = 1/5$.
Since $g_1^{\prime} (1) \ne f_1^{\prime} (1)$, by Remark \ref{rem:odds} (3) and Theorem \ref{thm:sing} (2), $\dim_H \mu_{\varphi} < 1$, and $\mu_{\varphi}$ is singular with respect to the Lebesgue measure.

We compute $\dim_H \mu_{\varphi}$ numerically by Theorem \ref{thm:dim} (2).
For $n \ge 0$, we consider the ratio
\[ \frac{\sum_{\sigma \in (\mathcal{I}_2)^n} \mu_{\varphi} (J(\sigma)) \, S(y_{\sigma}^g)}{\sum_{\sigma \in (\mathcal{I}_2)^n} \mu_{\varphi} (J(\sigma)) \, \Lambda (y_{\sigma}^g, y_{\sigma}^f)}, \]
which is the right hand side of the formula for $\dim_H \mu_{\varphi}$ in Theorem \ref{thm:dim} (2) in which $\nu$ is replaced by the distribution of the points $(y_{\sigma}^g, y_{\sigma}^f)$, where the words $\sigma$ of length $n$ are chosen with the probabilities $\mu_{\varphi} (J(\sigma))$.
These sums are computed by \eqref{eq:next-letter} and Lemma \ref{lem:id} (1).
Thus this computation uses only the transition kernel in \eqref{eq:stationary}, and it does not use a density of $\nu$; in this example, $\nu$ is concentrated on a curve (Remark \ref{rem:curve} (2)), and hence it is singular with respect to the Lebesgue measure on $Y^g \times Y^f$.
The values of this ratio for $n = 16$, $20$ and $24$ are approximately $0.976734$, $0.976754$ and $0.976758$, respectively, which suggests that $\dim_H \mu_{\varphi}$ is approximately $0.9768$.
The numerator for $n = 24$ is approximately $0.531748$, and, by Remark \ref{rem:two} (1), the first limit in Theorem \ref{thm:dim} (2) is $2 \left( L(3/4) + L(1/5) - \pi^2 / 6 \right) / \log (3/4)$, which is approximately $0.531747$.
By the proof of Lemma \ref{lem:stationary} and the uniqueness of $\nu$ (Remark \ref{rem:Blackwell} (2)), the averages of these distributions over the lengths $0, 1, \dots, n-1$ converge weakly to $\nu$ as $n \to \infty$.
We do not show the convergence of the above ratio in this paper, and we do not estimate the errors of these values.

For each of the three LF systems $h$ in (1) and (2), $h_0^{\prime} (0) + h_1^{\prime} (1) \ne 1$, and hence, by Remark \ref{rem:two} in which $f$ is replaced by $h$, there are no $s > -1$ such that $\psi_s \circ h_i \circ \psi_s^{-1}$, $i \in \{0, 1\}$, are affine maps.
\end{Exa}

\section{The parabolic case}\label{sec:para}

In this section, we show Theorem \ref{thm:para} and Corollary \ref{cor:positive}.
{\it Throughout this section, except for the proof of Corollary \ref{cor:positive}, let $f$ and $g$ be LF systems on $[0,1]$ such that $g$ is parabolic, and let $\varphi$ be the solution of \eqref{eq:conj}.}
Thus $g_0^{\prime} (0) = 1$ or $g_{N-1}^{\prime} (1) = 1$.
We say that a letter $i \in \mathcal{I}_N$ is parabolic if $i = 0$ and $g_0^{\prime} (0) = 1$, or $i = N-1$ and $g_{N-1}^{\prime} (1) = 1$.

Let
\[ \theta_0 \coloneqq  \frac{1 - b_1^g}{b_1^g}, \ \ \ \theta_1 \coloneqq  \frac{b_{N-1}^g}{1 - b_{N-1}^g}. \]
By Lemma \ref{lem:basic} (2), $\theta_0 > 0$ and $\theta_1 > 0$.
For words $\sigma = (\sigma_1, \dots, \sigma_n)$ and $\sigma^{\prime} = (\sigma_1^{\prime}, \dots, \sigma_{n^{\prime}}^{\prime})$, we denote the word $(\sigma_1, \dots, \sigma_n, \sigma_1^{\prime}, \dots, \sigma_{n^{\prime}}^{\prime})$ by $\sigma \sigma^{\prime}$.
For $i \in \mathcal{I}_N$ and $m \ge 0$, we denote by $i^m$ the word of length $m$ all of whose entries are $i$.
Thus, for a word $\sigma$, $\sigma i^m$ is the word obtained by adding $m$ letters $i$ to the end of $\sigma$.

\subsection{Repetitions of parabolic letters}

\begin{Lem}\label{lem:run}
Let $\sigma$ be a word and $m \ge 1$. Then,\\
(1) If $g_0^{\prime} (0) = 1$, then $H_0^g (y) = y + \theta_0$ for $y > -1$, and
\[ 1 + y_{\sigma 0^m}^g = 1 + y_{\sigma}^g + m \theta_0, \ \ \ \mu_{\varphi} (J(\sigma 0^m)) = \frac{1 + y_{\sigma}^g}{1 + y_{\sigma}^g + m \theta_0} \, \mu_{\varphi} (J(\sigma)). \]
(2) If $g_{N-1}^{\prime} (1) = 1$, then $\left( 1 + H_{N-1}^g (y) \right)^{-1} = (1+y)^{-1} + \theta_1$ for $y > -1$, and
\begin{align*}
\frac{1}{1 + y_{\sigma (N-1)^m}^g} &= \frac{1}{1 + y_{\sigma}^g} + m \theta_1, \\
\mu_{\varphi} (J(\sigma (N-1)^m)) &= \frac{1}{1 + m \theta_1 (1 + y_{\sigma}^g)} \, \mu_{\varphi} (J(\sigma)).
\end{align*}
(3) There exists a constant $C_1 \ge 1$ depending only on $g$ such that
\[ \mu_{\varphi} (J(\sigma i^m)) \ge \frac{\mu_{\varphi} (J(\sigma))}{1 + C_1 m} \]
for some parabolic letter $i$.
\end{Lem}

In the proof of (1) and (2), we use that $G_0^g (y) = (1+y) / (1 + H_0^g (y))$ if $g_0^{\prime} (0) = 1$, and that $G_{N-1}^g (y) = (1 + H_{N-1}^g (y)) / (1+y)$ if $g_{N-1}^{\prime} (1) = 1$, by which the products of these values along the repetitions of a parabolic letter are telescoping.
Assertion (3) follows from (1) and (2), where we use a lower bound or an upper bound for $y_{\sigma}^g$ if only one of the letters $0$ and $N-1$ is parabolic.

\begin{proof}
We omit the superscript $g$.
By Lemma \ref{lem:conj} (2) and Lemma \ref{lem:id} (6), for every word $\sigma^{\prime}$ and every $i \in \mathcal{I}_N$,
\begin{equation}\label{eq:mass}
\mu_{\varphi} (J(\sigma^{\prime} i)) = G_i (y_{\sigma^{\prime}}) \, \mu_{\varphi} (J(\sigma^{\prime})).
\end{equation}

(1) By Lemma \ref{lem:basic} (4), $a_0 = 1$.
Since $b_0 = 0$ and $g_0 (1) = b_1$, we have that $c_0 + 1 = 1/b_1$, that is, $c_0 = \theta_0$.
Hence, $H_0 (y) = a_0 y + c_0 = y + \theta_0$.
Since $1/b_1 = 1 + \theta_0$, by the definition of $G_0$,
\[ G_0 (y) = \frac{b_1 (1+y)}{1 + b_1 y} = \frac{1+y}{1 + y + \theta_0} = \frac{1+y}{1 + H_0 (y)}. \]
By Lemma \ref{lem:id} (1), $y_{\sigma 0^{j+1}} = H_0 (y_{\sigma 0^j}) = y_{\sigma 0^j} + \theta_0$ for $j \ge 0$.
By these and \eqref{eq:mass},
\[ \frac{\mu_{\varphi} (J(\sigma 0^m))}{\mu_{\varphi} (J(\sigma))} = \prod_{j=0}^{m-1} G_0 (y_{\sigma 0^j}) = \prod_{j=0}^{m-1} \frac{1 + y_{\sigma 0^j}}{1 + y_{\sigma 0^{j+1}}} = \frac{1 + y_{\sigma}}{1 + y_{\sigma} + m \theta_0}. \]

(2) By Lemma \ref{lem:basic} (4), $c_{N-1} + 1 = 1 - b_{N-1}$.
By this and \eqref{eq:one-plus-H} for $i = N-1$,
\[ 1 + H_{N-1} (y) = \frac{(1 - b_{N-1})(1+y)}{1 + b_{N-1} y}, \]
and hence, by the definition of $G_{N-1}$,
\[ G_{N-1} (y) = \frac{1 - b_{N-1}}{1 + b_{N-1} y} = \frac{1 + H_{N-1} (y)}{1+y}. \]
Since $1 + b_{N-1} y = (1 - b_{N-1}) + b_{N-1} (1+y)$, we have that
\[ \frac{1}{1 + H_{N-1} (y)} = \frac{1}{1+y} + \theta_1. \]
By this and Lemma \ref{lem:id} (1), $\left( 1 + y_{\sigma (N-1)^j} \right)^{-1} = (1 + y_{\sigma})^{-1} + j \theta_1$ for $j \ge 0$.
By these and \eqref{eq:mass},
\begin{align*}
\frac{\mu_{\varphi} (J(\sigma (N-1)^m))}{\mu_{\varphi} (J(\sigma))} &= \prod_{j=0}^{m-1} G_{N-1} (y_{\sigma (N-1)^j}) = \prod_{j=0}^{m-1} \frac{1 + y_{\sigma (N-1)^{j+1}}}{1 + y_{\sigma (N-1)^j}} \\
&= \frac{1 + y_{\sigma (N-1)^m}}{1 + y_{\sigma}} = \frac{1}{1 + m \theta_1 (1 + y_{\sigma})}.
\end{align*}

(3) Assume that $g_0^{\prime} (0) = g_{N-1}^{\prime} (1) = 1$.
If $y_{\sigma} \ge 0$, then $(1 + y_{\sigma}) / (1 + y_{\sigma} + m \theta_0) \ge 1 / (1 + m \theta_0)$.
If $y_{\sigma} < 0$, then $1 / (1 + m \theta_1 (1 + y_{\sigma})) \ge 1 / (1 + m \theta_1)$.
Hence, by (1) and (2), the assertion holds for $C_1 = \max\{1, \theta_0, \theta_1\}$.

Assume that $g_0^{\prime} (0) = 1$ and $g_{N-1}^{\prime} (1) < 1$.
By Lemma \ref{lem:fixed}, for each $1 \le i \le N-1$, $H_i$ has a unique fixed point $e_i$ in $(-1, +\infty)$.
Let $\alpha_{-}$ be the minimum of $0, e_1, \dots, e_{N-1}$.
Then, $\alpha_{-} > -1$.
Let $y \ge \alpha_{-}$.
By (1), $H_0 (y) > y \ge \alpha_{-}$.
Since $H_i$ is increasing, by Lemma \ref{lem:fixed}, $H_i (y) \ge H_i (\alpha_{-}) \ge \alpha_{-}$ for $1 \le i \le N-1$.
Since $y_{\sigma^{\prime}} = 0 \ge \alpha_{-}$ for the empty word $\sigma^{\prime}$, by Lemma \ref{lem:id} (1), $y_{\sigma^{\prime}} \ge \alpha_{-}$ for every word $\sigma^{\prime}$.
Hence, by (1), the assertion holds for $C_1 = \max\{1, \theta_0 / (1 + \alpha_{-})\}$.

Assume that $g_0^{\prime} (0) < 1$ and $g_{N-1}^{\prime} (1) = 1$.
By Lemma \ref{lem:fixed}, for each $0 \le i \le N-2$, $H_i$ has a unique fixed point $e_i$ in $(-1, +\infty)$.
Let $\beta_{+}$ be the maximum of $0, e_0, \dots, e_{N-2}$.
Let $-1 < y \le \beta_{+}$.
By (2), $H_{N-1} (y) < y \le \beta_{+}$.
Since $H_i$ is increasing, by Lemma \ref{lem:fixed}, $H_i (y) \le H_i (\beta_{+}) \le \beta_{+}$ for $0 \le i \le N-2$.
Hence, in the same manner as above, $y_{\sigma^{\prime}} \le \beta_{+}$ for every word $\sigma^{\prime}$, and, by (2), the assertion holds for $C_1 = \max\{1, \theta_1 (1 + \beta_{+})\}$.
\end{proof}

By Lemma \ref{lem:run} (1) and (2), the dual maps $H_i^g$, $i \in \mathcal{I}_N$, of the parabolic LF system $g$ have no common fixed points in $(-1, +\infty)$.

\subsection{Proofs of Theorem \ref{thm:para} and Corollary \ref{cor:positive}}

For $n \ge 1$, let $\mathcal{F}_n$ be the $\sigma$-algebra on $[0,1)$ generated by $X_1, \dots, X_n$, and let $\mathcal{F}_0$ be the trivial $\sigma$-algebra on $[0,1)$.
By Lemma \ref{lem:coding} (1), for every $n \ge 0$, $\mathcal{F}_n$ is generated by the intervals $J(\sigma)$, $\sigma \in (\mathcal{I}_N)^n$.

\begin{Prop}\label{prop:parabolic}
It holds that
\[ \lim_{n \to \infty} \frac{- \log \mu_{\varphi} (J_n (x))}{n} = 0, \ \ \ \textup{$\mu_{\varphi}$-a.s.} \, x. \]
\end{Prop}

In the proof, we decompose the sequence $(X_n (x))_{n \ge 1}$ into blocks, in each of which a parabolic letter is repeated, and the letters separating them.
In Step 1, we show by Lemma \ref{lem:run} (3) that the number of the separating letters among the first $n$ letters divided by $n$ converges to $0$ as $n \to \infty$.
In Step 2, we count the words of length $n$ in which few pairs of adjacent letters differ.

\begin{proof}
For $x \in [0,1)$, let $t_0 (x) \coloneqq  0$.
For $j \ge 0$ with $t_j (x) < +\infty$, let $\xi_j (x)$ be the supremum of the integers $m \ge 0$ such that
\[ X_{t_j (x) + 1} (x) = \cdots = X_{t_j (x) + m} (x) = i \ \textup{ for some parabolic letter } i, \]
and let $t_{j+1} (x) \coloneqq  t_j (x) + \xi_j (x) + 1$.
For $j \ge 0$ with $t_j (x) = +\infty$, let $\xi_j (x) \coloneqq  0$ and $t_{j+1} (x) \coloneqq  +\infty$.
Thus, if $t_{j+1} (x) < +\infty$, then $\xi_j (x)$ is the length of the $j$-th block, which may be zero, and $t_{j+1} (x)$ is the position of the letter separating the $j$-th block and the $(j+1)$-th block.

Let $j \ge 0$, and assume that $t_j (x) < +\infty$ and $\xi_j (x) = +\infty$.
Then, by Lemma \ref{lem:coding} (1), $x \in J(\sigma i^m)$ for every $m \ge 1$, where $\sigma$ is the word $(X_1 (x), \dots, X_{t_j (x)} (x))$ and $i$ is a parabolic letter.
By Lemma \ref{lem:run} (1) and (2), for every word $\sigma$ and every parabolic letter $i$, $\mu_{\varphi} (J(\sigma i^m))$ converges to $0$ as $m \to \infty$.
Since the words are countably many, we have that $t_j < +\infty$ for every $j \ge 0$, $\mu_{\varphi}$-a.s.

For $v \ge 0$ and $m \ge 1$, let
\[ \mathcal{E}_{v,m} \coloneqq  \left\{ X_{v+1} = \cdots = X_{v+m} = i \ \textup{ for some parabolic letter } i \right\}. \]
Then, $\mathcal{E}_{v,m} \in \mathcal{F}_{v+m}$, and, by the definition of $\xi_j$,
\[ \{\xi_j \ge m\} \cap \{t_j = v\} = \mathcal{E}_{v,m} \cap \{t_j = v\}, \ \ \ j \ge 0. \]
Hence, for $j \ge 0$ and $n \ge 1$,
\[ \{t_{j+1} \le n\} = \bigcup_{v=0}^{n-1} \left( \{t_j = v\} \setminus \mathcal{E}_{v,n-v} \right). \]
By this and induction on $j$, each $t_j$ is a stopping time with respect to $(\mathcal{F}_n)_{n \ge 0}$.
Since $t_j \le t_{j+1}$, $t_j$ is $\mathcal{F}_{t_{j+1}}$-measurable, and so is $\xi_j$, which is $t_{j+1} - t_j - 1$ if $t_{j+1} < +\infty$, $+\infty$ if $t_j < t_{j+1} = +\infty$, and $0$ if $t_j = +\infty$.

\textit{Step 1.}
Let $M_n$ be the maximum of the integers $j \ge 0$ such that $t_j \le n$.
We show that
\begin{equation}\label{eq:sparse}
\lim_{n \to \infty} \frac{M_n}{n} = 0, \ \ \ \textup{$\mu_{\varphi}$-a.s.}
\end{equation}
Let $m \ge 1$ and $v \ge 0$.
For every word $\sigma$ of length $v$, $\mathcal{E}_{v,m} \cap J(\sigma)$ is the union of the intervals $J(\sigma i^m)$ over the parabolic letters $i$.
Hence, by Lemma \ref{lem:run} (3),
\[ \mu_{\varphi} \left( \mathcal{E}_{v,m} \, \middle| \, \mathcal{F}_v \right) \ge \frac{1}{1 + C_1 m}, \ \ \ \textup{$\mu_{\varphi}$-a.s.} \]
Let $j \ge 0$ and $A \in \mathcal{F}_{t_j}$.
Then, $A \cap \{t_j = v\} \in \mathcal{F}_v$ for every $v \ge 0$.
Since $t_j < +\infty$, $\mu_{\varphi}$-a.s., by the relation between $\xi_j$ and $\mathcal{E}_{v,m}$ above, we have that
\begin{align*}
\mu_{\varphi} (A \cap \{\xi_j \ge m\}) &= \sum_{v \ge 0} \mu_{\varphi} (A \cap \{t_j = v\} \cap \mathcal{E}_{v,m}) \\
&\ge \frac{1}{1 + C_1 m} \sum_{v \ge 0} \mu_{\varphi} (A \cap \{t_j = v\}) = \frac{\mu_{\varphi} (A)}{1 + C_1 m}.
\end{align*}
Hence, for every $j \ge 0$,
\[ \mu_{\varphi} \left( \xi_j \ge m \, \middle| \, \mathcal{F}_{t_j} \right) \ge \frac{1}{1 + C_1 m}, \ \ \ \textup{$\mu_{\varphi}$-a.s.} \]
Let $R \ge 1$ be an integer, and let
\[ c_R \coloneqq  \sum_{m=1}^{R} \frac{1}{1 + C_1 m}. \]
Then,
\[ E^{\mu_{\varphi}} \left[ \min\{\xi_j, R\} \, \middle| \, \mathcal{F}_{t_j} \right] = \sum_{m=1}^{R} \mu_{\varphi} \left( \xi_j \ge m \, \middle| \, \mathcal{F}_{t_j} \right) \ge c_R. \]
Let
\[ Z_j \coloneqq  \sum_{l=0}^{j-1} \left( \min\{\xi_l, R\} - E^{\mu_{\varphi}} \left[ \min\{\xi_l, R\} \, \middle| \, \mathcal{F}_{t_l} \right] \right), \ \ \ j \ge 0. \]
Then, $(Z_j)_{j \ge 0}$ is a martingale with respect to $(\mathcal{F}_{t_j})_{j \ge 0}$ and $|Z_{j+1} - Z_j| \le R$.
By Azuma's inequality (\cite{A67}, \cite[E14.2]{W91}), for every $\epsilon > 0$ and $j \ge 1$,
\[ \mu_{\varphi} (|Z_j| \ge \epsilon j) \le 2 \exp\left( - \frac{\epsilon^2 j}{2 R^2} \right). \]
By this and the Borel--Cantelli lemma, $Z_j / j$ converges to $0$ as $j \to \infty$, $\mu_{\varphi}$-a.s.
Since
\[ t_j \ge \sum_{l=0}^{j-1} \min\{\xi_l, R\} \ge Z_j + c_R \, j, \]
the limit inferior of $t_j / j$ as $j \to \infty$ is at least $c_R$, $\mu_{\varphi}$-a.s.
Since $c_R$ diverges to $+\infty$ as $R \to \infty$, $t_j / j$ diverges to $+\infty$ as $j \to \infty$, $\mu_{\varphi}$-a.s.
Since $M_n$ diverges to $+\infty$ as $n \to \infty$ and $n \ge t_{M_n}$, we obtain \eqref{eq:sparse}.

\textit{Step 2.}
For a word $\sigma = (\sigma_1, \dots, \sigma_n)$, let $V(\sigma)$ be the number of the integers $1 \le j \le n-1$ such that $\sigma_j \ne \sigma_{j+1}$, and let $V_n (x) \coloneqq  V((X_1 (x), \dots, X_n (x)))$.
If two adjacent letters differ, then at least one of them is a separating letter, since the same letter is repeated in each block.
Hence, $V_n \le 2 M_n$, and, by \eqref{eq:sparse}, $V_n / n$ converges to $0$ as $n \to \infty$, $\mu_{\varphi}$-a.s.

Let $\epsilon > 0$.
We take $\delta \in (0,1)$ such that $1 + (N-1) \delta < \exp (\epsilon / 4)$, and then $\delta^{\prime} > 0$ such that $\delta^{- \delta^{\prime}} < \exp (\epsilon / 4)$.
Let $n \ge 1$.
By taking the sum over the first letter and then over each of the following letters successively,
\[ \sum_{\sigma \in (\mathcal{I}_N)^n} \delta^{V(\sigma)} = N \left( 1 + (N-1) \delta \right)^{n-1}. \]
Each word $\sigma$ of length $n$ with $V(\sigma) \le \delta^{\prime} n$ contributes at least $\delta^{\delta^{\prime} n}$ to the above sum, and hence the number of such words is at most $N \exp (\epsilon n / 2)$.
By taking the sum of $\mu_{\varphi} (J(\sigma))$ over the words $\sigma$ among them such that $\mu_{\varphi} (J(\sigma)) \le \exp (- \epsilon n)$,
\[ \mu_{\varphi} \left( \left\{ x \in [0,1) : V_n (x) \le \delta^{\prime} n, \ \mu_{\varphi} (J_n (x)) \le \exp (- \epsilon n) \right\} \right) \le N \exp (- \epsilon n / 2). \]
By this, the Borel--Cantelli lemma and the convergence of $V_n / n$ to $0$, it holds $\mu_{\varphi}$-a.s.\ that $\mu_{\varphi} (J_n (x)) > \exp (- \epsilon n)$ for every sufficiently large $n$.
By this for $\epsilon = 1/k$, $k \ge 1$, the assertion holds.
\end{proof}

\begin{proof}[Proof of Theorem \ref{thm:para}]
Assertion (1) is Proposition \ref{prop:parabolic}.
We show (2).
Assume that $f$ is hyperbolic.
Let $n \ge 1$ and $x \in [0,1)$.
By \eqref{eq:length}, $\ell (J_n (x)) \le \left( G_{\max}^f \right)^n < 1$, and hence
\[ 0 \le \frac{\log \mu_{\varphi} (J_n (x))}{\log \ell (J_n (x))} \le \frac{- \log \mu_{\varphi} (J_n (x))}{- n \log G_{\max}^f}. \]
Hence, by (1),
\[ \lim_{n \to \infty} \frac{\log \mu_{\varphi} (J_n (x))}{\log \ell (J_n (x))} = 0, \ \ \ \textup{$\mu_{\varphi}$-a.s.} \, x. \]
By this and Lemma \ref{lem:dimension} (1) for $c = 0$, there exists a Borel set $K_0$ such that $\mu_{\varphi} (K_0) = 1$ and $\dim_H K_0 = 0$.
Hence, $\dim_H \mu_{\varphi} = 0$.
\end{proof}

\begin{proof}[Proof of Corollary \ref{cor:positive}]
Let $f$ be a hyperbolic LF system on $[0,1]$, and let $g$ be an LF system on $[0,1]$.
If $g$ is parabolic, then, by Theorem \ref{thm:para} (2), $\dim_H \mu_{\varphi} = 0$.
Assume that $g$ is hyperbolic, and let $\nu$ be as in Theorem \ref{thm:dim}.
By Lemma \ref{lem:id} (5), $0 < G_i^g (y) < 1$ for $y > -1$ and $i \in \mathcal{I}_N$.
Hence, $S$ is positive and continuous on $Y^g$, and, by Theorem \ref{thm:dim} (2), $\dim_H \mu_{\varphi} > 0$.
\end{proof}

\begin{Exa}\label{exa:Minkowski}
Let $N = 2$, and let
\[ f_0 (x) \coloneqq  \frac{x}{2}, \ \ \ f_1 (x) \coloneqq  \frac{x+1}{2}, \ \ \ g_0 (z) \coloneqq  \frac{z}{z+1}, \ \ \ g_1 (z) \coloneqq  \frac{1}{2-z}. \]
Then, $f$ is a hyperbolic LF system, and $g$ is an LF system for which
\[ (a_0^g, b_0^g, c_0^g) = (1, 0, 1), \ \ \ (a_1^g, b_1^g, c_1^g) = (0, 1/2, -1/2). \]
Since $g_0^{\prime} (0) = g_1^{\prime} (1) = 1$, $g$ is parabolic, both of the letters $0$ and $1$ are parabolic, and $\theta_0 = \theta_1 = 1$.
Minkowski's question-mark function is the solution of \eqref{eq:conj} in which $f$ and $g$ are exchanged (\cite[Section 1]{O26}), and hence, by \cite[Lemma 2.5]{O26}, $\varphi$ is the inverse function of Minkowski's question-mark function.
By Lemma \ref{lem:conj} (3), $\mu_{\varphi}$ is the image measure of $\ell$ under Minkowski's question-mark function.
By Theorem \ref{thm:para} (2), $\dim_H \mu_{\varphi} = 0$.
By Lemma \ref{lem:run} (1) for the empty word, $\mu_{\varphi} ([0, 2^{-m})) = 1/(m+1)$ for $m \ge 1$.

In \cite[Example 35]{O20}, it was shown that $\dim_H \mu_{\varphi} < 1$ for this pair, and it was not known whether $\dim_H \mu_{\varphi} > 0$.
The maps $g_0$ and $g_1$ are the linear fractional transformations given by the matrices
\[ \begin{pmatrix} 1 & 0 \\ 1 & 1 \end{pmatrix}, \ \ \ \begin{pmatrix} 0 & 1 \\ -1 & 2 \end{pmatrix}, \]
respectively, which are positive constant multiples of $A_0^g$ and $A_1^g$, and whose entries are integers and whose determinants are $1$.
Let $\sigma$ be a word, and let $a$, $b$ and $c$, $d$ be the entries of the first row and of the second row of the product of these matrices along $\sigma$, respectively.
Then, $a$, $b$, $c$ and $d$ are integers, $ad - bc = 1$, and, by Lemma \ref{lem:id} (1), $d > 0$ and $c + d > 0$.
We have that $g_{\sigma} (0) = b/d$, $g_{\sigma} (1) = (a+b)/(c+d)$ and $g_{\sigma} (1/2) = (a+2b)/(c+2d)$.
Since $(a+b) d - b (c+d) = 1$, the fractions $b/d$ and $(a+b)/(c+d)$ are in lowest terms, and $g_{\sigma} (1/2)$ is their mediant.
Since $g_0 (1) = g_1 (0) = 1/2$, by Lemma \ref{lem:cyl} (1) for $h = g$, the interval $[g_{\sigma} (0), g_{\sigma} (1))$ is divided at the point $g_{\sigma} (1/2)$ into the intervals $[g_{\sigma 0} (0), g_{\sigma 0} (1))$ and $[g_{\sigma 1} (0), g_{\sigma 1} (1))$.
Hence, for every $n \ge 1$, the intervals $\varphi (J(\sigma)) = [g_{\sigma} (0), g_{\sigma} (1))$, $\sigma \in (\mathcal{I}_2)^n$, are the Stern--Brocot intervals of order $n$ in \cite{KS07}.
Thus, for this pair, Theorem \ref{thm:para} (1) states that, for $\ell$-a.e.\ $z \in [0,1)$, the logarithm of the length of the Stern--Brocot interval of order $n$ containing $z$ divided by $n$ converges to $0$ as $n \to \infty$, which is known; see \cite[Section 1]{KS07}.
\end{Exa}

\section{Explicit formulas}\label{sec:tiling}

In this section, we show Theorem \ref{thm:tiling} and consider natural duals.
{\it Throughout this section, we assume that $f$ and $g$ are hyperbolic LF systems on $[0,1]$, and let $\varphi$ be the solution of \eqref{eq:conj}.}

\subsection{Proof of Theorem \ref{thm:tiling}}

Let
\[ \eta (t) \coloneqq  - t \log t, \ \ \ 0 < t \le 1. \]

\begin{Lem}\label{lem:Rogers}
$L$ is differentiable on $(0,1)$ and
\[ L^{\prime} (u) = \frac{\eta (u) + \eta (1-u)}{2 u (1-u)}, \ \ \ 0 < u < 1. \]
\end{Lem}

\begin{proof}
Since the derivative of the sum of $u^n / n^2$ over $n \ge 1$ is $- \log (1-u) / u$, we have that
\begin{align*}
L^{\prime} (u) &= - \frac{\log (1-u)}{u} + \frac{1}{2} \left( \frac{\log(1-u)}{u} - \frac{\log u}{1-u} \right) \\
&= - \frac{1}{2} \left( \frac{\log(1-u)}{u} + \frac{\log u}{1-u} \right) = \frac{\eta (u) + \eta (1-u)}{2u(1-u)}. \qedhere
\end{align*}
\end{proof}

\begin{Lem}\label{lem:chi}
Let $y > -1$.
Then, $0 < \chi_i (y) < 1$ for $1 \le i \le N-1$, and
\[ \frac{d}{dy} \left( 2 \sum_{i=1}^{N-1} L(\chi_i (y)) \right) = \frac{S(y)}{1+y}. \]
\end{Lem}

In the proof, we express $\chi_i (y)$ as the ratio of $\psi_y (b_i^g)$ to $\psi_y (b_{i+1}^g)$, and use the chain rule for the entropy.

\begin{proof}
For $0 \le i \le N$, let $P_i (y) \coloneqq  \psi_y (b_i^g)$.
By Lemma \ref{lem:id} (5), $G_i^g (y) = P_{i+1} (y) - P_i (y)$ for $i \in \mathcal{I}_N$, $P_1 (y) = G_0^g (y)$ and $P_N (y) = 1$.
By the definitions of $\psi_y$ and $\chi_i$, for $1 \le i \le N-1$,
\begin{equation}\label{eq:chi-P}
\chi_i (y) = \frac{P_i (y)}{P_{i+1} (y)}, \ \ \ 1 - \chi_i (y) = \frac{G_i^g (y)}{P_{i+1} (y)}.
\end{equation}
Hence, $0 < \chi_i (y) < 1$, and
\[ P_{i+1} (y) \left( \eta (\chi_i (y)) + \eta (1 - \chi_i (y)) \right) = \eta (P_i (y)) - \eta (P_{i+1} (y)) + \eta (G_i^g (y)). \]
By taking the sum over $1 \le i \le N-1$,
\begin{equation}\label{eq:chain}
\sum_{i=1}^{N-1} P_{i+1} (y) \left( \eta (\chi_i (y)) + \eta (1 - \chi_i (y)) \right) = \eta (G_0^g (y)) + \sum_{i=1}^{N-1} \eta (G_i^g (y)) = S(y).
\end{equation}
This is the chain rule for the entropy.
By \eqref{eq:chi-P},
\[ \frac{\chi_i (y)}{1 - \chi_i (y)} = \frac{P_i (y)}{G_i^g (y)} = \frac{b_i^g \left( 1 + b_{i+1}^g y \right)}{b_{i+1}^g - b_i^g}, \]
and hence,
\[ \frac{\chi_i^{\prime} (y)}{\chi_i (y) (1 - \chi_i (y))} = \frac{d}{dy} \log \frac{\chi_i (y)}{1 - \chi_i (y)} = \frac{b_{i+1}^g}{1 + b_{i+1}^g y} = \frac{P_{i+1} (y)}{1+y}. \]
By this and Lemma \ref{lem:Rogers},
\begin{align*}
\frac{d}{dy} \left( 2 L(\chi_i (y)) \right) &= \frac{\left( \eta (\chi_i (y)) + \eta (1 - \chi_i (y)) \right) \chi_i^{\prime} (y)}{\chi_i (y) (1 - \chi_i (y))} \\
&= \frac{P_{i+1} (y) \left( \eta (\chi_i (y)) + \eta (1 - \chi_i (y)) \right)}{1+y}.
\end{align*}
The assertion follows from this and \eqref{eq:chain}.
\end{proof}

The following shows that the integral of $S$ with respect to the first marginal of $\nu$ and the function $\rho$ in Subsection \ref{subsec:tau} are determined by a stationary measure for the dual maps of $g$.

\begin{Lem}\label{lem:marginal}
Let $\nu$ be a Borel probability measure on $Y^g \times Y^f$ satisfying \eqref{eq:stationary} for every continuous function $u$ on $Y^g \times Y^f$, let $\nu^g$ be the first marginal of $\nu$, and let $\nu^{\prime}$ be a Borel probability measure on $Y^g$ such that
\begin{equation}\label{eq:stationary-g}
\int_{Y^g} \sum_{i \in \mathcal{I}_N} G_i^g (y) \, v(H_i^g (y)) \, \nu^{\prime} (dy) = \int_{Y^g} v \, d\nu^{\prime}
\end{equation}
for every continuous function $v$ on $Y^g$. Then,\\
(1) For $z \in [0,1]$,
\[ \int_{Y^g} k(y,z) \, \nu^g (dy) = \int_{Y^g} k(y,z) \, \nu^{\prime} (dy). \]
(2) It holds that
\[ \int_{Y^g} S \, d\nu^g = \int_{Y^g} S \, d\nu^{\prime}. \]
\end{Lem}

\begin{proof}
Let $\rho$ and $m$ be as in Subsection \ref{subsec:tau}, and let $\widetilde{\rho}$ and $\widetilde{m}$ be defined in the same manner, where $\nu^g$ is replaced by $\nu^{\prime}$.
By \eqref{eq:stationary-g} for the functions $y \mapsto u(y,y)$, the image measure of $\nu^{\prime}$ under the map $y \mapsto (y,y)$ satisfies \eqref{eq:stationary} in which $f$ is replaced by $g$, and its first marginal is $\nu^{\prime}$.
Hence, Lemma \ref{lem:rho} and Propositions \ref{prop:tau} and \ref{prop:average} hold also in the case that $f$, $\nu$, $\rho$ and $m$ are replaced by $g$, this measure, $\widetilde{\rho}$ and $\widetilde{m}$, respectively, in which $\tau$ and $\Psi_g$ are not changed.
Thus $m$ and $\widetilde{m}$ are probability measures which are invariant and ergodic with respect to $\tau$ and are mutually absolutely continuous with $\ell$.
Since two ergodic invariant probability measures are equal to each other or mutually singular, $m = \widetilde{m}$.
Since $k$ is continuous on the compact set $Y^g \times [0,1]$, $\rho$ and $\widetilde{\rho}$ are continuous on $[0,1]$, and hence $\rho = \widetilde{\rho}$ on $[0,1]$, which is (1).
By \eqref{eq:lambda} and Proposition \ref{prop:average} in the two cases, both sides of (2) are equal to the integral of $\Psi_g$ with respect to $m$.
\end{proof}

By \eqref{eq:stationary-marginal} and Lemma \ref{lem:marginal} (1), the function $\rho$ and the measure $m$ in Subsection \ref{subsec:tau} and the measure $\widehat{\mu}$ in Subsection \ref{subsec:invariant} do not depend on the choice of $\nu$.

We now show Theorem \ref{thm:tiling}.
In the proof, by \eqref{eq:log-derivative} and the assumption, we show that the probability measure on $[e_{-}, e_{+}]$ whose density is proportional to $1/(1+y)$ satisfies \eqref{eq:stationary-g}.
Then, (1) follows from Lemma \ref{lem:marginal} (2) and Lemma \ref{lem:chi}, and (2) follows from Lemma \ref{lem:marginal} (1), by which the function $\rho$ is explicit.

\begin{proof}[Proof of Theorem \ref{thm:tiling}]
Since the dual maps of $g$ are continuous and strictly increasing on $(-1, +\infty)$, we have that
\[ H_i^g ([e_{-}, e_{+}]) = \left[ H_i^g (e_{-}), H_i^g (e_{+}) \right], \ \ \ i \in \mathcal{I}_N. \]
Since the union of these intervals is $[e_{-}, e_{+}]$, there exist $i, j \in \mathcal{I}_N$ such that $H_i^g (e_{-}) = e_{-}$ and $H_j^g (e_{+}) = e_{+}$.
Hence, by Lemma \ref{lem:fixed}, $e_{-} = e_i^g$ and $e_{+} = e_j^g$, and therefore $[e_{-}, e_{+}] \subset Y^g$.

Let $\zeta \coloneqq  \log (1 + e_{+}) - \log (1 + e_{-}) > 0$, and let $\nu^{\prime}$ be the Borel probability measure on $Y^g$ defined by
\[ \int_{Y^g} v \, d\nu^{\prime} \coloneqq  \frac{1}{\zeta} \int_{e_{-}}^{e_{+}} \frac{v(y)}{1+y} \, dy \]
for bounded Borel measurable functions $v$ on $Y^g$.
Let $v$ be a continuous function on $Y^g$ and $i \in \mathcal{I}_N$.
By \eqref{eq:log-derivative} and the change of variables $y^{\prime} = H_i^g (y)$,
\begin{align*}
\int_{e_{-}}^{e_{+}} \frac{G_i^g (y) \, v(H_i^g (y))}{1+y} \, dy &= \int_{e_{-}}^{e_{+}} v(H_i^g (y)) \, \frac{d}{dy} \log \left( 1 + H_i^g (y) \right) dy \\
&= \int_{H_i^g (e_{-})}^{H_i^g (e_{+})} \frac{v(y^{\prime})}{1+y^{\prime}} \, dy^{\prime}.
\end{align*}
By taking the sum over $i \in \mathcal{I}_N$ and by the assumption, $\nu^{\prime}$ satisfies \eqref{eq:stationary-g} for every continuous function $v$ on $Y^g$.

Let $\nu$ be a Borel probability measure on $Y^g \times Y^f$ satisfying \eqref{eq:stationary} for every continuous function $u$ on $Y^g \times Y^f$, which exists by Theorem \ref{thm:dim} (1), and let $\nu^g$, $\rho$ and $m$ be as in Subsection \ref{subsec:tau}.
By Lemma \ref{lem:marginal} (1), for $z \in [0,1]$,
\begin{equation}\label{eq:rho-explicit}
\rho (z) = \frac{1}{\zeta} \int_{e_{-}}^{e_{+}} \frac{dy}{(1+yz)^2} = \frac{e_{+} - e_{-}}{\zeta \, (1 + e_{-} z)(1 + e_{+} z)}.
\end{equation}
By Lemma \ref{lem:marginal} (2) and Lemma \ref{lem:chi},
\[ \int_{Y^g} S \, d\nu^g = \frac{1}{\zeta} \int_{e_{-}}^{e_{+}} \frac{S(y)}{1+y} \, dy = \frac{2}{\zeta} \sum_{i=1}^{N-1} \left( L(\chi_i (e_{+})) - L(\chi_i (e_{-})) \right). \]
Since the left hand side is the integral of $S(y)$ with respect to $\nu (dy \, dw)$, (1) follows from Theorem \ref{thm:dim} (2).

We show (2).
By Lemma \ref{lem:conj} (3), $\ell$ is the image measure of $\mu_{\varphi}$ under $\varphi$.
Hence, by \eqref{eq:lambda} and \eqref{eq:Psi-omega},
\[ \lambda_f = \int_0^1 \Psi_f (z) \, \rho (z) \, dz = - \int_{[0,1)} \omega_f (x) \, \rho (\varphi (x)) \, \mu_{\varphi} (dx). \]
By this, \eqref{eq:rho-explicit} and Proposition \ref{prop:limits}, (2) holds.
By Theorem \ref{thm:dim} (2), $\dim_H \mu_{\varphi}$ is the limit in (1) divided by the limit in (2), and hence (3) holds.
\end{proof}

\begin{Cor}\label{cor:common}
Under the assumption of Theorem \ref{thm:tiling}, assume furthermore that there exists $e > -1$ such that $H_i^f (e) = e$ for every $i \in \mathcal{I}_N$.
Then,
\[ \dim_H \mu_{\varphi} = \frac{2 \sum_{i=1}^{N-1} \left( L(\chi_i (e_{+})) - L(\chi_i (e_{-})) \right)}{\sum_{i \in \mathcal{I}_N} \left( \log (1 + H_i^g (e_{+})) - \log (1 + H_i^g (e_{-})) \right) \log \left( 1 / G_i^f (e) \right)}. \]
\end{Cor}

\begin{proof}
By Lemma \ref{lem:fixed}, $e = e_i^f$ for every $i \in \mathcal{I}_N$, and hence $e \in Y^f$.
Let $\zeta$ and $\nu^{\prime}$ be as in the proof of Theorem \ref{thm:tiling}, and let $\nu$ be the image measure of $\nu^{\prime}$ under the map $y \mapsto (y, e)$ from $Y^g$ to $Y^g \times Y^f$.
Since $H_i^f (e) = e$ for every $i \in \mathcal{I}_N$, by \eqref{eq:stationary-g} for the functions $v(y) \coloneqq  u(y,e)$, $\nu$ satisfies \eqref{eq:stationary} for every continuous function $u$ on $Y^g \times Y^f$.
By \eqref{eq:log-derivative}, for $i \in \mathcal{I}_N$,
\[ \int_{Y^g} G_i^g \, d\nu^{\prime} = \frac{1}{\zeta} \int_{e_{-}}^{e_{+}} \frac{G_i^g (y)}{1+y} \, dy = \frac{1}{\zeta} \left( \log (1 + H_i^g (e_{+})) - \log (1 + H_i^g (e_{-})) \right), \]
and hence
\begin{align*}
&\int_{Y^g \times Y^f} \Lambda (y,w) \, \nu(dy \, dw) \\
&= \frac{1}{\zeta} \sum_{i \in \mathcal{I}_N} \left( \log (1 + H_i^g (e_{+})) - \log (1 + H_i^g (e_{-})) \right) \log \frac{1}{G_i^f (e)}.
\end{align*}
The assertion follows from this, Theorem \ref{thm:dim} (2) and Theorem \ref{thm:tiling} (1).
\end{proof}

\begin{Rem}\label{rem:common}
(1) By \eqref{eq:one-plus-H},
\[ \log (1 + H_i^g (e_{+})) - \log (1 + H_i^g (e_{-})) = \log \frac{\left( 1 + b_{i+1}^g e_{+} \right) \left( 1 + b_i^g e_{-} \right)}{\left( 1 + b_{i+1}^g e_{-} \right) \left( 1 + b_i^g e_{+} \right)}, \ \ \ i \in \mathcal{I}_N. \]
(2) If $f$ consists of affine maps as in Remark \ref{rem:affine}, then the assumption for $f$ in Corollary \ref{cor:common} holds for $e = 0$, and $G_i^f (0) = r_i$.
In general, there exists $e > -1$ such that $H_i^f (e) = e$ for every $i \in \mathcal{I}_N$ if and only if there exists $s > -1$ such that $\psi_s \circ f_i \circ \psi_s^{-1}$ is an affine map for every $i \in \mathcal{I}_N$.
Indeed, if $H_i^f (e) = e$ for every $i \in \mathcal{I}_N$, then (iii) in Proposition \ref{prop:equiv} holds for $y = 0$ and $s = e$ and for the pair of $f$ and the LF system consisting of the affine maps with the slopes $G_i^f (e)$, $i \in \mathcal{I}_N$, in place of the pair of $f$ and $g$, and hence $\psi_e \circ f_i \circ \psi_e^{-1}$, $i \in \mathcal{I}_N$, are these affine maps.
Conversely, if $\psi_s \circ f_i \circ \psi_s^{-1}$ is an affine map for every $i \in \mathcal{I}_N$, then these maps form an LF system as in Remark \ref{rem:affine}, whose dual maps fix $0$, and, by (ii) in Proposition \ref{prop:equiv} for $y = 0$, $H_i^f (s) = s$ for every $i \in \mathcal{I}_N$.\\
(3) If the dual maps of $f$ have no common fixed points, then we do not know a general closed form for the integral in Theorem \ref{thm:tiling} (2).
It has a closed form in the following case.
If the assumption of Theorem \ref{thm:sing} (1) holds, then, by Lemma \ref{lem:psi} (1), the density $k(s_0, \cdot)$ of $\mu_{\varphi}$ with respect to $\ell$ is bounded away from $0$ and $+\infty$, and hence the two limits in Theorem \ref{thm:dim} (2) are equal to each other.
Hence, the limit in Theorem \ref{thm:tiling} (2) is equal to the right hand side of Theorem \ref{thm:tiling} (1), whether the dual maps of $f$ have a common fixed point or not.
This is the case if $f_i = g_i$ for every $i \in \mathcal{I}_N$, for which $\varphi$ is the identity map; the dual maps of the LF system in Example \ref{exa:non-natural} below have no common fixed points.\\
(4) Under the assumption of Theorem \ref{thm:tiling}, let $\zeta$ be as in the proof of Theorem \ref{thm:tiling}, and let $\widehat{\mu}$ be the measure in Subsection \ref{subsec:invariant}.
Then, by \eqref{eq:rho-explicit} and Proposition \ref{prop:hatmu} (1),
\[ \frac{d\widehat{\mu}}{d\mu_{\varphi}} (x) = \frac{e_{+} - e_{-}}{\zeta \, (1 + e_{-} \varphi (x))(1 + e_{+} \varphi (x))}, \ \ \ x \in [0,1), \]
and, by Remark \ref{rem:entropy}, the limit in Theorem \ref{thm:tiling} (2) is the Lyapunov exponent of $T$ with respect to $\widehat{\mu}$.
\end{Rem}

\subsection{Natural duals}\label{subsec:natural}

For a hyperbolic LF system $h$ on $[0,1]$, let
\[ \Theta_h (x) \coloneqq  \frac{(1 + e_0^h)(1 + e_{N-1}^h x)}{1 + e_0^h x} - 1, \ \ \ x \in [0,1], \]
and
\[ \delta_i^h \coloneqq  \left( (1 + e_0^h) \, e_{N-1}^h - e_0^h \right) b_i^h + e_0^h \left( 1 - a_i^h \right) - c_i^h, \ \ \ i \in \mathcal{I}_N. \]

\begin{Lem}\label{lem:Theta}
Let $h$ be a hyperbolic LF system on $[0,1]$. Then,\\
(1) $\Theta_h (0) = e_0^h$, $\Theta_h (1) = e_{N-1}^h$, and
\[ \frac{d}{dx} \log (1 + \Theta_h (x)) = \frac{e_{N-1}^h - e_0^h}{(1 + e_0^h x)(1 + e_{N-1}^h x)}, \ \ \ x \in [0,1]. \]
(2) $\delta_0^h = \delta_{N-1}^h = 0$.\\
(3) For $i \in \mathcal{I}_N$,
\begin{align*}
\Theta_h (b_i^h) - H_i^h (e_0^h) &= \frac{\delta_i^h}{b_i^h e_0^h + 1}, \\
\Theta_h (b_{i+1}^h) - H_i^h (e_{N-1}^h) &= \frac{\delta_i^h \left( 1 + b_{i+1}^h e_{N-1}^h \right)}{\left( 1 + b_{i+1}^h e_0^h \right) \left( b_i^h e_{N-1}^h + 1 \right)}.
\end{align*}
(4) If $i \in \mathcal{I}_N$ and $\delta_i^h = 0$, then $\Theta_h (h_i (x)) = H_i^h (\Theta_h (x))$ for $x \in [0,1]$.
\end{Lem}

\begin{proof}
We omit the superscript and the subscript $h$.
Let $\vartheta \coloneqq  (1 + e_0) \, e_{N-1} - e_0$.
Then, $\delta_i = \vartheta \, b_i + e_0 (1 - a_i) - c_i$ for $i \in \mathcal{I}_N$, and
\[ \Theta (x) = \frac{\vartheta x + e_0}{e_0 x + 1}, \ \ \ 1 + \Theta (x) = \frac{(1 + e_0)(1 + e_{N-1} x)}{1 + e_0 x}, \ \ \ x \in [0,1]. \]
By Lemma \ref{lem:fixed}, $e_0 > -1$ and $e_{N-1} > -1$, and hence $1 + e_0 x > 0$ and $1 + e_{N-1} x > 0$ for $x \in [0,1]$.

(1) This follows from the expression for $1 + \Theta (x)$.

(2) Since $b_0 = 0$, by the definition of $e_0$, $\delta_0 = e_0 (1 - a_0) - c_0 = 0$.
Since $a_{N-1} + b_{N-1} = c_{N-1} + 1$, we have that $1 - a_{N-1} = b_{N-1} - c_{N-1}$.
Hence, by Lemma \ref{lem:fixed},
\[ \delta_{N-1} = (\vartheta + e_0) \, b_{N-1} - (1 + e_0) \, c_{N-1} = (1 + e_0) (e_{N-1} b_{N-1} - c_{N-1}) = 0. \]

(3) By the definitions,
\[ \Theta (b_i) - H_i (e_0) = \frac{\vartheta \, b_i + e_0}{e_0 b_i + 1} - \frac{a_i e_0 + c_i}{b_i e_0 + 1} = \frac{\delta_i}{b_i e_0 + 1}. \]
By \eqref{eq:one-plus-H},
\[ \Theta (b_{i+1}) - H_i (e_{N-1}) = (1 + b_{i+1} e_{N-1}) \left( \frac{1 + e_0}{1 + b_{i+1} e_0} - \frac{c_i + 1}{b_i e_{N-1} + 1} \right). \]
Since $b_{i+1} (c_i + 1) = a_i + b_i$,
\begin{align*}
&(1 + e_0)(b_i e_{N-1} + 1) - (c_i + 1)(1 + b_{i+1} e_0) \\
&= (\vartheta + e_0) \, b_i + e_0 - c_i - e_0 (a_i + b_i) = \delta_i.
\end{align*}
The second identity follows from these.

(4) By the definitions of $\Theta$, $h_i$ and $H_i$,
\begin{align*}
\Theta (h_i (x)) &= \frac{(\vartheta \, a_i + e_0 c_i) x + \vartheta \, b_i + e_0}{(e_0 a_i + c_i) x + e_0 b_i + 1}, \\
H_i (\Theta (x)) &= \frac{(\vartheta \, a_i + e_0 c_i) x + e_0 a_i + c_i}{(\vartheta \, b_i + e_0) x + e_0 b_i + 1}.
\end{align*}
If $\delta_i = 0$, then $\vartheta \, b_i + e_0 = e_0 a_i + c_i$, and these are equal to each other.
\end{proof}

\begin{Prop}\label{prop:natural}
Let $h$ be a hyperbolic LF system on $[0,1]$.
Then, the following are equivalent:\\
(i) $H_{i-1}^h (e_{N-1}^h) = H_i^h (e_0^h)$ for every $1 \le i \le N-1$.\\
(ii) $\delta_i^h = 0$ for every $i \in \mathcal{I}_N$.\\
(iii) $\Theta_h (h_i (x)) = H_i^h (\Theta_h (x))$ for every $i \in \mathcal{I}_N$ and $x \in [0,1]$.
\end{Prop}

\begin{proof}
By Lemma \ref{lem:Theta} (4), (ii) implies (iii).

Assume that (iii) holds.
By Lemma \ref{lem:basic} (2) and Lemma \ref{lem:Theta} (1), by letting $x = 0$ and $x = 1$ in (iii), we have that $\Theta_h (b_i^h) = H_i^h (e_0^h)$ and $\Theta_h (b_{i+1}^h) = H_i^h (e_{N-1}^h)$ for every $i \in \mathcal{I}_N$.
Hence, $H_{i-1}^h (e_{N-1}^h) = \Theta_h (b_i^h) = H_i^h (e_0^h)$ for every $1 \le i \le N-1$.
Thus (i) holds.

Assume that (i) holds.
We show that $\delta_i^h = 0$ by induction on $i$.
By Lemma \ref{lem:Theta} (2), $\delta_0^h = 0$.
Let $1 \le i \le N-1$ and assume that $\delta_{i-1}^h = 0$.
Then, by Lemma \ref{lem:Theta} (3) and (i), $\Theta_h (b_i^h) = H_{i-1}^h (e_{N-1}^h) = H_i^h (e_0^h)$.
Hence, by Lemma \ref{lem:Theta} (3), $\delta_i^h = 0$.
Thus (ii) holds.
\end{proof}

By Lemma \ref{lem:Theta} (2), (ii) is equivalent to that $\delta_i^h = 0$ for every $1 \le i \le N-2$.
In particular, (i)--(iii) hold if $N = 2$.
If (i)--(iii) hold and $e_0^h \ne e_{N-1}^h$, then, by Lemma \ref{lem:Theta} (1), $\Theta_h$ is a linear fractional transformation which maps $[0,1]$ onto the closed interval with endpoints $e_0^h$ and $e_{N-1}^h$, and, by (iii), it conjugates $h_i$ to its dual map $H_i^h$ for every $i \in \mathcal{I}_N$.
The matrix $\begin{pmatrix} \vartheta & e_0^h \\ e_0^h & 1 \end{pmatrix}$ which gives $\Theta_h$ is symmetric, where $\vartheta$ is as in the proof of Lemma \ref{lem:Theta}.
See Remark \ref{rem:dual} below for the relation with natural duals in the sense of Schweiger \cite{S21}.

\begin{Cor}\label{cor:natural}
Assume that (i)--(iii) in Proposition \ref{prop:natural} hold for $h = g$. Then,\\
(1) If $e_0^g \ne e_{N-1}^g$, then the assumption of Theorem \ref{thm:tiling} holds for $e_{-} = \min\{e_0^g, e_{N-1}^g\}$ and $e_{+} = \max\{e_0^g, e_{N-1}^g\}$.\\
(2) If $e_0^g = e_{N-1}^g$, then $H_i^g (e_0^g) = e_0^g$ for every $i \in \mathcal{I}_N$, the function $\rho$ in Subsection \ref{subsec:tau} is $k(e_0^g, \cdot)$, and, for $\mu_{\varphi}$-a.s.\ $x$,
\begin{align*}
\lim_{n \to \infty} \frac{- \log \mu_{\varphi} (J_n (x))}{n} &= S(e_0^g), \\
\lim_{n \to \infty} \frac{- \log \ell (J_n (x))}{n} &= - \int_{[0,1)} k\left( e_0^g, \varphi (x^{\prime}) \right) \omega_f (x^{\prime}) \, \mu_{\varphi} (dx^{\prime}).
\end{align*}
\end{Cor}

\begin{proof}
(1) By Lemma \ref{lem:Theta} (1), $\Theta_g$ is a strictly monotone bijection from $[0,1]$ onto $[e_{-}, e_{+}]$.
By (iii) and Lemma \ref{lem:basic} (2),
\[ H_i^g ([e_{-}, e_{+}]) = H_i^g (\Theta_g ([0,1])) = \Theta_g (g_i ([0,1])) = \Theta_g \left( [b_i^g, b_{i+1}^g] \right), \ \ \ i \in \mathcal{I}_N. \]
Since the intervals $[b_i^g, b_{i+1}^g]$, $i \in \mathcal{I}_N$, have pairwise disjoint interiors and their union is $[0,1]$, the assumption of Theorem \ref{thm:tiling} holds.

(2) By Lemma \ref{lem:Theta} (1), $\Theta_g (x) = e_0^g$ for every $x \in [0,1]$.
Hence, by (iii), $H_i^g (e_0^g) = e_0^g$ for every $i \in \mathcal{I}_N$.
Then, by Lemma \ref{lem:id} (5), the Dirac measure at $e_0^g$ satisfies \eqref{eq:stationary-g} for every continuous function $v$ on $Y^g$.
Hence, by Lemma \ref{lem:marginal}, $\rho = k(e_0^g, \cdot)$, and the integral of $S$ with respect to the first marginal of a measure $\nu$ as in Theorem \ref{thm:dim} is $S(e_0^g)$.
The two limits are shown in the same manner as in the proof of Theorem \ref{thm:tiling}.
\end{proof}

\begin{Rem}[The case that $N = 2$]\label{rem:two}
Let $N = 2$.
For $h \in \{f, g\}$, let $\kappa_0^h \coloneqq  h_0^{\prime} (0)$ and $\kappa_1^h \coloneqq  h_1^{\prime} (1)$.
Then, $0 < \kappa_0^h < 1$, $0 < \kappa_1^h < 1$, and (i)--(iii) in Proposition \ref{prop:natural} hold.
Since $H_0^h (e_0^h) = e_0^h$ and $H_1^h (e_1^h) = e_1^h$, by \eqref{eq:one-plus-H} and Lemma \ref{lem:basic},
\[ G_0^h (e_0^h) = b_1^h \left( c_0^h + 1 \right) = \kappa_0^h, \ \ \ G_1^h (e_1^h) = \frac{1 - b_1^h}{c_1^h + 1} = \kappa_1^h. \]
By Lemma \ref{lem:id} (5) and Lemma \ref{lem:psi} (3), $G_0^h (y) = \psi_y (b_1^h)$, $G_1^h (y) = 1 - G_0^h (y)$, and
\[ \frac{G_0^h (y)}{G_1^h (y)} = \frac{b_1^h (1+y)}{1 - b_1^h}, \ \ \ y > -1. \]
Hence, $G_0^h$ is strictly increasing, $e_0^h = e_1^h$ if and only if $\kappa_0^h + \kappa_1^h = 1$, and
\begin{equation}\label{eq:zeta-two}
\log (1 + e_0^h) - \log (1 + e_1^h) = \log \frac{\kappa_0^h}{1 - \kappa_0^h} + \log \frac{\kappa_1^h}{1 - \kappa_1^h}.
\end{equation}
(1) We have that $\chi_1 = G_0^g$, and hence $\chi_1 (e_0^g) = \kappa_0^g$ and $\chi_1 (e_1^g) = 1 - \kappa_1^g$.
It holds that $L(u) + L(1-u) = \pi^2 / 6$ for $0 < u < 1$ (see \cite{K95, Z07}).
Therefore, by Corollary \ref{cor:natural}, Theorem \ref{thm:tiling} (1) and \eqref{eq:zeta-two} for $h = g$, for $\mu_{\varphi}$-a.s.\ $x$,
\[ \lim_{n \to \infty} \frac{- \log \mu_{\varphi} (J_n (x))}{n} = \frac{2 \left( L(\kappa_0^g) + L(\kappa_1^g) - \pi^2 / 6 \right)}{\log (\kappa_0^g / (1 - \kappa_0^g)) + \log (\kappa_1^g / (1 - \kappa_1^g))} \]
if $\kappa_0^g + \kappa_1^g \ne 1$, and this limit is $\eta (\kappa_0^g) + \eta (\kappa_1^g)$ if $\kappa_0^g + \kappa_1^g = 1$.
Thus, if $N = 2$, then the first limit in Theorem \ref{thm:dim} (2) is determined by the derivatives of $g_0$ and $g_1$ at their fixed points $0$ and $1$, respectively.\\
(2) Assume that $\kappa_0^f + \kappa_1^f = 1$.
Then, $e_0^f = e_1^f$, and the assumption for $f$ in Corollary \ref{cor:common} holds for $e = e_0^f$, for which $G_0^f (e) = \kappa_0^f$ and $G_1^f (e) = \kappa_1^f$.
Assume that $\kappa_0^g + \kappa_1^g \ne 1$, and let $\nu^{\prime}$ be the measure in the proof of Theorem \ref{thm:tiling}.
Since $G_0^g (y) / (1+y) = b_1^g / (1 + b_1^g y)$ and $1 + b_1^g y = (1 - b_1^g) / G_1^g (y)$ for $y > -1$, by \eqref{eq:zeta-two} for $h = g$,
\begin{align*}
\int_{Y^g} G_0^g \, d\nu^{\prime} &= \frac{\log (1 + b_1^g e_0^g) - \log (1 + b_1^g e_1^g)}{\log (1 + e_0^g) - \log (1 + e_1^g)} \\
&= \frac{\log (\kappa_1^g / (1 - \kappa_0^g))}{\log (\kappa_0^g / (1 - \kappa_0^g)) + \log (\kappa_1^g / (1 - \kappa_1^g))}.
\end{align*}
Hence, as in the proof of Corollary \ref{cor:common},
\[ \dim_H \mu_{\varphi} = \frac{2 \left( L(\kappa_0^g) + L(\kappa_1^g) - \pi^2 / 6 \right)}{\log (\kappa_1^g / (1 - \kappa_0^g)) \log (1 / \kappa_0^f) + \log (\kappa_0^g / (1 - \kappa_1^g)) \log (1 / \kappa_1^f)}. \]
If $\kappa_0^g + \kappa_1^g = 1$, then the Dirac measure at $(e_0^g, e_0^f)$ satisfies \eqref{eq:stationary} for every continuous function $u$ on $Y^g \times Y^f$, and, by Theorem \ref{thm:dim} (2),
\[ \dim_H \mu_{\varphi} = \frac{S(e_0^g)}{\Lambda (e_0^g, e_0^f)} = \frac{\kappa_0^g \log \kappa_0^g + \kappa_1^g \log \kappa_1^g}{\kappa_0^g \log \kappa_0^f + \kappa_1^g \log \kappa_1^f}. \]
Thus, if $N = 2$ and $\kappa_0^f + \kappa_1^f = 1$, then $\dim_H \mu_{\varphi}$ is determined by the four numbers $\kappa_0^f$, $\kappa_1^f$, $\kappa_0^g$ and $\kappa_1^g$.
By Remark \ref{rem:common} (2), $\kappa_0^f + \kappa_1^f = 1$ if and only if there exists $s > -1$ such that $\psi_s \circ f_i \circ \psi_s^{-1}$, $i \in \{0, 1\}$, are affine maps.
This is the case if $f_0 (x) = x/2$ and $f_1 (x) = (x+1)/2$; see \cite{O14} for this case.
\end{Rem}

\begin{Rem}[On the assumption of Theorem \ref{thm:tiling}]\label{rem:tiling}
(1) If the assumption of Theorem \ref{thm:tiling} holds, then $e_{-}$ and $e_{+}$ are the minimum and the maximum of $e_i^g$, $i \in \mathcal{I}_N$, respectively.
Indeed, by the proof of Theorem \ref{thm:tiling}, $e_{-}$ and $e_{+}$ are fixed points of some of the dual maps of $g$.
For every $i \in \mathcal{I}_N$, $H_i^g$ is a continuous map from $[e_{-}, e_{+}]$ to itself, and hence it has a fixed point in $[e_{-}, e_{+}]$, which is $e_i^g$ by Lemma \ref{lem:fixed}.\\
(2) Let $e_{-}$ and $e_{+}$ be the minimum and the maximum of $e_i^g$, $i \in \mathcal{I}_N$, respectively, and assume that $e_{-} < e_{+}$.
By Lemma \ref{lem:fixed}, for $y > -1$ and $i \in \mathcal{I}_N$, $H_i^g (y)$ is between $y$ and $e_i^g$, and hence $H_i^g ([e_{-}, e_{+}]) \subset [e_{-}, e_{+}]$.
By \eqref{eq:log-derivative} and Lemma \ref{lem:id} (5),
\begin{align*}
\sum_{i \in \mathcal{I}_N} \left( \log (1 + H_i^g (e_{+})) - \log (1 + H_i^g (e_{-})) \right) &= \int_{e_{-}}^{e_{+}} \frac{dy}{1+y} \\
&= \log (1 + e_{+}) - \log (1 + e_{-}).
\end{align*}
Thus the sum of the lengths of the images of the intervals $H_i^g ([e_{-}, e_{+}])$, $i \in \mathcal{I}_N$, under the map $y \mapsto \log (1+y)$ is equal to the length of the image of $[e_{-}, e_{+}]$.
Hence, the assumption of Theorem \ref{thm:tiling} holds if and only if the intervals $H_i^g ([e_{-}, e_{+}])$, $i \in \mathcal{I}_N$, have pairwise disjoint interiors.\\
(3) If $e_i^g = e_j^g$ for some $i \ne j$, then the assumption of Theorem \ref{thm:tiling} does not hold.
Indeed, assume that it holds for $e_{-}$ and $e_{+}$, and let $e \coloneqq  e_i^g = e_j^g$.
By (1), $e \in [e_{-}, e_{+}]$.
Since $H_i^g$ and $H_j^g$ are strictly increasing and fix $e$, the closed intervals $H_i^g ([e_{-}, e_{+}])$ and $H_j^g ([e_{-}, e_{+}])$ have positive lengths, and $e$ is an interior point of both of them if $e_{-} < e < e_{+}$, the left endpoint of both of them if $e = e_{-}$, and the right endpoint of both of them if $e = e_{+}$.
In each case, their interiors intersect.
In particular, if $N = 2$ and $e_0^g = e_1^g$, then Theorem \ref{thm:tiling} does not apply, and Corollary \ref{cor:natural} (2) applies.
\end{Rem}

\begin{Exa}\label{exa:non-natural}
Let $N = 3$, and let
\[ g_0 (z) \coloneqq  \frac{2z}{z+5}, \ \ \ g_1 (z) \coloneqq  \frac{z+1}{3}, \ \ \ g_2 (z) \coloneqq  \frac{5z+4}{3z+6}, \ \ \ z \in [0,1]. \]
Then, $g = \{g_i\}_{i \in \mathcal{I}_3}$ is a hyperbolic LF system on $[0,1]$ for which $(a_0^g, b_0^g, c_0^g) = (2/5, 0, 1/5)$, $(a_1^g, b_1^g, c_1^g) = (1/3, 1/3, 0)$ and $(a_2^g, b_2^g, c_2^g) = (5/6, 2/3, 1/2)$, and
\[ H_0^g (y) = \frac{2y+1}{5}, \ \ \ H_1^g (y) = \frac{y}{y+3}, \ \ \ H_2^g (y) = \frac{5y+3}{4y+6}, \ \ \ y > -1. \]
We have that $e_0^g = 1/3$, $e_1^g = 0$ and $e_2^g = 3/4$.
Since $H_0^g (e_2^g) = 1/2$ and $H_1^g (e_0^g) = 1/10$, (i) in Proposition \ref{prop:natural} does not hold for $h = g$.
On the other hand,
\[ H_1^g ([0, 3/4]) = [0, 1/5], \ \ \ H_0^g ([0, 3/4]) = [1/5, 1/2], \ \ \ H_2^g ([0, 3/4]) = [1/2, 3/4], \]
and hence the assumption of Theorem \ref{thm:tiling} holds for $e_{-} = 0$ and $e_{+} = 3/4$.
Since $\chi_1 (y) = (3 + 2y)/(6 + 2y)$ and $\chi_2 (y) = (2 + 2y)/(3 + 2y)$, by Theorem \ref{thm:tiling} (1), for every hyperbolic LF system $f$ on $[0,1]$ with $N = 3$,
\[ \lim_{n \to \infty} \frac{- \log \mu_{\varphi} (J_n (x))}{n} = \frac{2 \left( L(3/5) - L(1/2) + L(7/9) - L(2/3) \right)}{\log (7/4)}, \ \ \ \textup{$\mu_{\varphi}$-a.s.} \, x. \]
The right hand side is approximately $1.0833$.
This is the first limit in Theorem \ref{thm:dim} (2), and it is not a Hausdorff dimension.
For example, if $f_i (x) = (x+i)/3$ for $i \in \mathcal{I}_3$, then $\ell (J_n (x)) = 3^{-n}$, and, by Theorem \ref{thm:dim} (2), $\dim_H \mu_{\varphi}$ is this limit divided by $\log 3$, which is approximately $0.9861$.
\end{Exa}

\begin{Rem}[Dual maps]\label{rem:dual}
Let $\tau$ be the map in Section \ref{sec:pre}.
This is a piecewise linear fractional map with $N$ increasing branches whose inverse branches are $g_i$, $i \in \mathcal{I}_N$.
The construction of the invariant density $\rho$ of $\tau$ in Subsection \ref{subsec:tau} follows the method of dual algorithms in Schweiger \cite{S83, S21}.
In \cite{S21}, piecewise linear fractional maps $\mathcal{T}$ on $[0,1]$ each of whose branches is onto $(0,1)$ are considered, and a piecewise linear fractional map $\mathcal{T}^{*}$ is called a dual map of $\mathcal{T}$ if its inverse branches are given by the transposed matrices of the inverse branches of $\mathcal{T}$.
It is stated there that if there exists a dual map forming a fibred system on a set $B^{*}$, then the function
\begin{equation}\label{eq:Schweiger}
x \mapsto \int_{B^{*}} \frac{dy}{(1+xy)^2}
\end{equation}
is an invariant density of $\mathcal{T}$.
It is also stated there that if $\mathcal{T}$ has two branches, then there exists a linear fractional transformation $M$ given by a symmetric matrix such that $M \circ \mathcal{T} = \mathcal{T}^{*} \circ M$, for which $B^{*} = M([0,1])$, and that this is not always the case for maps with three or more branches.
Such a dual map is called a natural dual in \cite{S21}.
In this terminology, the functions $H_i^g$, $i \in \mathcal{I}_N$, which we call the dual maps of $g$, correspond to the inverse branches of a dual map of $\tau$, as $g_i$, $i \in \mathcal{I}_N$, are the inverse branches of $\tau$.\\
(1) Under the assumption of Theorem \ref{thm:tiling}, let $B^{*} \coloneqq  [e_{-}, e_{+}]$.
Since the interiors of the intervals $H_i^g (B^{*})$, $i \in \mathcal{I}_N$, are pairwise disjoint, a map $\mathcal{T}^{*}$ on the union of these interiors is defined by
\[ \mathcal{T}^{*} (H_i^g (y)) \coloneqq  y, \ \ \ e_{-} < y < e_{+}, \ \ i \in \mathcal{I}_N, \]
that is, $\mathcal{T}^{*}$ is the inverse function of $H_i^g$ on the interior of $H_i^g (B^{*})$ for each $i \in \mathcal{I}_N$.
By \eqref{eq:rho-explicit}, $\rho$ is equal to the function \eqref{eq:Schweiger} for this set $B^{*}$ up to a multiplicative constant.
The image measure of the measure $\nu^{\prime}$ in the proof of Theorem \ref{thm:tiling} under the map $y \mapsto \log (1+y)$ is the uniform distribution on $[\log (1 + e_{-}), \log (1 + e_{+})]$.\\
(2) If (i)--(iii) in Proposition \ref{prop:natural} hold for $h = g$ and $e_0^g \ne e_{N-1}^g$, then $\Theta_g$ is a linear fractional transformation given by a symmetric matrix which maps $[0,1]$ onto the closed interval $B^{*}$ with endpoints $e_0^g$ and $e_{N-1}^g$, and $\Theta_g \circ g_i = H_i^g \circ \Theta_g$ for every $i \in \mathcal{I}_N$.
Hence, the dual map is a natural dual in the sense of \cite{S21}.
This is always the case if $N = 2$ and $e_0^g \ne e_1^g$.
If (i)--(iii) in Proposition \ref{prop:natural} hold for $h = g$ and $e_0^g = e_{N-1}^g$, then the determinant $(1 + e_0^g)(e_{N-1}^g - e_0^g)$ of the symmetric matrix which gives $\Theta_g$ is $0$, $\Theta_g$ is the constant function whose value is $e_0^g$, and $e_0^g$ is a common fixed point of $H_i^g$, $i \in \mathcal{I}_N$ (Corollary \ref{cor:natural} (2)).
In this case, the set $B^{*}$ consists of the single point $e_0^g$.
In Schweiger \cite{S24}, a dual for which the set $B^{*}$ consists of a single point is regarded as a natural dual given by a singular matrix $M$, and it is called a singular dual.
This is the case if $N = 2$, $g_0 (z) = z/2$ and $g_1 (z) = (z+1)/2$, for which $e_0^g = e_1^g = 0$; see also Remark \ref{rem:tiling} (3).
Example \ref{exa:non-natural} shows that, for $N = 3$, the assumption of Theorem \ref{thm:tiling} may hold even if (i)--(iii) in Proposition \ref{prop:natural} fail.
\end{Rem}

\begin{Prop}\label{prop:nu-natural}
Assume that (i)--(iii) in Proposition \ref{prop:natural} hold for $h = f$ and for $h = g$, and let $m$ be the measure in Subsection \ref{subsec:tau}.
Then, the image measure $\nu$ of $m$ under the map $z \mapsto \left( \Theta_g (z), \Theta_f \left( \varphi^{-1} (z) \right) \right)$ is a Borel probability measure on $Y^g \times Y^f$ satisfying \eqref{eq:stationary} for every continuous function $u$ on $Y^g \times Y^f$.
\end{Prop}

In the proof, by \eqref{eq:log-derivative} and (iii) for $h = g$, we first show an identity for the functions $G_i^g \circ \Theta_g$ and the function $\rho$, which is \eqref{eq:reversible} below.
Then, we obtain \eqref{eq:stationary} by the change of variables $z^{\prime} = g_i (z)$.

\begin{proof}
By Lemma \ref{lem:Theta} (1), for $h \in \{f, g\}$, $\Theta_h$ is continuous and monotone on $[0,1]$, and hence $\Theta_h ([0,1])$ is the closed interval with endpoints $e_0^h$ and $e_{N-1}^h$, which is contained in $Y^h$.
Hence, $\nu$ is a Borel probability measure on $Y^g \times Y^f$.
Let $\rho$ be the function in Subsection \ref{subsec:tau}.
We first show that
\begin{equation}\label{eq:reversible}
G_i^g (\Theta_g (z)) \, \rho (z) = \rho (g_i (z)) \, g_i^{\prime} (z), \ \ \ z \in [0,1], \ \ i \in \mathcal{I}_N.
\end{equation}
Assume that $e_0^g \ne e_{N-1}^g$, and let
\[ \zeta^{\prime} \coloneqq  \log (1 + e_{N-1}^g) - \log (1 + e_0^g), \]
which is not $0$ and has the same sign as $e_{N-1}^g - e_0^g$.
Then, by Corollary \ref{cor:natural} (1), \eqref{eq:rho-explicit} and Lemma \ref{lem:Theta} (1),
\[ \rho (z) = \frac{e_{N-1}^g - e_0^g}{\zeta^{\prime} \, (1 + e_0^g z)(1 + e_{N-1}^g z)} = \frac{1}{\zeta^{\prime}} \, \frac{d}{dz} \log (1 + \Theta_g (z)), \ \ \ z \in [0,1]. \]
By \eqref{eq:log-derivative}, the chain rule and (iii) for $h = g$,
\begin{align*}
G_i^g (\Theta_g (z)) \, \frac{d}{dz} \log (1 + \Theta_g (z)) &= \frac{d}{dz} \log \left( 1 + H_i^g (\Theta_g (z)) \right) \\
&= \frac{d}{dz} \log \left( 1 + \Theta_g (g_i (z)) \right) = \zeta^{\prime} \rho (g_i (z)) \, g_i^{\prime} (z).
\end{align*}
Thus \eqref{eq:reversible} holds.
Assume that $e_0^g = e_{N-1}^g$.
Then, by Corollary \ref{cor:natural} (2), $\Theta_g (z) = e_0^g$, $\rho (z) = k(e_0^g, z)$ and $H_i^g (e_0^g) = e_0^g$, and \eqref{eq:reversible} follows from Lemma \ref{lem:id} (3) for $h = g$.

Let $u$ be a continuous function on $Y^g \times Y^f$, and let
\[ U(z) \coloneqq  u \! \left( \Theta_g (z), \Theta_f \left( \varphi^{-1} (z) \right) \right), \ \ \ z \in [0,1]. \]
By (iii) for $h = g$ and for $h = f$, and by $f_i \circ \varphi^{-1} = \varphi^{-1} \circ g_i$, for $z \in [0,1]$ and $i \in \mathcal{I}_N$,
\begin{align*}
u \! \left( H_i^g (\Theta_g (z)), H_i^f \left( \Theta_f \left( \varphi^{-1} (z) \right) \right) \right) &= u \! \left( \Theta_g (g_i (z)), \Theta_f \left( f_i \left( \varphi^{-1} (z) \right) \right) \right) \\
&= U(g_i (z)).
\end{align*}
Hence, by the definition of $\nu$, \eqref{eq:reversible}, the change of variables $z^{\prime} = g_i (z)$ and Lemma \ref{lem:basic} (2),
\begin{align*}
&\int_{Y^g \times Y^f} \sum_{i \in \mathcal{I}_N} G_i^g (y) \, u \! \left( H_i^g (y), H_i^f (w) \right) \nu(dy \, dw) \\
&= \sum_{i \in \mathcal{I}_N} \int_0^1 G_i^g (\Theta_g (z)) \, U(g_i (z)) \, \rho (z) \, dz = \sum_{i \in \mathcal{I}_N} \int_0^1 U(g_i (z)) \, \rho (g_i (z)) \, g_i^{\prime} (z) \, dz \\
&= \sum_{i \in \mathcal{I}_N} \int_{b_i^g}^{b_{i+1}^g} U(z^{\prime}) \, \rho (z^{\prime}) \, dz^{\prime} = \int_0^1 U(z) \, \rho (z) \, dz = \int_{Y^g \times Y^f} u \, d\nu. \qedhere
\end{align*}
\end{proof}

\begin{Rem}\label{rem:curve}
Assume that (i)--(iii) in Proposition \ref{prop:natural} hold for $h = f$ and for $h = g$, and let $\nu$ be the measure in Proposition \ref{prop:nu-natural}.\\
(1) By Theorem \ref{thm:dim} (2),
\[ \dim_H \mu_{\varphi} = \frac{\int_0^1 S(\Theta_g (z)) \, \rho (z) \, dz}{\int_0^1 \Lambda \left( \Theta_g (z), \Theta_f \left( \varphi^{-1} (z) \right) \right) \rho (z) \, dz}, \]
where $\rho$ is given by \eqref{eq:rho-explicit} if $e_0^g \ne e_{N-1}^g$, and $\rho = k(e_0^g, \cdot)$ if $e_0^g = e_{N-1}^g$.\\
(2) The measure $\nu$ is concentrated on the curve consisting of the points $\left( \Theta_g (\varphi (x)), \Theta_f (x) \right)$, $x \in [0,1]$.
Assume that $e_0^g \ne e_{N-1}^g$, and let $\Theta \coloneqq  \Theta_f \circ \varphi^{-1} \circ \Theta_g^{-1}$, which is a continuous function on the closed interval with endpoints $e_0^g$ and $e_{N-1}^g$.
Then, this curve is the graph of $\Theta$, and, by (iii) for $h = g$ and for $h = f$, and by $f_i \circ \varphi^{-1} = \varphi^{-1} \circ g_i$,
\[ H_i^f \circ \Theta = \Theta \circ H_i^g, \ \ \ i \in \mathcal{I}_N. \]
Thus $\Theta$ is a solution of the conjugate equation for the dual maps of $g$ and of $f$, and $\nu$ is the image measure of its first marginal under the map $y \mapsto (y, \Theta (y))$.\\
(3) Assume that $e_0^g \ne e_{N-1}^g$ and the assumption of Theorem \ref{thm:sing} (1) holds.
Then, the function $\Theta$ in (2) is the restriction of $\gamma_{s_0}$.
Indeed, by (ii) in Proposition \ref{prop:equiv} and Lemma \ref{lem:fixed}, $\gamma_{s_0} (e_i^g) = e_i^f$ for every $i \in \mathcal{I}_N$.
Since $1 + e \, \psi_{s_0} (x) = (1 + \gamma_{s_0} (e) \, x) / (1 + s_0 x)$ for $e > -1$ and $x \in [0,1]$, by the definitions of $\Theta_f$ and $\Theta_g$, we have that $1 + \Theta_f (x) = (1 + s_0) \left( 1 + \Theta_g (\psi_{s_0} (x)) \right)$, that is, $\Theta_f = \gamma_{s_0} \circ \Theta_g \circ \psi_{s_0}$.
Since $\varphi = \psi_{s_0}$, we have that $\Theta = \gamma_{s_0}$ on the interval in (2).
\end{Rem}

\section{Remarks}\label{sec:rem}

\begin{Rem}[Affine maps]\label{rem:affine}
Let $r_i$, $i \in \mathcal{I}_N$, be positive numbers such that $r_0 + \cdots + r_{N-1} = 1$, and let
\[ f_i (x) \coloneqq  r_i x + (r_0 + \cdots + r_{i-1}), \ \ \ x \in [0,1], \ \ i \in \mathcal{I}_N. \]
Then, $f = \{f_i\}_{i \in \mathcal{I}_N}$ is a hyperbolic LF system for which $a_i^f = r_i$, $b_i^f = r_0 + \cdots + r_{i-1}$ and $c_i^f = 0$.
We have that $H_i^f (0) = 0$, $G_i^f (0) = r_i$ and $e_i^f = 0$ for every $i \in \mathcal{I}_N$, and hence $Y^f = \{0\}$.
Let $g$ be a hyperbolic LF system on $[0,1]$.
Then, \eqref{eq:conj} is the equation \cite[(41)]{O20} considered in \cite[Subsection 3.6]{O20}, in which the numbers $r_i$ and the maps $f_i$ are denoted by $q_i$ and $S_i$, respectively.\\
(1) A Borel probability measure $\nu$ on $Y^g \times \{0\}$ is identified with a Borel probability measure on $Y^g$, and \eqref{eq:stationary} is the equation that
\[ \int_{Y^g} \sum_{i \in \mathcal{I}_N} G_i^g (y) \, v(H_i^g (y)) \, \nu(dy) = \int_{Y^g} v \, d\nu \]
for every continuous function $v$ on $Y^g$.
Furthermore,
\[ \Lambda (y,0) = - \sum_{i \in \mathcal{I}_N} G_i^g (y) \log r_i, \ \ \ y \in Y^g, \]
and $- s_N \left( \mathbf{G}^g (y) | \mathbf{G}^f (0) \right)$ is the relative entropy of $\mathbf{G}^g (y)$ with respect to $(r_i)_{i \in \mathcal{I}_N}$.
Thus, in Theorem \ref{thm:dim}, the numbers $r_i$ appear as the values $G_i^f (0)$, $i \in \mathcal{I}_N$, at the common fixed point $0$ of $H_i^f$, $i \in \mathcal{I}_N$.\\
(2) The assumption of Theorem \ref{thm:sing} (1) holds if and only if there exists $e > -1$ such that $H_i^g (e) = e$ and $G_i^g (e) = r_i$ for every $i \in \mathcal{I}_N$, and then
\[ \varphi (x) = \frac{x}{1 + e - e x}, \ \ \ x \in [0,1]. \]
Indeed, assume that there exists such a number $e$, and let $s > -1$ be the number such that $(1+s)(1+e) = 1$.
Then, $\gamma_s (e) = 0$, $G_i^g (e) = r_i = G_i^f (\gamma_s (e))$ and $\gamma_s \left( H_i^g (e) \right) = 0 = H_i^f (\gamma_s (e))$ for every $i \in \mathcal{I}_N$.
Hence, (iii) in Proposition \ref{prop:equiv} holds for $y = e$, and, by Proposition \ref{prop:equiv}, $g_i \circ \psi_s = \psi_s \circ f_i$ for every $i \in \mathcal{I}_N$.
Then, as is stated after Proposition \ref{prop:equiv}, $s = s_0$.
Thus the assumption of Theorem \ref{thm:sing} (1) holds, and $\varphi (x) = \psi_s (x) = x / (1 + e - ex)$.
Conversely, assume that the assumption of Theorem \ref{thm:sing} (1) holds, and let $e > -1$ be the number such that $(1 + s_0)(1+e) = 1$.
Then, $\gamma_{s_0} (e) = 0$, and, by (ii) in Proposition \ref{prop:equiv} for $s = s_0$ and $y = e$, $G_i^g (e) = G_i^f (0) = r_i$ and $\gamma_{s_0} \left( H_i^g (e) \right) = H_i^f (0) = 0 = \gamma_{s_0} (e)$, and hence $H_i^g (e) = e$, for every $i \in \mathcal{I}_N$.

The existence of such a number $e$ is the condition that both of (42) and (43) in \cite[Theorem 43]{O20} hold.
Thus Theorem \ref{thm:sing} corresponds to \cite[Theorem 43]{O20} for hyperbolic LF systems $g$.\\
(3) If the assumption of Theorem \ref{thm:tiling} holds, then, by Corollary \ref{cor:common} for $e = 0$,
\[ \dim_H \mu_{\varphi} = \frac{2 \sum_{i=1}^{N-1} \left( L(\chi_i (e_{+})) - L(\chi_i (e_{-})) \right)}{\sum_{i \in \mathcal{I}_N} \left( \log (1 + H_i^g (e_{+})) - \log (1 + H_i^g (e_{-})) \right) \log (1/r_i)}. \]
\end{Rem}

\begin{Rem}[The inverse function]\label{rem:inverse}
Let $f$ and $g$ be LF systems on $[0,1]$, and let $\varphi$ be the solution of \eqref{eq:conj}.
By \cite[Lemma 2.5]{O26}, the inverse function $\varphi^{-1}$ of $\varphi$ is the solution of \eqref{eq:conj} in which $f$ and $g$ are exchanged.
Let $\mu_{\varphi^{-1}}$ be the Borel probability measure on $[0,1]$ whose distribution function is $\varphi^{-1}$.
By Lemma \ref{lem:conj} (3) in which $f$ and $g$ are exchanged, $\mu_{\varphi^{-1}} (B) = \ell \left( \varphi^{-1} (B) \right)$ for every Borel subset $B$ of $[0,1]$.
Thus $\mu_{\varphi^{-1}}$ is the image measure of $\ell$ under $\varphi$, whereas, by Lemma \ref{lem:conj} (3), $\ell$ is the image measure of $\mu_{\varphi}$ under $\varphi$.\\
(1) Assume that $f$ and $g$ are hyperbolic.
Then, Theorem \ref{thm:dim} in which $f$ and $g$ are exchanged gives $\dim_H \mu_{\varphi^{-1}}$.
In this case, the measure in Theorem \ref{thm:dim} is a stationary distribution of the transition kernel on $Y^f \times Y^g$ by which $(w,y)$ moves to $\left( H_i^f (w), H_i^g (y) \right)$ with probability $G_i^f (w)$, $i \in \mathcal{I}_N$.\\
(2) $\mu_{\varphi}$ is singular with respect to the Lebesgue measure if and only if $\mu_{\varphi^{-1}}$ is singular with respect to the Lebesgue measure.
Here, $f$ and $g$ are not assumed to be hyperbolic.
Indeed, if $B$ is a Borel subset of $[0,1]$ such that $\ell (B) = 0$ and $\mu_{\varphi} (B) = 1$, then, by the above, $\ell (\varphi (B)) = \mu_{\varphi} (B) = 1$ and $\mu_{\varphi^{-1}} (\varphi (B)) = \ell (B) = 0$.
The converse is shown in the same manner.
See also the statement after \cite[Lemma 2.5]{O26}.
For hyperbolic LF systems $f$ and $g$, this agrees with Theorem \ref{thm:sing} and Theorem \ref{thm:sing} in which $f$ and $g$ are exchanged.
Indeed, let $s_0^{\prime} > -1$ be the number such that $(1+s_0)(1+s_0^{\prime}) = 1$.
Since $\psi_{s_0^{\prime}}$ is the inverse function of $\psi_{s_0}$, we have that $\psi_{s_0^{\prime}} (g_0 (1)) = f_0 (1)$, and $g_i \circ \psi_{s_0} = \psi_{s_0} \circ f_i$ for every $i \in \mathcal{I}_N$ if and only if $f_i \circ \psi_{s_0^{\prime}} = \psi_{s_0^{\prime}} \circ g_i$ for every $i \in \mathcal{I}_N$.\\
(3) See Remarks \ref{rem:Furstenberg} and \ref{rem:projective} for the measure $\mu_{\varphi^{-1}}$ in the case that $f$ consists of affine maps.
\end{Rem}

\begin{Rem}[On the assumptions]\label{rem:assumption}
(1) Let $f$ and $g$ be hyperbolic LF systems on $[0,1]$.
For $n \ge 0$, $x \in [0,1)$ and $h \in \{f, g\}$, let $y_n^h (x) \coloneqq  y_{\sigma}^h$ for $\sigma = (X_1 (x), \dots, X_n (x))$.
By Lemma \ref{lem:id} (1) and (6), Lemma \ref{lem:conj} and Lemma \ref{lem:coding}, under $\mu_{\varphi}$, $\left( (y_n^g, y_n^f) \right)_{n \ge 0}$ is a Markov chain with the initial value $(0,0)$ and with the transition kernel by which $(y,w)$ moves to $\left( H_i^g (y), H_i^f (w) \right)$ with probability $G_i^g (y)$, $i \in \mathcal{I}_N$, and, for $n \ge 1$ and $x \in [0,1)$,
\[ \mu_{\varphi} (J_n (x)) = \prod_{j=1}^{n} G_{X_j (x)}^g \left( y_{j-1}^g (x) \right), \ \ \ \ell (J_n (x)) = \prod_{j=1}^{n} G_{X_j (x)}^f \left( y_{j-1}^f (x) \right). \]
The measure $\nu_n$ in the proof of Lemma \ref{lem:stationary} is the average of the distributions of $(y_j^g, y_j^f)$, $0 \le j \le n-1$, under $\mu_{\varphi}$.
The sets $Y^g$ and $Y^f$ are intervals for the variables $y$ and $w$, respectively, and no relation between them is assumed in Theorems \ref{thm:dim} and \ref{thm:sing}.\\
(2) The hyperbolicity of $g$ is used for the compactness of $Y^g$, for Lemma \ref{lem:rho} (1) and the ergodicity in Proposition \ref{prop:tau}, and for the first limit in Proposition \ref{prop:limits}.
The hyperbolicity of $f$ is used for the compactness of $Y^f$, for Lemma \ref{lem:dimension}, and for the second limit in Proposition \ref{prop:limits}.\\
(3) In Theorem \ref{thm:para} (2), the hyperbolicity of $f$ is used only for Lemma \ref{lem:dimension} (1) and \eqref{eq:length}.
We do not consider the Hausdorff dimension of $\mu_{\varphi}$ in the case that $f$ is parabolic, and we only give two remarks for it.
First, if $f_i = g_i$ for every $i \in \mathcal{I}_N$, then $\varphi (x) = x$ and $\mu_{\varphi} = \ell$, whether $f$ is hyperbolic or not.
Hence, if both of $f$ and $g$ are parabolic, then $\dim_H \mu_{\varphi} = 0$ does not hold in general, in contrast to Theorem \ref{thm:para} (2).
Second, Minkowski's question-mark function is the solution of \eqref{eq:conj} for $f_0 (x) = x/(x+1)$, $f_1 (x) = 1/(2-x)$, $g_0 (x) = x/2$ and $g_1 (x) = (x+1)/2$ (\cite[Section 1]{O26}).
In this case, $f = \{f_0, f_1\}$ is a parabolic LF system, and $g = \{g_0, g_1\}$ is a hyperbolic LF system.
See Example \ref{exa:Minkowski} for the inverse function of Minkowski's question-mark function.
For this pair, Kinney \cite{K60} showed that
\[ \dim_H \mu_{\varphi} = \frac{\log 2}{2 \int_{[0,1]} \log (1+x) \, \mu_{\varphi} (dx)}. \]
(4) In \cite[Section 3]{O26}, $g$ is an LF system and $f$ is a weakly contractive D-system, and, in \cite[Subsection 3.6]{O20}, $g$ is an LF system and $f$ consists of affine maps, whereas both of $f$ and $g$ are LF systems in this paper.
The functions $G_i^g$, $H_i^g$ and $\mathbf{G}^g$ are the functions $G_i$, $H_i$ and $\mathbf{G}$ in \cite[Section 3]{O26}.
Let $h = \{h_i\}_{i \in \mathcal{I}_N}$ be an LF system on $[0,1]$.
If the strict inequality holds in the condition (2) of LF systems for every $i \in \mathcal{I}_N$, which is the condition (2') in \cite[Definition 2.8]{O26}, then $h_i^{\prime} < 1$ on $[0,1]$ for every $i \in \mathcal{I}_N$, and hence $h$ is hyperbolic.
\end{Rem}

\begin{Rem}[Place-dependent probabilities]\label{rem:Blackwell}
Let $f$ and $g$ be hyperbolic LF systems on $[0,1]$.\\
(1) The equation \eqref{eq:stationary} states that $\nu$ is a stationary distribution for the maps $(y,w) \mapsto \left( H_i^g (y), H_i^f (w) \right)$, $i \in \mathcal{I}_N$, which are chosen with the probabilities $G_i^g (y)$, $i \in \mathcal{I}_N$, depending on the place $y$.
Blackwell \cite{B57} expressed the entropy of a function of a stationary Markov chain with finitely many states as an integral with respect to a stationary distribution of this kind on a simplex, which is called a Blackwell measure; see also \cite{BPS12}.
The first limit in Theorem \ref{thm:dim} (2) is of the same form.
Indeed, by Remark \ref{rem:assumption} (1), $G_i^g (y_n^g (x))$ is the conditional probability under $\mu_{\varphi}$ of the event that $X_{n+1} = i$ given $X_1, \dots, X_n$, and $S(y)$ is the entropy of the probability vector $\mathbf{G}^g (y)$.\\
(2) The measure $\nu$ in Theorem \ref{thm:dim} is unique.
Indeed, by \eqref{eq:log-derivative}, for $h \in \{f, g\}$ and $i \in \mathcal{I}_N$, the derivative of the function $\theta \mapsto \log \left( 1 + H_i^h (\exp (\theta) - 1) \right)$ at $\theta \in \mathbb{R}$ is $G_i^h (\exp (\theta) - 1)$.
Hence, with respect to the metric on $Y^g \times Y^f$ for which the distance between $(y,w)$ and $(y^{\prime}, w^{\prime})$ is
\[ \left| \log (1+y) - \log (1+y^{\prime}) \right| + \left| \log (1+w) - \log (1+w^{\prime}) \right|, \]
the Lipschitz constants of the maps in (1) are at most $\max\{G_{\max}^g, G_{\max}^f\}$, which is smaller than $1$ by Lemma \ref{lem:hyp} (4).
The functions $G_i^g$, $i \in \mathcal{I}_N$, are continuously differentiable on the compact interval $Y^g$, which is contained in $(-1, +\infty)$.
Hence, the probabilities $G_i^g (y)$, $i \in \mathcal{I}_N$, are Lipschitz continuous functions of $(y,w)$ with respect to this metric, and their values are at least $G_{\min}^g > 0$.
Hence, the uniqueness follows from \cite[Theorem 2]{B57}.
\end{Rem}

\begin{Rem}[The measure $\mu_{\varphi^{-1}}$ and Furstenberg measures]\label{rem:Furstenberg}
In this remark and in Remark \ref{rem:projective}, we assume that the source system $f$ of \eqref{eq:conj} consists of affine maps as in Remark \ref{rem:affine}, and that the target system $g$ is an LF system on $[0,1]$, which may be parabolic.
The reason for this assumption on $f$ is stated at the end of (1).
Let $\mu_{\varphi^{-1}}$ be as in Remark \ref{rem:inverse}, which is the image measure of $\ell$ under $\varphi$.
For a word $\sigma = (\sigma_1, \dots, \sigma_n)$, let $r_{\sigma} \coloneqq  r_{\sigma_1} \cdots r_{\sigma_n}$, which is $1$ if $n = 0$, and let $r_{\max}$ be the maximum of $r_i$, $i \in \mathcal{I}_N$.
Then, $\ell (J(\sigma)) = r_{\sigma}$, and hence $\mu_{\varphi^{-1}} (\varphi (J(\sigma))) = r_{\sigma}$, for every word $\sigma$.\\
(1) For a subset $B$ of $[0,1]$ and $i \in \mathcal{I}_N$, we denote by $g_i^{-1} (B)$ the set of the points $z \in [0,1]$ such that $g_i (z) \in B$, and we define $f_i^{-1} (B)$ in the same manner.
Then, for every Borel subset $B$ of $[0,1]$,
\begin{equation}\label{eq:F-stationary}
\mu_{\varphi^{-1}} (B) = \sum_{i \in \mathcal{I}_N} r_i \, \mu_{\varphi^{-1}} \left( g_i^{-1} (B) \right).
\end{equation}
Indeed, by \eqref{eq:conj}, $f_i^{-1} \left( \varphi^{-1} (B) \right) = \varphi^{-1} \left( g_i^{-1} (B) \right)$ for $i \in \mathcal{I}_N$.
Since $f_i$ is an affine map with the slope $r_i$ from $[0,1]$ onto $[b_i^f, b_{i+1}^f]$,
\[ \ell \left( \varphi^{-1} (B) \cap [b_i^f, b_{i+1}^f] \right) = r_i \, \ell \left( \varphi^{-1} \left( g_i^{-1} (B) \right) \right) = r_i \, \mu_{\varphi^{-1}} \left( g_i^{-1} (B) \right), \ \ \ i \in \mathcal{I}_N. \]
By taking the sum over $i \in \mathcal{I}_N$, we obtain \eqref{eq:F-stationary}.

The equation \eqref{eq:F-stationary} states that $\mu_{\varphi^{-1}}$ is a stationary measure on $[0,1]$ for the probability measure on the set of $2 \times 2$ real matrices which assigns the mass $r_i$ to $A_i^g$ for each $i \in \mathcal{I}_N$, where the matrix $A_i^g$ acts on $[0,1]$ by $g_i$; see \cite[Section 4]{O14}.
Stationary measures of this kind are called Furstenberg measures; see \cite{BPS12, HS17}.
Since $\ell (J(\sigma)) = r_{\sigma}$ for every word $\sigma$, by Lemma \ref{lem:coding} (1), under $\ell$, $X_n$, $n \ge 1$, are independent, and each of them takes the value $i$ with probability $r_i$, $i \in \mathcal{I}_N$.
For $x \in [0,1)$ and $n \ge 1$, by Lemma \ref{lem:coding} (1) and Lemma \ref{lem:conj} (2), $\varphi (x)$ is in the interval $\varphi (J_n (x)) = [g_{\sigma} (0), g_{\sigma} (1))$, where $\sigma = (X_1 (x), \dots, X_n (x))$.
Since $\ell (J_n (x)) \le r_{\max}^n$ and $\varphi$ is uniformly continuous, the length of this interval converges to $0$ as $n \to \infty$.
Hence, $\varphi (x)$ is the limit of $g_{X_1 (x)} \circ \cdots \circ g_{X_n (x)} (0)$ as $n \to \infty$, and $\mu_{\varphi^{-1}}$ is the distribution of the limit of the compositions of the maps $g_i$, $i \in \mathcal{I}_N$, chosen independently with probabilities $r_i$, $i \in \mathcal{I}_N$.

Furthermore, $\mu_{\varphi^{-1}}$ is the unique Borel probability measure on $[0,1]$ satisfying \eqref{eq:F-stationary} for every Borel subset $B$ of $[0,1]$.
Indeed, let $\mu$ be such a measure.
By Lemma \ref{lem:basic} (1) and (2), the set $g_i^{-1} (\{1\})$ is $\{1\}$ if $i = N-1$, and it is empty otherwise.
Hence, $\mu (\{1\}) = r_{N-1} \, \mu (\{1\})$, and, since $r_{N-1} < 1$, we have that $\mu (\{1\}) = 0$.
Let $\sigma = (\sigma_1, \dots, \sigma_n)$ be a word of length $n \ge 1$, and let $\sigma^{\prime} \coloneqq  (\sigma_2, \dots, \sigma_n)$.
By Lemma \ref{lem:conj} (2), $\varphi (J(\sigma)) = [g_{\sigma} (0), g_{\sigma} (1)) = g_{\sigma_1} (\varphi (J(\sigma^{\prime})))$, and this interval is contained in $[b_{\sigma_1}^g, b_{\sigma_1 + 1}^g)$.
Since $g_{\sigma_1}$ is injective, $g_{\sigma_1}^{-1} (\varphi (J(\sigma))) = \varphi (J(\sigma^{\prime}))$, and, since $g_i ([0,1)) = [b_i^g, b_{i+1}^g)$, we have that $g_i^{-1} (\varphi (J(\sigma))) \subset \{1\}$ for $i \ne \sigma_1$.
Hence, by \eqref{eq:F-stationary}, $\mu (\varphi (J(\sigma))) = r_{\sigma_1} \, \mu (\varphi (J(\sigma^{\prime})))$, and, by induction on $n$, $\mu (\varphi (J(\sigma))) = r_{\sigma} = \mu_{\varphi^{-1}} (\varphi (J(\sigma)))$ for every word $\sigma$.
By Lemma \ref{lem:cyl} (3) for $h = g$, two of the intervals $\varphi (J(\sigma))$ are disjoint or one of them contains the other, and hence the family consisting of these intervals and the empty set is closed under finite intersections.
Since the lengths of $J(\sigma)$, $\sigma \in (\mathcal{I}_N)^n$, are at most $r_{\max}^n$ and $\varphi$ is a homeomorphism, this family generates the Borel $\sigma$-algebra on $[0,1)$.
Hence, by the Dynkin $\pi$-$\lambda$ theorem, $\mu = \mu_{\varphi^{-1}}$.
If $N = 2$, $r_0 = r_1 = 1/2$ and the strict inequality holds in the condition (2) of LF systems for $g$, then \eqref{eq:F-stationary} and this uniqueness are shown in \cite[Section 4]{O14}.

If the source system does not consist of affine maps, then the constants $r_i$ in \eqref{eq:F-stationary} are replaced by functions.
Indeed, let $f$ be an LF system on $[0,1]$.
Then, for $i \in \mathcal{I}_N$ and every Borel subset $B^{\prime}$ of $[0,1]$, $\ell \left( B^{\prime} \cap [b_i^f, b_{i+1}^f] \right)$ is the integral of $f_i^{\prime}$ over $f_i^{-1} (B^{\prime})$ with respect to the Lebesgue measure, and hence, in the same manner as above,
\[ \mu_{\varphi^{-1}} (B) = \sum_{i \in \mathcal{I}_N} \int_{g_i^{-1} (B)} f_i^{\prime} \left( \varphi^{-1} (z) \right) \mu_{\varphi^{-1}} (dz) \]
for every Borel subset $B$ of $[0,1]$.
The sum of the weights $f_i^{\prime} \left( \varphi^{-1} (z) \right)$ over $i \in \mathcal{I}_N$ is $1$ for every $z \in [0,1]$ if and only if $f$ consists of affine maps.
Indeed, by Lemma \ref{lem:basic} (3) and Lemma \ref{lem:psi} (1), $f_i^{\prime} (x) = (b_{i+1}^f - b_i^f) \, k(c_i^f, x)$ for $x \in [0,1]$ and $i \in \mathcal{I}_N$.
If $c_i^f = 0$, then $f_i^{\prime}$ is the constant $b_{i+1}^f - b_i^f$, and, if $c_i^f \ne 0$, then $f_i^{\prime}$ is strictly convex on $[0,1]$.
Since the sum of $b_{i+1}^f - b_i^f$ over $i \in \mathcal{I}_N$ is $1$, the sum of $f_i^{\prime}$ over $i \in \mathcal{I}_N$ is the constant $1$ if $c_i^f = 0$ for every $i \in \mathcal{I}_N$, and it is strictly convex on $[0,1]$ otherwise.
Furthermore, these weights cannot be replaced by other constants.
Indeed, assume that there exist real numbers $p_i$, $i \in \mathcal{I}_N$, such that
\[ \mu_{\varphi^{-1}} (B) = \sum_{i \in \mathcal{I}_N} p_i \, \mu_{\varphi^{-1}} \left( g_i^{-1} (B) \right) \]
for every Borel subset $B$ of $[0,1]$.
Let $i \in \mathcal{I}_N$, let $E$ be a Borel subset of $(0,1)$, and let $B \coloneqq  g_i (\varphi (E))$.
Then, $B$ is contained in $(b_i^g, b_{i+1}^g)$, $g_i^{-1} (B) = \varphi (E)$, and $g_j^{-1} (B)$ is empty for $j \ne i$.
By \eqref{eq:conj}, $B = \varphi (f_i (E))$.
Hence,
\[ \ell (f_i (E)) = \mu_{\varphi^{-1}} (B) = p_i \, \mu_{\varphi^{-1}} (\varphi (E)) = p_i \, \ell (E). \]
Since the left hand side is the integral of $f_i^{\prime}$ over $E$ with respect to the Lebesgue measure and $f_i^{\prime}$ is continuous, we have that $f_i^{\prime} = p_i$ on $[0,1]$.
Thus $f$ consists of affine maps, and $p_i$ is the slope of $f_i$ for every $i \in \mathcal{I}_N$.
Thus the description of $\mu_{\varphi^{-1}}$ by maps chosen independently with constant probabilities is obtained in the case that the source side of \eqref{eq:conj} is linear, and we restrict ourselves to this case in this remark and in Remark \ref{rem:projective}.\\
(2) By Remark \ref{rem:inverse} (2), $\mu_{\varphi^{-1}}$ is singular with respect to the Lebesgue measure if and only if $\mu_{\varphi}$ is singular with respect to the Lebesgue measure.
Hence, if $g$ is parabolic, then, by Theorem \ref{thm:para} (2), $\mu_{\varphi^{-1}}$ is singular with respect to the Lebesgue measure.
If $g$ is hyperbolic, then, by Theorem \ref{thm:sing} and Remark \ref{rem:affine} (2), the following holds: if there exists $e > -1$ such that $H_i^g (e) = e$ and $G_i^g (e) = r_i$ for every $i \in \mathcal{I}_N$, then
\[ \varphi^{-1} (z) = \frac{(1 + e) z}{1 + e z}, \ \ \ z \in [0,1], \]
and $\mu_{\varphi^{-1}}$ is absolutely continuous with respect to the Lebesgue measure, and otherwise $\mu_{\varphi^{-1}}$ is singular with respect to the Lebesgue measure.
If $N = 2$, $r_0 = r_1 = 1/2$ and the strict inequality holds in the condition (2) of LF systems for $g$, then this corresponds to \cite[Theorem 4.3]{O14}.\\
(3) Assume that $g$ is hyperbolic.
Then, by Remark \ref{rem:inverse} (1) and Remark \ref{rem:affine},
\[ \dim_H \mu_{\varphi^{-1}} = \frac{\sum_{i \in \mathcal{I}_N} r_i \log r_i}{\int_{Y^g} \sum_{i \in \mathcal{I}_N} r_i \log G_i^g (y) \, \nu_r (dy)} \]
for every Borel probability measure $\nu_r$ on $Y^g$ such that
\[ \int_{Y^g} \sum_{i \in \mathcal{I}_N} r_i \, v(H_i^g (y)) \, \nu_r (dy) = \int_{Y^g} v \, d\nu_r \]
for every continuous function $v$ on $Y^g$, and such a measure $\nu_r$ exists.
Thus $\nu_r$ is a stationary measure for the dual maps $H_i^g$, $i \in \mathcal{I}_N$, chosen independently with the constant probabilities $r_i$, $i \in \mathcal{I}_N$.
In contrast, the measure $\nu$ in Theorem \ref{thm:dim} is a stationary measure for the dual maps chosen with the probabilities $G_i^g (y)$, $i \in \mathcal{I}_N$, which depend on the place $y$ (Remark \ref{rem:Blackwell}).
\end{Rem}

\begin{Rem}[Stationary measures on the projective line]\label{rem:projective}
Let $f$, $g$, $\mu_{\varphi^{-1}}$ and $r_{\sigma}$ be as in Remark \ref{rem:Furstenberg}; in particular, $f$ consists of affine maps.
Let $\mathbb{P}^1 (\mathbb{R}) \coloneqq  \mathbb{R} \cup \{\infty\}$ be the one-point compactification of $\mathbb{R}$.
For a $2 \times 2$ real matrix $A = \begin{pmatrix} a & b \\ c & d \end{pmatrix}$ with $\det A \ne 0$, the map $z \mapsto (az + b)/(cz + d)$ extends to a homeomorphism of $\mathbb{P}^1 (\mathbb{R})$, which maps $\infty$ to $a/c$ and $-d/c$ to $\infty$ if $c \ne 0$, and which fixes $\infty$ if $c = 0$.
We denote the extensions of $g_i$ and $H_i^g$, which are given by $A_i^g$ and its transposed matrix, by $\overline{g}_i$ and $\overline{H}_i$, respectively, and let $\overline{g}_{\sigma} \coloneqq  \overline{g}_{\sigma_1} \circ \cdots \circ \overline{g}_{\sigma_n}$ for a word $\sigma = (\sigma_1, \dots, \sigma_n)$.
If we identify $z \in \mathbb{R}$ and $\infty$ with the lines in $\mathbb{R}^2$ spanned by the vectors ${}^t (z, 1)$ and ${}^t (1, 0)$, respectively, then $\mathbb{P}^1 (\mathbb{R})$ is the real projective line, and $\overline{g}_i$ is the map induced by the linear map given by $A_i^g$.
We say that a Borel probability measure $\mu$ on $\mathbb{P}^1 (\mathbb{R})$ is stationary if
\[ \mu (B) = \sum_{i \in \mathcal{I}_N} r_i \, \mu \left( \overline{g}_i^{-1} (B) \right) \]
for every Borel subset $B$ of $\mathbb{P}^1 (\mathbb{R})$.
By \eqref{eq:F-stationary}, $\mu_{\varphi^{-1}}$, regarded as a measure on $\mathbb{P}^1 (\mathbb{R})$, is stationary.\\
(1) Let $A$ be as above.
A point $z \in \mathbb{R}$ is a fixed point of the extension of the map given by $A$ if and only if $c z^2 + (d - a) z - b = 0$, and $\infty$ is a fixed point of it if and only if $c = 0$.
By this for $A$ and for its transposed matrix, a point $z \in \mathbb{P}^1 (\mathbb{R})$ is a fixed point of the extension of the map given by $A$ if and only if $-1/z$ is a fixed point of the extension of the map given by the transposed matrix of $A$, where we let $-1/0 \coloneqq  \infty$ and $-1/\infty \coloneqq  0$.
By \cite[Lemma 2.9]{O26}, $g_0$ is a weak contraction on $[0,1]$ in the sense of \cite{O26}.
Since $g_0 (0) = 0$, $0$ is the unique fixed point of $g_0$ in $[0,1]$, and, by Lemma \ref{lem:basic} (2), $g_{N-1} (0) = b_{N-1}^g > 0$.
Hence, $\overline{g}_i$, $i \in \mathcal{I}_N$, have no common fixed points in $[0,1]$, and therefore $\overline{H}_i$, $i \in \mathcal{I}_N$, have no common fixed points in $\mathbb{P}^1 (\mathbb{R}) \setminus (-1, +\infty)$.
Thus $\overline{g}_i$, $i \in \mathcal{I}_N$, have a common fixed point in $\mathbb{P}^1 (\mathbb{R})$ if and only if $H_i^g$, $i \in \mathcal{I}_N$, have a common fixed point in $(-1, +\infty)$.
If $e \in (-1, +\infty)$ is a common fixed point of $H_i^g$, $i \in \mathcal{I}_N$, then $-1/e$ is a common fixed point of $\overline{g}_i$, $i \in \mathcal{I}_N$, which is not in $[0,1]$, and the Dirac measure at $-1/e$ is stationary.
Thus stationary measures on $\mathbb{P}^1 (\mathbb{R})$ are not unique in this case.
For example, if $g$ consists of affine maps, then $H_i^g (0) = c_i^g = 0$ and $\overline{g}_i (\infty) = \infty$ for every $i \in \mathcal{I}_N$, and the Dirac measure at $\infty$ is stationary.
This is the case if $g_i = f_i$ for every $i \in \mathcal{I}_N$, for which $\mu_{\varphi^{-1}} = \ell$.\\
(2) Assume that $H_i^g$, $i \in \mathcal{I}_N$, have no common fixed points in $(-1, +\infty)$.
Then, $\mu_{\varphi^{-1}}$ is the unique stationary measure on $\mathbb{P}^1 (\mathbb{R})$.
Indeed, we first show that for every $z \in \mathbb{P}^1 (\mathbb{R})$, there exists a word $\sigma$ such that $\overline{g}_{\sigma} (z) \in [0,1]$.
Since $b_0^g = 0$,
\[ \frac{1}{\overline{g}_0 (z)} = \frac{1}{a_0^g} \cdot \frac{1}{z} + \frac{c_0^g}{a_0^g}, \ \ \ z \in \mathbb{P}^1 (\mathbb{R}), \]
where we let $1/0 \coloneqq  \infty$ and $1/\infty \coloneqq  0$.
Thus, in the coordinate $1/z$, $\overline{g}_0$ is the affine map with the slope $1/a_0^g$ and the intercept $c_0^g / a_0^g$, which fixes $\infty$.
By Lemma \ref{lem:basic}, $a_0^g = g_0^{\prime} (0) = b_1^g (c_0^g + 1)$.
Hence, if $g_0^{\prime} (0) < 1$, then the slope is larger than $1$, and, if $g_0^{\prime} (0) = 1$, then the slope is $1$ and $c_0^g = (1 - b_1^g)/b_1^g > 0$.
In both cases, this affine map has at most one fixed point in $\mathbb{R}$, and the absolute values of its $n$-fold compositions at every other point of $\mathbb{R}$ diverge to $+\infty$ as $n \to \infty$.
Hence, $\overline{g}_0$ has at most two fixed points, one of which is $0$, and, if $z$ is not a fixed point of $\overline{g}_0$, then $\overline{g}_{0^n} (z)$ converges to $0$ as $n \to \infty$.
By Lemma \ref{lem:basic} (2), $\overline{g}_1 (0) = b_1^g \in (0,1)$.
Since $\overline{g}_1$ is continuous, if $z = 0$ or $z$ is not a fixed point of $\overline{g}_0$, then $\overline{g}_1 \left( \overline{g}_{0^n} (z) \right) \in [0,1]$ for some $n \ge 0$.
Let $z \ne 0$ be a fixed point of $\overline{g}_0$.
By (1) and the assumption, $z$ is not a common fixed point of $\overline{g}_i$, $i \in \mathcal{I}_N$, and hence there exists $i \in \mathcal{I}_N$ such that $\overline{g}_i (z) \ne z$.
Then, $\overline{g}_i (z)$ is $0$ or it is not a fixed point of $\overline{g}_0$, and hence $\overline{g}_1 \left( \overline{g}_{0^n} \left( \overline{g}_i (z) \right) \right) \in [0,1]$ for some $n \ge 0$.

Let $\mu$ be a stationary measure on $\mathbb{P}^1 (\mathbb{R})$, and let $\mathcal{O} \coloneqq  \mathbb{P}^1 (\mathbb{R}) \setminus [0,1]$.
By induction on $n$,
\[ \mu (B) = \sum_{\sigma \in (\mathcal{I}_N)^n} r_{\sigma} \, \mu \left( \overline{g}_{\sigma}^{-1} (B) \right) \]
for every $n \ge 0$ and every Borel subset $B$ of $\mathbb{P}^1 (\mathbb{R})$.
Since $g_i ([0,1]) \subset [0,1]$ for every $i \in \mathcal{I}_N$, we have that $\overline{g}_{\sigma}^{-1} (\mathcal{O}) \subset \mathcal{O}$ for every word $\sigma$.
Since $r_{\sigma} > 0$ and the sum of $r_{\sigma}$ over $\sigma \in (\mathcal{I}_N)^n$ is $1$, by the above for $B = \mathcal{O}$, we have that $\mu \left( \overline{g}_{\sigma}^{-1} (\mathcal{O}) \right) = \mu (\mathcal{O})$, that is, $\mu \left( \mathcal{O} \setminus \overline{g}_{\sigma}^{-1} (\mathcal{O}) \right) = 0$, for every word $\sigma$.
By the claim above, $\mathcal{O}$ is the union of the sets $\mathcal{O} \setminus \overline{g}_{\sigma}^{-1} (\mathcal{O})$ over the words $\sigma$, which are countably many.
Hence, $\mu (\mathcal{O}) = 0$.
Therefore, the restriction of $\mu$ to $[0,1]$ satisfies \eqref{eq:F-stationary} for every Borel subset $B$ of $[0,1]$, and, by Remark \ref{rem:Furstenberg} (1), $\mu = \mu_{\varphi^{-1}}$.\\
(3) By (1) and (2), $\mu_{\varphi^{-1}}$ is the unique stationary measure on $\mathbb{P}^1 (\mathbb{R})$ if and only if $H_i^g$, $i \in \mathcal{I}_N$, have no common fixed points in $(-1, +\infty)$.
By Lemma \ref{lem:run} (1) and (2), this is the case if $g$ is parabolic.
In \cite{HS17}, the Furstenberg measure is the unique stationary measure on the projective line for a probability measure on $SL_2 (\mathbb{R})$ whose support generates an unbounded group which does not preserve any nonempty finite subset of the projective line.
Let $\widehat{A}_i \coloneqq  (\det A_i^g)^{-1/2} A_i^g$ for $i \in \mathcal{I}_N$, which induces the map $\overline{g}_i$.
If $H_i^g$, $i \in \mathcal{I}_N$, have no common fixed points in $(-1, +\infty)$, then the probability measure which assigns the mass $r_i$ to the matrix $\widehat{A}_i$ for each $i \in \mathcal{I}_N$ satisfies this assumption, and hence $\mu_{\varphi^{-1}}$ is the Furstenberg measure in the sense of \cite{HS17}.
Indeed, if a nonempty finite subset of $\mathbb{P}^1 (\mathbb{R})$ were mapped onto itself by $\overline{g}_i$ for every $i \in \mathcal{I}_N$, then the uniform distribution on it would be stationary, which contradicts (2), since $\mu_{\varphi^{-1}}$ has no atoms.
Furthermore, the lower right entry of the $n$-th power of $\widehat{A}_0$ is $(a_0^g)^{-n/2}$, and, if $a_0^g = 1$, then its lower left entry is $n c_0^g$, where $c_0^g > 0$.
Hence, these powers are unbounded.
It is shown in \cite{HS17} that, for such probability measures supported on finitely many matrices with algebraic entries, the dimension of the Furstenberg measure is the minimum of $1$ and the random walk entropy divided by twice the Lyapunov exponent.

Assume that $g$ is hyperbolic.
Then, the formula for $\dim_H \mu_{\varphi^{-1}}$ in Remark \ref{rem:Furstenberg} (3) is of this form.
Indeed, for a word $\sigma = (\sigma_1, \dots, \sigma_n)$, let $\widehat{A}_{\sigma} \coloneqq  \widehat{A}_{\sigma_1} \cdots \widehat{A}_{\sigma_n}$.
By Lemma \ref{lem:cyl} (2) for $h = g$, the maps $g_{\sigma}$, $\sigma \in (\mathcal{I}_N)^n$, are distinct, and hence so are the matrices $\widehat{A}_{\sigma}$, $\sigma \in (\mathcal{I}_N)^n$.
Hence, the entropy of the distribution of the product of $n$ independent random matrices each of which is $\widehat{A}_i$ with probability $r_i$, $i \in \mathcal{I}_N$, is $n$ times the entropy of the probability vector $(r_i)_{i \in \mathcal{I}_N}$, and the numerator in Remark \ref{rem:Furstenberg} (3) multiplied by $-1$ is the random walk entropy.
Let $\| \cdot \|$ be the Euclidean norm on $\mathbb{R}^2$, and we denote the operator norm of a $2 \times 2$ real matrix $A$ with respect to it by $\| A \|$.
For $z \in [0,1]$ and a word $\sigma$, by the proof of Lemma \ref{lem:id} (2), $A_{\sigma}^g \, {}^t (z, 1) = (c_{\sigma}^g z + d_{\sigma}^g) \, {}^t (g_{\sigma} (z), 1)$ and $g_{\sigma}^{\prime} (z) = \det A_{\sigma}^g / (c_{\sigma}^g z + d_{\sigma}^g)^2$, and hence
\[ \left\| \widehat{A}_{\sigma} \, {}^t (z, 1) \right\|^2 = \frac{1 + g_{\sigma} (z)^2}{g_{\sigma}^{\prime} (z)}. \]
By this and Lemma \ref{lem:hyp} (3) for $h = g$,
\[ \left| \log \left( g_{\sigma} (1) - g_{\sigma} (0) \right) + 2 \log \left\| \widehat{A}_{\sigma} \, {}^t (z, 1) \right\| \right| \le \log (2 C_g). \]
We use this for $z = 0$ and for $z = 1$.
Since ${}^t (0, 1)$ is a unit vector, the lower bound for the norm of $\widehat{A}_{\sigma} \, {}^t (0, 1)$ is also a lower bound for $\| \widehat{A}_{\sigma} \|$.
Since ${}^t (0, 1)$ and ${}^t (1, 1)$ form a basis of $\mathbb{R}^2$, $\| \widehat{A}_{\sigma} \|$ is at most a constant multiple of the sum of the norms of $\widehat{A}_{\sigma} \, {}^t (0, 1)$ and $\widehat{A}_{\sigma} \, {}^t (1, 1)$.
Hence, there exists a constant $C_2 \ge 1$ depending only on $g$ such that
\begin{equation}\label{eq:norm}
C_2^{-1} \le \left( g_{\sigma} (1) - g_{\sigma} (0) \right) \| \widehat{A}_{\sigma} \|^2 \le C_2
\end{equation}
for every word $\sigma$.

Theorem \ref{thm:dim} in which $f$ and $g$ are exchanged is an assertion for $\mu_{\varphi^{-1}}$-a.s.\ points of the target side and for the intervals $[g_{\sigma} (0), g_{\sigma} (1))$, and we rewrite it in terms of the points $x$ of the source side.
Let $x \in [0,1)$, $n \ge 1$ and $\sigma \coloneqq  (X_1 (x), \dots, X_n (x))$.
By Lemma \ref{lem:coding} (1) and Lemma \ref{lem:conj} (2), the interval of the form $[g_{\sigma^{\prime}} (0), g_{\sigma^{\prime}} (1))$, $\sigma^{\prime} \in (\mathcal{I}_N)^n$, which contains $\varphi (x)$ is $\varphi (J_n (x)) = [g_{\sigma} (0), g_{\sigma} (1))$.
Since $\mu_{\varphi^{-1}}$ is the image measure of $\ell$ under $\varphi$, by the second limit in Theorem \ref{thm:dim} (2) in which $f$ and $g$ are exchanged, for $\ell$-a.s.\ $x$, $- \log \left( g_{\sigma} (1) - g_{\sigma} (0) \right)$ for $\sigma = (X_1 (x), \dots, X_n (x))$ divided by $n$ converges, as $n \to \infty$, to the denominator in Remark \ref{rem:Furstenberg} (3) multiplied by $-1$.
Hence, by \eqref{eq:norm}, for $\ell$-a.s.\ $x$, $\log \| \widehat{A}_{X_1 (x)} \cdots \widehat{A}_{X_n (x)} \|$ divided by $n$ converges to half of it.

Since $\det \widehat{A}_i = 1$ for every $i \in \mathcal{I}_N$, the operator norm of every product of $n$ of these matrices is at least $1$, and it is at most the $n$-th power of the maximum of $\| \widehat{A}_i \|$, $i \in \mathcal{I}_N$.
Hence, for every $x \in [0,1)$ and every $n \ge 1$,
\[ 0 \le \frac{1}{n} \log \left\| \widehat{A}_{X_1 (x)} \cdots \widehat{A}_{X_n (x)} \right\| \le \max_{i \in \mathcal{I}_N} \log \| \widehat{A}_i \|, \]
and, by the bounded convergence theorem, the expectation of the middle term with respect to $\ell$ converges to the same limit as $n \to \infty$.
For each $n \ge 1$, since $X_1, \dots, X_n$ are independent and identically distributed under $\ell$, the distribution of the matrix $\widehat{A}_{X_1} \cdots \widehat{A}_{X_n}$ is equal to that of $\widehat{A}_{X_n} \cdots \widehat{A}_{X_1}$, and hence this expectation does not depend on the order of the product.
We use this identity of the distributions only for each fixed $n$, and not for the almost sure convergence.
The Lyapunov exponent of the products of the independent random matrices is the limit of these expectations as $n \to \infty$.
Thus, if $f$ consists of affine maps and $g$ is hyperbolic, then $\dim_H \mu_{\varphi^{-1}}$ is the random walk entropy divided by twice the Lyapunov exponent.
We do not consider $\dim_H \mu_{\varphi^{-1}}$ in this paper in the case that $g$ is parabolic.\\
(4) Assume that $H_i^g$, $i \in \mathcal{I}_N$, have a common fixed point $e$ in $(-1, +\infty)$.
Then, $g$ is hyperbolic by Lemma \ref{lem:run}, and $e \in Y^g$ by Lemma \ref{lem:fixed}.
The Dirac measure at $e$ satisfies the equation for $\nu_r$ in Remark \ref{rem:Furstenberg} (3), and the Dirac measure at $(e, 0)$ satisfies \eqref{eq:stationary}.
Hence, by Remark \ref{rem:Furstenberg} (3) and Theorem \ref{thm:dim} (2),
\[ \dim_H \mu_{\varphi^{-1}} = \frac{\sum_{i \in \mathcal{I}_N} r_i \log r_i}{\sum_{i \in \mathcal{I}_N} r_i \log G_i^g (e)}, \ \ \ \dim_H \mu_{\varphi} = \frac{\sum_{i \in \mathcal{I}_N} G_i^g (e) \log G_i^g (e)}{\sum_{i \in \mathcal{I}_N} G_i^g (e) \log r_i}. \]
Thus the roles of $r_i$ and $G_i^g (e)$ are exchanged.
By Jensen's inequality for the logarithm, each of these two dimensions is $1$ if and only if $G_i^g (e) = r_i$ for every $i \in \mathcal{I}_N$, which agrees with Remark \ref{rem:Furstenberg} (2), and they are different in general.
For example, let $N = 3$, $r_0 = r_1 = r_2 = 1/3$, $g_0 (z) = z/(z+2)$, $g_1 (z) = (2z+1)/(2z+3)$ and $g_2 (z) = (5z+3)/(3z+5)$.
Then, $g$ is a hyperbolic LF system, $H_i^g (1) = 1$ for every $i \in \mathcal{I}_3$, $(G_0^g (1), G_1^g (1), G_2^g (1)) = (1/2, 1/4, 1/4)$, and
\[ \dim_H \mu_{\varphi} = \frac{3 \log 2}{2 \log 3} = 0.9463\cdots, \ \ \ \dim_H \mu_{\varphi^{-1}} = \frac{3 \log 3}{5 \log 2} = 0.9509\cdots. \]
\end{Rem}

\begin{Rem}[Second iterates and the results in \cite{JKPS09}]\label{rem:JKPS}
(1) An LF system $h$ on $[0,1]$ is hyperbolic if and only if
\[ \max_{i, j \in \mathcal{I}_N} \max_{x \in [0,1]} \left( h_i \circ h_j \right)^{\prime} (x) < 1. \]
We show this.
Let $i \in \mathcal{I}_N$.
By Lemma \ref{lem:basic} (1) and (3) and Lemma \ref{lem:psi} (1), $h_i^{\prime} (x) = (b_{i+1}^h - b_i^h) \, k(c_i^h, x) \le 1$ for $x \in [0,1]$.
If $c_i^h = 0$, then $h_i^{\prime}$ is the constant $b_{i+1}^h - b_i^h$, which is smaller than $1$ by Lemma \ref{lem:basic} (2), and, if $c_i^h \ne 0$, then $h_i^{\prime}$ is strictly monotone on $[0,1]$.
Hence, $h_i^{\prime} (x) < 1$ for $0 < x < 1$.
Let $i, j \in \mathcal{I}_N$ and $x \in [0,1]$, and assume that $(h_i \circ h_j)^{\prime} (x) = 1$.
Since $(h_i \circ h_j)^{\prime} (x) = h_i^{\prime} (h_j (x)) \, h_j^{\prime} (x)$, we have that $h_j^{\prime} (x) = 1$ and $h_i^{\prime} (h_j (x)) = 1$, and hence, by the above, $x \in \{0,1\}$ and $h_j (x) \in \{0,1\}$.
By Lemma \ref{lem:basic} (2), this holds only if $j = 0$ and $x = 0$, or $j = N-1$ and $x = 1$, and then $h_0^{\prime} (0) = 1$ or $h_{N-1}^{\prime} (1) = 1$, respectively.
Hence, if $h$ is hyperbolic, then $(h_i \circ h_j)^{\prime} < 1$ on $[0,1]$ for all $i, j \in \mathcal{I}_N$, and, since these derivatives are continuous on $[0,1]$, the above inequality holds.
If $h$ is parabolic, then $(h_0 \circ h_0)^{\prime} (0) = h_0^{\prime} (0)^2 = 1$ or $(h_{N-1} \circ h_{N-1})^{\prime} (1) = h_{N-1}^{\prime} (1)^2 = 1$, and the above inequality fails.

The maps $h_i \circ h_j$, $i, j \in \mathcal{I}_N$, are the inverse branches of the second iterate of the map on $[0,1)$ whose inverse branches are $h_i$, $i \in \mathcal{I}_N$.
Thus $h$ is hyperbolic if and only if this second iterate is uniformly expanding.
If $h$ is parabolic, then no iterate of this map is uniformly expanding, since $h_{0^n}^{\prime} (0) = 1$ for every $n \ge 1$, or $h_{(N-1)^n}^{\prime} (1) = 1$ for every $n \ge 1$, where $i^n$ denotes the word of length $n$ all of whose entries are $i$.
If $h$ is hyperbolic, then the map itself is not uniformly expanding in general.
For example, if $N = 2$, $h_0 (x) = x/2$ and $h_1 (x) = (3x+1)/(2x+2)$, then $h = \{h_0, h_1\}$ is a hyperbolic LF system on $[0,1]$ for which $(a_0^h, b_0^h, c_0^h) = (1/2, 0, 0)$ and $(a_1^h, b_1^h, c_1^h) = (3/2, 1/2, 1)$, and $h_0^{\prime} (0) = 1/2$, $h_1^{\prime} (1) = 1/4$ and $h_1^{\prime} (0) = 1$.
If the strict inequality holds in the condition (2) of LF systems for every $i \in \mathcal{I}_N$, then $h_i^{\prime} < 1$ on $[0,1]$ for every $i \in \mathcal{I}_N$ (Remark \ref{rem:assumption} (4)), and the map itself is uniformly expanding.\\
(2) Let $f$ and $g$ be hyperbolic LF systems on $[0,1]$, and let $\varphi$ be the solution of \eqref{eq:conj}.
In \cite{JKPS09}, the two maps are expanding, and the logarithms of the derivatives of their inverse branches are negative functions on $[0,1]$, that is, every inverse branch is a uniform contraction.
By (1), this may fail for $T$ and for $\tau$, and it holds for $T^2$ and $\tau^2$, whose inverse branches are $f_i \circ f_j$ and $g_i \circ g_j$, $i, j \in \mathcal{I}_N$, respectively.
These are restrictions of linear fractional transformations which are differentiable on neighborhoods of $[0,1]$, and, by Lemma \ref{lem:conj} (1), $g_i \circ g_j \circ \varphi = \varphi \circ f_i \circ f_j$ for all $i, j \in \mathcal{I}_N$.
Thus $T^2$ and $\tau^2$ are expanding piecewise $C^{1+\epsilon}$ maps with $N^2$ increasing full inverse branches, and $\varphi$ is the conjugacy between them in the sense of \cite{JKPS09}.
Hence, by \cite[Theorem 1.2]{JKPS09}, either $\varphi$ is a $C^{1+\epsilon}$ diffeomorphism of $[0,1]$, or $\varphi$ is a singular function and $\dim_H \mathcal{D}_{\infty} < 1$, where $\mathcal{D}_{\infty}$ is the set of the points of $[0,1]$ at which the derivative of $\varphi$ is $+\infty$.
The dichotomy in Theorem \ref{thm:sing} follows from this and the following two arguments.

First, assume that $\varphi$ is a singular function, that is, $\mu_{\varphi}$ is singular with respect to the Lebesgue measure.
Then, by Remark \ref{rem:inverse} (2), $\mu_{\varphi^{-1}}$ is also singular with respect to the Lebesgue measure, and hence there exists a Borel subset $B$ of $(0,1)$ such that $\ell (B) = 1$ and the derivative of $\varphi^{-1}$ is $0$ at every point of $B$.
Since $\varphi$ is a strictly increasing homeomorphism of $[0,1]$, the derivative of $\varphi$ at $\varphi^{-1} (z)$ is $+\infty$ for every $z \in B$, that is, $\varphi^{-1} (B) \subset \mathcal{D}_{\infty}$.
By Lemma \ref{lem:conj} (3), $\mu_{\varphi} \left( \varphi^{-1} (B) \right) = \ell (B) = 1$.
Hence, $\dim_H \mu_{\varphi} \le \dim_H \mathcal{D}_{\infty} < 1$.

Second, assume that $\varphi$ is differentiable at $0$ and $0 < \varphi^{\prime} (0) < +\infty$, which is the case if $\varphi$ is a $C^{1+\epsilon}$ diffeomorphism of $[0,1]$.
For $h \in \{f, g\}$, let
\[ u_h (x) \coloneqq  \frac{x}{1 + e_0^h x}, \ \ \ x \in [0,1]. \]
Since $e_0^h > -1$, $u_h$ is continuous and strictly increasing on $[0,1]$, $u_h (0) = 0$ and $u_h^{\prime} (0) = 1$.
Since $h_0 (x) = a_0^h x / (c_0^h x + 1)$ and $c_0^h + a_0^h e_0^h = e_0^h$, we have that $u_h (h_0 (x)) = a_0^h \, u_h (x)$ for $x \in [0,1]$.
Let $F \coloneqq  u_g \circ \varphi \circ u_f^{-1}$, which is a function on $[0, u_f (1)]$.
Then, by \eqref{eq:conj} for $i = 0$,
\[ F \! \left( a_0^f \, t \right) = a_0^g \, F(t), \ \ \ t \in [0, u_f (1)]. \]
Furthermore, $F(0) = 0$, $F$ is differentiable at $0$, and $F^{\prime} (0) = \varphi^{\prime} (0)$.
By dividing both sides of the above identity by $t$ and letting $t \to +0$, we have that $a_0^f = a_0^g$.
Hence, for $t \in (0, u_f (1)]$ and $n \ge 1$,
\[ \frac{F(t)}{t} = \frac{F \! \left( (a_0^f)^n \, t \right)}{(a_0^f)^n \, t}, \]
and, since $a_0^f < 1$, the right hand side converges to $F^{\prime} (0)$ as $n \to \infty$.
Thus $F(t) = F^{\prime} (0) \, t$ for $t \in [0, u_f (1)]$.
Since $u_g^{-1} (t) = t / (1 - e_0^g t)$, we have that $\varphi (x) = \lambda x / (1 + s x)$ for $x \in [0,1]$, where $\lambda \coloneqq  \varphi^{\prime} (0)$ and $s \coloneqq  e_0^f - \lambda e_0^g$.
Since $\varphi (1) = 1$, $\lambda = 1 + s$, and hence $s > -1$ and $\varphi = \psi_s$.
Since $\psi_s (f_0 (1)) = \varphi (f_0 (1)) = g_0 (\varphi (1)) = g_0 (1)$, by Lemma \ref{lem:psi} (4) and the definition of $s_0$, $s = s_0$, and the assumption of Theorem \ref{thm:sing} (1) holds.\\
(3) We compare our results with those in \cite{JKPS09}.
The proofs in \cite{JKPS09} are based on the thermodynamic formalism and the multifractal analysis on the shift space.
Our proofs of Theorems \ref{thm:dim} and \ref{thm:sing} do not use them or the results in \cite{JKPS09}, and they are based on the stationary distribution $\nu$.
By \cite[Theorem 1.2]{JKPS09}, the conjugacy is a $C^{1+\epsilon}$ diffeomorphism if and only if the two maps are cohomologically dependent in the sense of \cite{JKPS09}, that is, a non-trivial linear combination of the logarithms of the derivatives of the inverse branches of the two maps, which are regarded as functions on the shift space, is the difference of the composition of a continuous function on the shift space with the shift map and this function.
In general, it is not easy to verify this condition for given maps.
By (2), Theorem \ref{thm:sing} and Proposition \ref{prop:equiv}, for hyperbolic LF systems $f$ and $g$, the maps $T^2$ and $\tau^2$ are cohomologically dependent in the sense of \cite{JKPS09} if and only if the $2N-2$ equations for the coefficients stated after Theorem \ref{thm:sing} hold.
Thus, for hyperbolic LF systems, this condition is decided by finitely many explicit equations.
The Hausdorff dimension considered in \cite{JKPS09} is $\dim_H \mathcal{D}_{\infty}$, which is given in \cite[Theorem 1.1]{JKPS09} in terms of a pressure function, whereas Theorem \ref{thm:dim} gives $\dim_H \mu_{\varphi}$, which is at most $\dim_H \mathcal{D}_{\infty}$ by (2) if $\varphi$ is a singular function.
If $g$ is parabolic, then, by (1), no iterate of $\tau$ is uniformly expanding, and the results in \cite{JKPS09} do not apply to $\varphi$; Theorem \ref{thm:para} is for this case.\\
(4) In the proof of \cite[Theorem 1.2]{JKPS09} in \cite[Section 3]{JKPS09}, it is noted that if the two maps are cohomologically dependent in the sense of \cite{JKPS09}, then the derivatives of their first inverse branches at the point $0$ coincide, and the derivatives of their last inverse branches at the point $1$ coincide.
Let $f$ and $g$ be hyperbolic LF systems on $[0,1]$, and let $\varphi$ be the solution of \eqref{eq:conj}.
The first inverse branches of $T^2$ and $\tau^2$ are $f_0 \circ f_0$ and $g_0 \circ g_0$, respectively, and the last inverse branches of them are $f_{N-1} \circ f_{N-1}$ and $g_{N-1} \circ g_{N-1}$, respectively.
Since $(h_0 \circ h_0)^{\prime} (0) = h_0^{\prime} (0)^2$ and $(h_{N-1} \circ h_{N-1})^{\prime} (1) = h_{N-1}^{\prime} (1)^2$ for $h \in \{f, g\}$, for $T^2$ and $\tau^2$, this necessary condition for the cohomological dependence is that
\[ f_0^{\prime} (0) = g_0^{\prime} (0), \ \ \ f_{N-1}^{\prime} (1) = g_{N-1}^{\prime} (1). \]
By Remark \ref{rem:odds} (3), these two equations are also necessary for the assumption of Theorem \ref{thm:sing} (1), and, if $N = 2$, then they are also sufficient for it.
Hence, by (3), if $N = 2$, then $T^2$ and $\tau^2$ are cohomologically dependent in the sense of \cite{JKPS09} if and only if these two equations hold.
If $N \ge 3$, then these two equations are not sufficient for the cohomological dependence of $T^2$ and $\tau^2$.
Indeed, let $N \ge 3$, let $f_i (x) \coloneqq  (x+i)/N$, $x \in [0,1]$, for $i \in \mathcal{I}_N$, let $g_i \coloneqq  f_i$ for $i \in \mathcal{I}_N$ with $i \ne 1$, and let
\[ g_1 (z) \coloneqq  \frac{3z+1}{N (z+1)}, \ \ \ z \in [0,1]. \]
Then, $f$ and $g$ are hyperbolic LF systems on $[0,1]$, where $f$ is as in Remark \ref{rem:affine} for $r_i = 1/N$, $i \in \mathcal{I}_N$, and $(a_1^g, b_1^g, c_1^g) = (3/N, 1/N, 1)$.
Since $g_0 = f_0$ and $g_{N-1} = f_{N-1}$, the above two equations hold.
Since $g_0 (1) = f_0 (1)$, we have that $s_0 = 0$, and $\psi_{s_0}$ is the identity map.
Since $g_1 \ne f_1$, the assumption of Theorem \ref{thm:sing} (1) does not hold.
Hence, by Theorem \ref{thm:sing} (2), $\varphi$ is a singular function, and, by (3), $T^2$ and $\tau^2$ are not cohomologically dependent in the sense of \cite{JKPS09}.
\end{Rem}

\begin{Rem}[The Gibbs property]\label{rem:Gibbs-property}
Let $f$ and $g$ be LF systems on $[0,1]$, and let $\varphi$ be the solution of \eqref{eq:conj}.
We regard $\mu_{\varphi}$ as a Borel probability measure on $[0,1)$, as in Subsection \ref{subsec:invariant}.
We say that a Borel probability measure $\mu$ on $[0,1)$ has the Gibbs property for a real-valued function $\omega$ on $[0,1)$ if there exist constants $C \ge 1$ and $P \in \mathbb{R}$ such that
\begin{equation}\label{eq:Gibbs}
C^{-1} \le \frac{\mu (J_n (x))}{\exp \left( - n P + \sum_{j=0}^{n-1} \omega (T^j (x)) \right)} \le C
\end{equation}
for every $x \in [0,1)$ and every $n \ge 1$.
This follows the definition of Gibbs measures on shift spaces in Bowen \cite[Theorem 1.2]{B75}, where the intervals $J_n (x)$ correspond to the cylinder sets by the map $x \mapsto (X_n (x))_{n \ge 1}$.\\
(1) If $g$ is hyperbolic, then $\mu_{\varphi}$ has the Gibbs property for the function
\[ \omega_g (x) \coloneqq  \log g_{X_1 (x)}^{\prime} (\varphi (T(x))), \ \ \ x \in [0,1). \]
Indeed, let $x \in [0,1)$, $n \ge 1$ and $\sigma \coloneqq  (X_1 (x), \dots, X_n (x))$.
In the same manner as in the proof of Proposition \ref{prop:limits}, by Lemma \ref{lem:coding} and the chain rule, the sum of $\omega_g (T^j (x))$ over $0 \le j \le n-1$ is $\log g_{\sigma}^{\prime} (\varphi (T^n (x)))$.
Hence, by Lemma \ref{lem:coding} (1), Lemma \ref{lem:conj} (2) and Lemma \ref{lem:id} (2) for $h = g$,
\[ \frac{\mu_{\varphi} (J_n (x))}{\exp \left( \sum_{j=0}^{n-1} \omega_g (T^j (x)) \right)} = \frac{1}{k\left( y_{\sigma}^g, \varphi (T^n (x)) \right)}, \]
and, by Lemma \ref{lem:hyp} (1) and (2) for $h = g$, the right hand side is at least $C_g^{-1}$ and at most $C_g$.
Thus \eqref{eq:Gibbs} holds for $\mu = \mu_{\varphi}$, $\omega = \omega_g$, $C = C_g$ and $P = 0$.\\
(2) If $g$ is parabolic, then there are no real-valued functions $\omega$ on $[0,1)$ such that $\mu_{\varphi}$ has the Gibbs property for $\omega$.
Indeed, assume that $\mu_{\varphi}$ has the Gibbs property for a real-valued function $\omega$ on $[0,1)$, let $C$ and $P$ be the constants in \eqref{eq:Gibbs}, and, for $x \in [0,1)$ and $n \ge 1$, let $\Xi_n (x)$ be the denominator in \eqref{eq:Gibbs}.
Then, $\Xi_{2n} (x) = \Xi_n (x) \, \Xi_n (T^n (x))$.
We use the notation in Section \ref{sec:para}.
Let $i$ be a parabolic letter, and let $\theta \coloneqq  \theta_0$ if $i = 0$, and $\theta \coloneqq  \theta_1$ if $i = N-1$.
Let $n \ge 1$ and $x \in J(i^{2n})$.
By Lemma \ref{lem:coding} (1) and Lemma \ref{lem:cyl} (2) for $h = f$, $X_j (x) = i$ for $1 \le j \le 2n$, and hence $J_{2n} (x) = J(i^{2n})$ and $J_n (x) = J_n (T^n (x)) = J(i^n)$.
Hence, by \eqref{eq:Gibbs},
\[ \mu_{\varphi} (J(i^{2n})) \le C \, \Xi_{2n} (x) = C \, \Xi_n (x) \, \Xi_n (T^n (x)) \le C^3 \, \mu_{\varphi} (J(i^n))^2. \]
By Lemma \ref{lem:run} (1) and (2) for the empty word, $\mu_{\varphi} (J(i^n)) = 1/(1 + n \theta)$ for every $n \ge 1$.
Thus $(1 + n \theta)^2 \le C^3 (1 + 2 n \theta)$ for every $n \ge 1$, which is a contradiction.
In particular, the measure $\mu_{\varphi}$ in Example \ref{exa:Minkowski} does not have the Gibbs property for any real-valued function on $[0,1)$.\\
(3) Assume that $f$ and $g$ are hyperbolic, and let $\widehat{\mu}$ be the measure in Subsection \ref{subsec:invariant}.
By (1) and Proposition \ref{prop:hatmu}, $\widehat{\mu}$ is $T$-invariant and has the Gibbs property for $\omega_g$.
This may open a way to apply known results for Gibbs measures to $\mu_{\varphi}$; see \cite{PW97} for their multifractal analysis, \cite{AJ26, KS09} for the H\"older regularity of their distribution functions, and \cite{JS16} for the Fourier decay of certain Gibbs measures for the Gauss map.
\end{Rem}

\begin{Rem}[Related results]\label{rem:related}
(1) Let $f$ and $g$ be hyperbolic LF systems on $[0,1]$, and let $m$, $\Psi_g$, $\lambda_g$, $\lambda_f$ and $\widehat{\mu}$ be as in Section \ref{sec:dim}.
By Proposition \ref{prop:hatmu} and Remark \ref{rem:entropy}, $\widehat{\mu}$ is $T$-invariant, $T$-ergodic and mutually absolutely continuous with $\mu_{\varphi}$, and $\lambda_g$ and $\lambda_f$ are the entropy and the Lyapunov exponent of $T$ with respect to $\widehat{\mu}$, respectively.
Thus, for $T$ and $\widehat{\mu}$, the formula in Theorem \ref{thm:dim} (2) is of the same form as the formula of Hofbauer and Raith \cite{HR92} stated in Section \ref{sec:intro}.
Our proof of Theorem \ref{thm:dim} does not use the result in \cite{HR92}.
Since the restriction of $\varphi$ to $[0,1)$ is a homeomorphism of $[0,1)$ such that $\varphi \circ T = \tau \circ \varphi$ (Lemma \ref{lem:coding} (2)), and $m$ is the image measure of $\widehat{\mu}$ under $\varphi$, $\lambda_g$ is also the entropy of $\tau$ with respect to $m$.
Furthermore, since $\Psi_g$ is the logarithm of the derivative of $\tau$ on each open interval $(b_i^g, b_{i+1}^g)$, $i \in \mathcal{I}_N$, by \eqref{eq:lambda}, $\lambda_g$ is the Lyapunov exponent of $\tau$ with respect to $m$.
Thus $\dim_H \mu_{\varphi}$ is the ratio of the Lyapunov exponents of $\tau$ and of $T$ along the orbits corresponding to each other by $\varphi$.\\
(2) A dichotomy of the same type as that in Theorem \ref{thm:sing} is known for two lattice embeddings into $PSL_2 (\mathbb{R})$ which are conjugate by a homeomorphism of the upper half-plane: the homeomorphism of the boundary induced by it is either a linear fractional transformation or singular, and, in the latter case, it maps a set of Hausdorff dimension smaller than $1$ onto a set whose complement has Hausdorff dimension smaller than $1$; see Bishop and Steger \cite{BS93} and the references therein.
See also \cite[Section 5]{AJ26} for the H\"older spectra of conjugacies between expanding piecewise $C^{1+\epsilon}$ maps of $[0,1]$, and \cite{KS08} for the sets of non-differentiability points of Minkowski's question-mark function.\\
(3) For the continued fraction map, which has infinitely many decreasing branches and is not in our class, the objects corresponding to those in Subsection \ref{subsec:tau} are classical: $k(y, \cdot)$ is the density in the Brod\'en--Borel--L\'evy formula, the dual maps are the inverse branches of the map itself, and $\nu^g$ corresponds to the Gauss measure; see \cite{IK02}, and see \cite{N81} for natural extensions.\\
(4) From the viewpoint of iterated function systems, $\mu_{\varphi}$ is an eigenmeasure of the dual operator of the Ruelle operator for the maps $f_i$, $i \in \mathcal{I}_N$, with the weights $g_i^{\prime} \circ \varphi$, $i \in \mathcal{I}_N$, as in Fan and Lau \cite{FL99}, that is,
\[ \int_{[0,1)} \sum_{i \in \mathcal{I}_N} g_i^{\prime} (\varphi (x)) \, v(f_i (x)) \, \mu_{\varphi} (dx) = \int_{[0,1)} v \, d\mu_{\varphi} \]
for every bounded Borel measurable function $v$ on $[0,1)$, which follows from Lemma \ref{lem:conj} (3), \eqref{eq:conj} and the change of variables $z^{\prime} = g_i (z)$; see \cite{B15, FL99, JR08} for dimension results, under their respective assumptions, for such measures and for invariant measures of iterated function systems with place-dependent probabilities.\\
(5) Arcozzi \cite{A26} related the equation in \cite{O14}, which is \eqref{eq:conj} for $N = 2$, $f_0 (x) = x/2$ and $f_1 (x) = (x+1)/2$, to the dyadic potential theory on the binary tree.
It is shown there that, for a large class of parameters, the solution is the normalized cumulative capacitary function of a multiplicative dyadic capacity \cite[Theorem 1.2]{A26}, and that the upper bound for $\dim_H \mu_{\varphi}$ in \cite{O14} is sharpened at its endpoint in terms of dyadic Riesz capacities \cite[Theorem 1.3]{A26}; see Remark \ref{rem:Arcozzi} for a problem posed there.
See also Mantica and Totik \cite{MT19} and \cite{O26b} for the regularity, in the sense of logarithmic potential theory, of measures of this kind, which are called conjugacy measures in \cite{O26b}.
\end{Rem}

\begin{Rem}[A problem in \cite{A26}]\label{rem:Arcozzi}
In \cite{A26}, the equation \eqref{eq:conj} is considered for $N = 2$, $f_0 (x) = x/2$, $f_1 (x) = (x+1)/2$ and
\[ g_0 (z) = \frac{z}{(\gamma - \lambda) z + \lambda}, \ \ \ g_1 (z) = \frac{(q-1) z + 1}{(q - \gamma) z + \gamma}, \ \ \ z \in [0,1], \]
where $\lambda > 0$, $\gamma > 1$ and $q > 0$.
For these maps, $g_0^{\prime} (0) = 1/\lambda$, $g_0^{\prime} (1) = \lambda / \gamma^2$, $g_1^{\prime} (0) = q (\gamma - 1) / \gamma^2$ and $g_1^{\prime} (1) = (\gamma - 1)/q$.
In this remark, we consider \eqref{eq:conj} for these maps.\\
(1) Assume that $g = \{g_0, g_1\}$ is an LF system on $[0,1]$.
In \cite[Subsection 3.3.1]{A26}, it is posed as a problem to determine the infimum of the numbers $s \ge 0$ for which there exists a Borel subset $K$ of $[0,1]$ such that $\mu_{\varphi} (K) = 1$ and the dyadic Riesz $s$-capacity of the preimage of $K$ under the binary coding map is $0$.
By the relations between these capacities and Hausdorff dimensions stated in \cite[Subsection 2.3]{A26}, this infimum is equal to $\dim_H \mu_{\varphi}$.
Since $\ell (J_n (x)) = 2^{-n}$ for every $n \ge 1$ and every $x \in [0,1)$, if $g$ is hyperbolic, then, by Theorem \ref{thm:dim} (2), $\dim_H \mu_{\varphi}$ is the limit in Remark \ref{rem:two} (1) divided by $\log 2$, and, if $g$ is parabolic, then $\dim_H \mu_{\varphi} = 0$ by Theorem \ref{thm:para} (2).
Thus our results give this infimum.
By \cite[Theorem 3.3]{A26}, its value is known in the case that corresponds to $g_0^{\prime} (0) + g_1^{\prime} (1) = 1$, and it agrees with the above value.\\
(2) Let $0 < \sigma < 1$.
By \cite[Theorems 1.1 and 1.2]{A26}, the normalized cumulative capacitary function of the dyadic Riesz $\sigma$-capacity in \cite{A26} is the solution $\varphi$ of \eqref{eq:conj} for $(\lambda, \gamma, q) = (2^{\sigma}, 1 + 2^{\sigma - 1}, 2)$, that is, for
\[ g_0 (z) = \frac{2z}{(2 - 2^{\sigma}) z + 2^{\sigma + 1}}, \ \ \ g_1 (z) = \frac{2z + 2}{(2 - 2^{\sigma}) z + 2 + 2^{\sigma}}, \ \ \ z \in [0,1]. \]
Then, $g = \{g_0, g_1\}$ is a hyperbolic LF system on $[0,1]$ for which
\[ (a_0^g, b_0^g, c_0^g) = \left( 2^{-\sigma}, 0, 2^{-\sigma} - \frac{1}{2} \right), \ \ \ (a_1^g, b_1^g, c_1^g) = \left( \frac{2}{2 + 2^{\sigma}}, \frac{2}{2 + 2^{\sigma}}, \frac{2 - 2^{\sigma}}{2 + 2^{\sigma}} \right), \]
and, by Lemma \ref{lem:basic} (4), $g_0^{\prime} (0) = 2^{-\sigma}$ and $g_1^{\prime} (1) = 2^{\sigma - 2}$, whose sum is larger than $1$.
Hence, by (1) and Remark \ref{rem:two} (1),
\[ \dim_H \mu_{\varphi} = \frac{2 \left( L(2^{-\sigma}) + L(2^{\sigma - 2}) - \pi^2 / 6 \right)}{\left( \sigma \log 2 - \log (2^{\sigma} - 1) - \log (4 - 2^{\sigma}) \right) \log 2}, \]
which is approximately $0.9063$ if $\sigma = 1/2$.
This solves the problem in \cite[Subsection 3.3.1 (i)]{A26}.\\
(3) By \cite[Theorem 1.1 (i)]{A26}, in the trapping range of \cite{A26}, that is, if $\lambda \ge 1$ and $q \ge \gamma - 1$, \eqref{eq:conj} has a unique continuous and strictly increasing solution $\varphi$ such that $\varphi (0) = 0$ and $\varphi (1) = 1$.
The problem in (1) is posed in this range, which is the set of the parameters such that $g_0^{\prime} (0) \le 1$ and $g_1^{\prime} (1) \le 1$.
In this range, $\{g_0, g_1\}$ is an LF system on $[0,1]$ if and only if $g_0^{\prime} (1) \le 1$ and $g_1^{\prime} (0) \le 1$, that is, $\lambda \le \gamma^2$ and $q (\gamma - 1) \le \gamma^2$.
\end{Rem}

\vspace{1pc}

\noindent{\bf Use of generative AI} \ During the preparation of this work, the author used Anthropic's Claude Fable 5.2 and OpenAI's ChatGPT 6.2 Astra to assist with exploration, structuring, drafting, and refining the text. 
The author takes full responsibility for the content and integrity of the publication. 

\bibliographystyle{amsplain}
\bibliography{dimension_formula_lft_conj}

\end{document}